\documentclass[10pt,leqno]{article}
\usepackage{amssymb,amsbsy,amsmath,amsfonts,amssymb,amscd, mathrsfs}

\usepackage[english]{babel}
\usepackage[T1]{fontenc}
\usepackage{indentfirst}
\usepackage[all]{xy}

\makeatletter
\@addtoreset{equation}{section}
\makeatother

\newtheorem{statement}{}[section]
\newtheorem{theoreme}[statement]{Theorem}
\newtheorem{lemme}[statement]{Lemma}

\newtheorem{proposition}[statement]{Proposition}
\newtheorem{definition}[statement]{Definition}
\newtheorem{corollaire}[statement]{Corollary}

\newcommand\C{\mathbb C}
\newcommand\N{\mathbb N}
\newcommand\R{\mathbb R}
\newcommand\T{\mathbb T}
\newcommand\D{\mathbb D}
\newcommand\Z{\mathbb Z}
\newcommand\e{{\rm e}}
\renewcommand\P{\mathbb P}
\newcommand\eps{\varepsilon}
\newcommand\ind{{\rm 1\kern-.30em I}}

\title{Composition operators on Hardy-Orlicz spaces}
\author{Pascal Lef\`evre, Daniel Li,\\ Herv\'e Queff\'elec, Luis Rodr{\'\i}guez-Piazza}

\date{\footnotesize \today}

\begin{document}

\maketitle

\noindent{\bf Abstract.} \emph{We investigate composition operators on Hardy-Orlicz spaces when the 
Orlicz function $\Psi$ grows rapidly: compactness, weak compactness, to be $p$-summing, order 
bounded, \dots, and show how these notions behave according to the growth of $\Psi$. We introduce an 
adapted version of Carleson measure. We construct various examples showing that our results are essentially 
sharp. In the last part, we study the case of Bergman-Orlicz spaces.}
\medskip

\noindent{\bf Mathematics Subject Classification.} Primary: 47 B 33 -- 46 E 30; Secondary: \par
\smallskip
\noindent{\bf Key-words.} \emph{Bergman-Orlicz space -- Carleson measure -- composition operator -- 
Hardy-Orlicz space}

\section{Introduction.}
Composition operators on the classical Hardy spaces $H^p$ have been widely studied 
(see \cite{Shap}, and \cite{Co-McC}, and references therein; see also \cite{Jar1} and \cite{Jar2}, and 
\cite{Choa}, \cite{Cima}, \cite{Jaoua2}, \cite{Mata}, \cite{Sarason}, \cite{Shap2}, \cite{Shap-SS} for some 
more recent works), but it seems that one has not paid much  attention to the Hardy-Orlicz spaces 
(in \cite{Stromb1} and \cite{Stromb2}, J.-O. Str\"omberg studied Hardy-Orlicz spaces in the case when the 
Orlicz function $\Psi$ increases smoothly;  see also \cite{Liu} for composition operators). We shall 
investigate what happens when the Orlicz function grows more rapidly than a power function.\par
Recall that, given an analytic self-map $\phi \colon \D \to \D$ of the unit disk $\D$, the 
\emph{composition operator} associated to $\phi$ is the map $C_\phi \colon f \mapsto f\circ \phi$. This map 
may operate on various Banach spaces $X$ of analytic functions on $\D$ (Hardy spaces, Bergman spaces, 
\dots, and their weighted versions (see \cite{Zo} for instance), Bloch spaces ${\mathscr B}$ and 
${\mathscr B}_0$, $BMOA$ and $VMOA$, Dirichlet spaces (see \cite{AFP}), or some more general spaces as 
Nevanlinna or Smirnov classes: \cite{Choa2}, \cite{ChoaShap}, \cite{Jaoua1}, \cite{JaXi}; see also 
\cite{Ba1}, \cite{Ba2}, \cite{GoHe} for composition operators on ${\mathscr H}^p$ spaces of Dirichlet series, 
though they are not induced by an analytic self-map of $\D$). The goal is to link properties of the composition 
operator $C_\phi\colon X \to X$ (compactness, strong or weak, for example) to  properties of the 
\emph{symbol} $\phi$ (essentially its behaviour near the frontier of $\D$). 
For that study, one can roughly speaking (see \cite{Co-McC}, Chapter 4, though their notions are different 
from ours) distinguish two kind of spaces.\par
1) The \emph{small spaces} $X$; those spaces are in a sense close to the Hardy space $H^\infty$: the 
compactness of $C_\phi\colon X \to X$ is very restrictive and it imposes severe restrictions on $\phi$. 
For  example, if $X=H^\infty$, a theorem of J. Schwartz (\cite{Schw}) implies that 
$\C_\phi \colon H^\infty \to H^\infty$ is compact if and only if $\|\phi\|_\infty <1$ (weakly compact 
suffices: see \cite{Pascal}), which in turn implies reinforced compactness properties for $C_\phi$. For 
example, $C_\phi \colon H^\infty \to H^\infty$ is nuclear and $1$-summing as soon as it is compact.\par
2) The \emph{large spaces} $X$; those spaces are in a sense close to the Hardy space $H^1$: the compactness 
of $C_\phi \colon X \to X$ can take place fairly often, and in general implies no self-improvement. For 
example, for $X= H^2$, $C_\phi \colon H^2 \to H^2$ can be compact without being 
Hilbert-Schmidt, even if $\phi$ is injective (\cite{Shap-T}, Theorem 6.2). Another formulation (which lends 
better to generalizations in the non-Hilbertian case) is that $C_\phi$ can be non-order bounded (see 
\cite{HJ}, and our Section 3) and yet compact.
\medskip

In this paper, we shall rather be on the small space side, since we shall work in spaces associated to a very 
large Orlicz function $\Psi$ (typically: $\Psi(x)= \Psi_2(x) = \e^{x^2} -1$), and the previous situation will 
not take place: our operators will be {\it e.g.} order-bounded as soon as they are (weakly) compact, even 
if the situation is not so extreme as for $H^\infty$. However, for slightly smaller Orlicz functions (for instance 
$\Psi(x)= \exp\big[\big(\log (x+1)\big)^2\big] - 1$), the situation is closer to the $H^2$ case: the 
composition operators may be compact on $H^\Psi$, but not order-bounded (Theorem~\ref{suite exemple}).
\medskip

This paper is divided into five parts. Section 1 is this Introduction. In Section 2, which is essentially notational, 
we recall some more or less standard facts on Orlicz functions $\Psi$, on associated Orlicz spaces $L^\Psi$, 
and the ``little'' Orlicz space $M^\Psi$, and their banachic properties, associated with various slow growth 
conditions (indicated by subscripts: $\Delta_1$, $\Delta_2$, \dots) or fast growth conditions (indicated 
by superscripts: $\Delta^0$, $\Delta^1$, $\Delta^2$, \dots).
\smallskip

In Section 3, we introduce the Hardy-Orlicz space $H^\Psi$, and its (or his) ``little brother'' $HM^\Psi$. These 
spaces have already been studied (see \cite{Ha}, \cite{Ros-Rov}), but 
rather for slowly growing functions $\Psi$ (having $\Delta_2$ most of the time), and their definition is 
not so clearly outlined, so we give a detailed exposition of the equivalence of the two natural definitions that 
one has in mind (if one wants to extend the case of Hardy spaces $H^p$ associated to $\Psi(x)=x^p$), as 
well as of the automatic boundedness (through the Littlewood subordination principle and the case of inner 
self-maps of the disk) of composition operators on those spaces. Two of the main theorems are Theorem 
\ref{corollaire Delta2} and Theorem \ref{exemple}. Roughly speaking, Theorem \ref{corollaire Delta2} 
says the following: if $\Psi$ is very fast growing (having $\Delta^2$ more precisely), $H^\Psi$ is a small 
space, the (weak) compactness of $C_\phi$ is very restrictive, and even if the situation is not so 
extreme as for $H^\infty$ ($\|\phi\|_\infty <1$), $\phi$ has to tend to the boundary very slowly, and 
$C_\phi$ is automatically order-bounded into $M^\Psi$. However, Theorem \ref{exemple} shows the limits 
of this self-improvement: $C_\phi$ may be order-bounded into $M^\Psi$ (and hence compact), but 
$p$-summing for no finite $p$. We also show that, when $\Psi$ has $\Delta^2$ growth,  there are always 
symbols $\phi$ inducing compact composition operators on $H^2$ (even Hilbert-Schmidt), but not 
compact on $H^\Psi$.
\smallskip

Section 4 is devoted to the use of Carleson measures. The usefulness of those 
measures in the study of composition operators is well-known (see \cite{Du}, \cite{Co-McC}, \cite{BlasJa}) 
and, to our knowledge, first explicitly used for compactness in \cite{McC}. In particular, we recall the following 
necessary and sufficient condition for $C_\phi \colon H^2 \to H^2$ to be compact: if $h>0$ and 
$w\in \partial \D$, consider the \emph{Carleson window}: 
\begin{displaymath}
W(w,h)= \{z\in \D\,;\ |z | \geq 1 - h\hskip 2mm \text{and}\hskip 2mm |\arg (z \bar{w})| \leq h\}.  
\end{displaymath}
If $\phi$ is an analytic self-map of $\D$ with boundary values $\phi^\ast$, and 
$\mu_\phi = \phi^\ast(m)$ denotes the image under $\phi^\ast$ of the normalized Lebesgue measure 
(Haar measure) on $\T=\partial \D$, the measure $\mu_\phi$ is always a Carleson measure, {\it i.e.}:
\begin{displaymath}
\sup_{w\in \partial \D} \mu_\phi\big( W(w,h)\big) = O\,(h).
\end{displaymath}
Now, we can state:

\begin{theoreme} [B. MacCluer \cite{McC}] \label{Barbara}
The composition operator $C_\phi\colon H^2 \to H^2$ is compact if and only if $\mu_\phi$ 
satisfies the ``\emph{little-oh}'' condition, {\it i.e.} if and only if:
\begin{equation}\label{Carleson petito}
\sup_{w\in \partial \D} \mu_\phi\big( W(w,h)\big) = o\,(h)\hskip 2mm \text{as}\hskip 2mm h\to 0.
\end{equation}
\end{theoreme}

There is another famous necessary and sufficient compactness condition, due to J. Shapiro 
(\cite{Shap: Annals}): Let us denote by $N_\phi$ the Nevanlinna counting function of $\phi$, {\it i.e.}:
\begin{displaymath}
N_\phi(w)= \left\{
\begin{array}{cl}
\sum\limits_{\phi(z)=w} \log \frac{1}{|z|} & \text{if}\ w\not= \phi(0)\ \text{and}\ w\in \phi(\D) \\
0 & \text{if}\  w\notin \phi(\D).
\end{array}
\right.
\end{displaymath}
By Littlewood's inequality, one always has (see \cite{Co-McC}, page 33):
\begin{displaymath}
N_\phi(w)= O\,(1- |w|).
\end{displaymath}
Now, Shapiro's Theorem reads:

\begin{theoreme}[J. Shapiro \cite{Shap: Annals}]\label{Shapiro}
The composition operator $C_\phi\colon H^2 \to H^2$ is compact if and only if  $N_\phi$ satisfies 
the ``\emph{little-oh}'' condition, {\it i.e.}:
\begin{equation}\label{Nevanlinna petito}
N_\phi(w)= o\,(1- |w|) \hskip 2mm \text{as} \hskip 2mm |w|\stackrel{<}{\to} 1.
\end{equation}
\end{theoreme}

Theorem \ref{Shapiro} is very elegant, and probably more ``popular'' than Theorem \ref{Barbara}. Yet, it is 
difficult to apply because the assumption \eqref{Nevanlinna petito} is difficult to check. Here, we shall appeal 
to Theorem \ref{Barbara} to prove (Theorem \ref{meme module}) that the compactness of 
$C_\phi \colon H^2 \to H^2$ cannot be read on $|\phi^\ast|$ when $\phi$ is not finitely valent; more 
precisely, there are two analytic self-maps $\phi_1$ and $\phi_2\colon \D \to \D$ such that:
$|\phi_1^\ast| = |\phi_2^\ast|$ $m$-{\it a.e.}, but $C_{\phi_1}\colon H^2 \to H^2$ is not compact, though 
$C_{\phi_2}\colon H^2 \to H^2$ is compact.\par

We show then that every composition operator which is compact on $H^\Psi$ is necessarily compact on $H^p$ 
for all $p <\infty$. However, there exist (see above, or Section 3) symbols $\phi$ inducing compact composition 
operators on $H^2$ but which are not compact on $H^\Psi$, when $\Psi$ has $\Delta^2$ growth. Hence 
condition \eqref{Carleson petito} does not suffice to characterize the compact composition operators on 
$H^\Psi$. We have to replace Carleson measures and condition \eqref{Carleson petito} by what we may 
call ``$\Psi$-Carleson measures'', and an adaptated ``\emph{little-oh}'' condition, which allows us to 
characterize compactness for composition operators. It follows that if $\Psi \in \Delta^0$, then the weak 
compactness of $C_\phi \colon H^\Psi \to H^\Psi$ implies its compactness.
\par
We also show that the above example $\phi_2$ induces, for an Orlicz function $\Psi$  which does not satisfy 
$\Delta^2$, but which satisfies $\Delta^1$, a composition operator on 
$H^\Psi$ which is compact, but not order bounded into $M^\Psi(\T)$ (Theorem \ref{suite exemple}), 
showing that the assumption that $\Psi\in \Delta^2$ in Theorem \ref{corollaire Delta2} is not 
only a technical assumption.
\smallskip

In Section 5, we introduce the Bergman-Orlicz spaces. Let us remind that, in the Hilbertian case, the study of 
compactness of composition operators is simpler for the Bergman space ${\mathscr B}^2$ than for the 
Hardy space $H^2$. For example, we have the following:

\begin{theoreme} {\rm (see \cite{Shap})} \hskip 5mm \par
i) $C_\phi \colon {\mathscr B}^2 \to {\mathscr B}^2$ is compact if and only if 
\begin{equation}\label{derivee angulaire}
\lim_{|z|\stackrel{<}{\to} 1} \frac{ 1 - |\phi(z)|}{1 - |z|} =+\infty.
\end{equation}\par
ii) \eqref{derivee angulaire} is always necessary for $C_\phi \colon H^2 \to H^2$ to be compact, and 
it is sufficient when $\phi$ is injective, or only finitely valent.\par
iii) There are Blaschke products $\phi$ satisfying \eqref{derivee angulaire} for which 
$C_\phi \colon H^2 \to H^2$ is (in an obvious manner) non-compact.
\end{theoreme}

We perform here a similar study for the Bergman-Orlicz space ${\mathscr B}^\Psi$, and compare the situation 
with that of the Hardy-Orlicz space $H^\Psi$. We are naturally led to a reinforcement of 
\eqref{derivee angulaire} under the form:
\begin{equation}\label{renforcement}
\frac{\Psi^{-1}\bigg[{\displaystyle \frac{1}{(1 - |\phi(a)|)^2}}\bigg]}
{\Psi^{-1}\bigg[{\displaystyle \frac{1}{(1 - |a|)^2}}\bigg]} 
\mathop{\longrightarrow}_{|a| \mathop{\to}\limits^< 1} 0
\end{equation}
(always necessary, and sufficient when $\Psi$ is $\Delta^2$), which reads, in the case 
$\Psi(x)=\Psi_2(x)=\e^{x^2}-1$:
\begin{equation}\label{renforcement particulier}
1 -|\phi(z)| \geq  c_\eps (1 -|z|)^\eps \hskip 5mm \text{for all}\hskip 2mm \eps>0.
\end{equation}
In \cite{Shap}, the construction of a Blaschke product satisfying \eqref{derivee angulaire} is fairly 
delicate, and appeals to Frostman's Lemma and Julia-Caratheodory's Theorem on non-angular 
derivatives at the boundary. Here, we can no longer use these tools for the reinforcement 
\eqref{renforcement particulier}, so we do make a direct construction, using the Parseval 
formula for finite groups. In passing, the construction gives a simpler proof for 
\eqref{derivee angulaire}. Otherwise, the theorem which we obain is similar to Shapiro's one, 
if one ignores some technical difficulties due to the non-separability of ${\mathscr B}^\Psi$: 
we have to ``transit'' by the smaller Bergman-Morse-Transue space ${\mathscr B}M^\Psi$, 
which is the closure of $H^\infty$ in ${\mathscr B}^\Psi$, is separable, and has 
${\mathscr B}^\Psi$ as its bidual.

\section{Notation}

Let $\D$ be the open unit disk of the complex plane, that is the set of complex numbers 
with modulus strictly less than $1$, and $\T$ the unit circle, \emph{i.e.} the set of 
complex numbers with modulus $1$.
\par\medskip

We shall consider in this paper Orlicz spaces defined on a probability space 
$(\Omega, \P)$, which will be the unit circle $\T$, with its (normalized) Haar measure $m$ 
(most often identified with the normalized Lebesgue measure $dx/2\pi$ on the interval 
$[0, 2\pi]$), or the open unit disk $\D$, provided with the normalized area measure 
${\mathscr A}$.
\par\medskip

By an Orlicz function, we shall understand that $\Psi\colon [0, \infty] \to [0,\infty]$
is a non-decreasing convex function such that $\Psi(0)=0$ and $\Psi(\infty)=\infty$. To avoid
pathologies, we shall assume that we work with an Orlicz function $\Psi$ having the
following additional properties: $\Psi$ is continuous at $0$, strictly convex (hence increasing), and such that
\begin{displaymath}
\frac{\Psi(x)}{x}\mathop{\longrightarrow}_{x\to \infty} \infty.
\end{displaymath}
This is essentially to exclude the case of $\Psi(x) = ax$.\par
If $\Psi'$ is the left (or instead, if one prefers, the right) derivative of $\Psi$, one has
$\Psi(x)=\int_0^x \Psi'(t)\,dt$ for every $x>0$.
\par\smallskip
The Orlicz
space $L^\Psi(\Omega)$ is the space of all (equivalence classes of) measurable
functions $f\colon \Omega\to\C$ for which there is a constant $C > 0$ such
that
\begin{displaymath}
\int_\Omega \Psi\Big(\frac{\vert f(t)\vert} {C}\Big)\,d\P(t) < +\infty
\end{displaymath}
and then $\Vert f\Vert_\Psi$ (the \emph{Luxemburg norm}) is the infimum of all possible
constants $C$ such that this integral is $\leq 1$. The Morse-Transue space $M^\Psi(\Omega)$ 
is the subspace generated by $L^\infty(\Omega)$, or, equivalently, the subspace of all  
functions $f$ for which the above integral is finite for all $C>0$.\par
To every Orlicz function is associated the complementary Orlicz function
$\Phi=\Psi^\ast\colon [0,\infty] \to [0,\infty]$ defined by:
\begin{displaymath}
\Phi(x)=\sup_{y\geq 0} \big(xy - \Psi(y)\big),
\end{displaymath}
\smallskip

The extra assumptions on $\Psi$ ensure that $\Phi$ is itself strictly convex.\par

When $\Phi$ satisfies the $\Delta_2$ condition (see the definition below), $L^\Psi$ is (isomorphically, if $L^\Phi$ is 
itself normed by the Luxemburg norm) the dual space of $L^\Phi$, which is, in turn, the dual of $M^\Psi$.
\bigskip

\subsection{Growth conditions}

We shall have to use various growth conditions for the Orlicz function $\Psi$. These conditions 
are usually denoted as $\Delta$-conditions. Our interest is in Orlicz functions which have a 
somewhat fast growth. Usually, some of these conditions are defined through a moderate 
growth condition on the complementary function $\Phi$ of $\Psi$, and the condition $\Delta$ 
for the Orlicz function is translated as a $\nabla$-condition for the complementary function.
So we shall distinguish between \emph{moderate growth} conditions, that we shall define for 
the complementary Orlicz function, and \emph{fast growth} conditions. To emphasize this 
distinction, we shall denote, sometimes in changing the usual notation 
(see \cite{Krasno, Rao}), the moderate growth conditions with a subscript, and the
fast growth conditions with a superscript.\par
\medskip\goodbreak

\noindent{\bf Moderate growth conditions}\par\smallskip
\begin{itemize}
\item The Orlicz function $\Phi$ satifies the \emph{$\Delta_1$-condition} ($\Phi\in \Delta_1$) 
if, for some constant $c>0$, one has:
\begin{displaymath}
\Phi(xy)\leq c\,\Phi(x)\Phi(y)
\end{displaymath}
for $x,y$ large enough.\par
This is equivalent to say that
\begin{displaymath}
\Phi(axy) \leq \Phi(x)\Phi(y)
\end{displaymath}
for some constant $a>0$ and $x,y$ large enough.\par
This condition is usually denoted by $\Delta'$ (see \cite{Rao}, page 28).
\item $\Phi$ satisfies the \emph{$\Delta_2$-condition} ($\Phi\in \Delta_2$) if
\begin{displaymath}
\Phi(2x)\leq K\,\Phi(x)
\end{displaymath}
for some constant $K>1$ and $x$ large enough.
\end{itemize}

One has:
\begin{displaymath}
\Phi\in \Delta_1 \hskip 5mm \Rightarrow \hskip 5mm \Phi\in \Delta_2.
\end{displaymath}
\goodbreak

\noindent{\bf Fast growth conditions}\par\smallskip

\begin{itemize}
\item The Orlicz function $\Psi$ satisfies the \emph{$\Delta^0$-condition} 
($\Psi\in \Delta^0$) if (see \cite{LLQR}), for some $\beta >1$:
\begin{displaymath}
\lim_{x\to +\infty} \frac{\Psi(\beta x)}{\Psi(x)}= +\infty.
\end{displaymath}
A typical example is 
$\Psi (x) = \exp\big[\log (x+2) \log \log (x+2)\big] - 2^{\log\log 2}$; another is 
$\Psi (x) = \exp \big[\big(\log (x +1)^{3/2}\big)\big] -1$.
\item The Orlicz function $\Psi$ satisfies the \emph{$\Delta^1$-condition} 
($\Psi\in \Delta^1$) if there is
some $\beta>1$ such that:
\begin{displaymath}
x\Psi(x)\leq \Psi(\beta x)
\end{displaymath}
for $x$ large enough.
\end{itemize}
Note that this latter condition is usually written as $\Delta_3$-condition, with a subscript 
(see \cite{Rao}, \S 2.5).
This notation fits better with our convention, and the superscript $1$ agrees with the fact 
that this $\Delta^1$-condition is between the $\Delta^0$-condition and the following 
$\Delta^2$-condition. $\Psi\in \Delta^1$ implies that
\begin{displaymath}
\Psi(x) \geq \exp\big(\alpha\, (\log x)^2\big)
\end{displaymath}
for some $\alpha>0$ and $x$ large enough (see \cite{Rao}, Proposition 2, page 37). A typical 
example is $\Psi(x)= \e^{(\log( x +1))^2} -1$.
\par
\begin{itemize}
\item The Orlicz function $\Psi \colon [0,\infty) \to [0,\infty)$ is said to satisfy 
the \emph{$\Delta^2$-condition}
($\Psi\in \Delta^2$) if there exists some $\alpha>1$ such that:
\begin{displaymath}
\Psi(x)^2 \leq \Psi(\alpha x)
\end{displaymath}
for $x$ large enough.
\end{itemize}
This implies that
\begin{displaymath}
\Psi(x)\geq \exp(x^\alpha)
\end{displaymath}
for some $\alpha>0$ and $x$ large enough (\cite{Rao}, Proposition 6, page 40). A typical example 
is $\Psi(x)= \Psi_2(x)=\e^{x^2}-1$.
\medskip\goodbreak

\noindent{\bf Conditions of  regularity}\par\smallskip

\begin{itemize}
\item The Orlicz function $\Psi$ satisfies the \emph{$\nabla_2$-condition} 
($\Psi\in \nabla_2$) if its complementary
function $\Phi$ satisfies the $\Delta_2$-condition.\\
This is equivalent to say that for some constant $\beta>1$ and some $x_0>0$, one has
$\Psi(\beta x) \geq 2\beta\,\Psi(x)$ for $x\geq x_0$, and that implies that
$\frac{\Psi(x)}{x}\mathop{\longrightarrow}\limits_{x\to \infty} \infty$. In particular,
this excludes the case $L^\Psi = L^1$.\par

\item The Orlicz function $\Psi$ satisfies the \emph{$\nabla_1$-condition} 
($\Psi\in \nabla_1$) if its complementary
function $\Phi$ satisfies the $\Delta_1$-condition.\\
This is equivalent to say that
\begin{displaymath}
\Psi (x)\Psi (y)\leq \Psi (bxy)
\end{displaymath}
for some constant $b>0$ and $x,y$ large enough.\\ 
All power functions $\Psi (x) = x^p$ satisfy $\nabla_1$, but $\Psi (x) = x^p\log (x +1)$ does not.
\end{itemize}

One has (see \cite{Rao}, page 43):
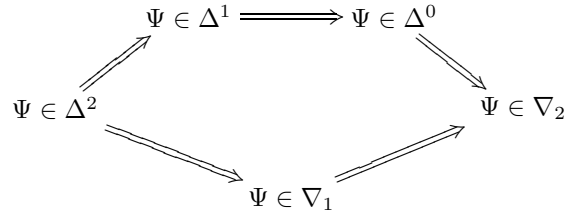
\begin{displaymath}
\xymatrix@R=20pt@C=0pt{
 & \quad \Psi\in \Delta^1  \ar@2{->}[rr]  & &  \Psi\in \Delta^0  \ar@2{->}[dr]  &  \\
\Psi\in \Delta^2\,  \ar@2{ ->}[ur]  \ar@2{->}[drr] &   & &  &   \quad\Psi\in \nabla_2  \\
 &  & \Psi\in \nabla_1 \ar@2{->}[urr] &  &  &  \\
 }
\end{displaymath}

But $\Delta^1$ does not imply $\nabla_1$. That $\nabla_1$ does not even imply $\Delta^0$ is clear 
since any power function $\Psi(x)=x^p$ ($p\geq 1$) is in $\nabla_1$.
\medskip

\subsection{Some specific functions}
In this paper, we shall make a repeated use of the following functions:\par
\begin{itemize}
\item If $\Psi$ is an Orlicz function, we set, for every $K>0$:
\begin{equation}
\chi_K(x)=\Psi\big(K \Psi^{-1}(x)\big), \hskip 3mm x>0.
\end{equation}
For example, if $\Psi (x) =\e^x -1$, then $\Psi^{-1} (x) =\log (1+x)$, and $\chi_K (x) = (1+x)^K -1$.
\par

Note that:
\begin{itemize}
\item $\Psi\in \Delta^0$ means that 
$\displaystyle \frac{\chi_\beta (u)}{u} \mathop{\longrightarrow}_{u \to \infty} +\infty$, 
for some $\beta >1$.
\item $\Psi \in \Delta^1$ means that $\chi_\beta (u) \geq u \Psi^{-1}(u) $ for $u$ large enough, for some 
$\beta >1$.
\item $\Psi \in \Delta^2$ means that $\chi_\alpha (u) \geq u^2 $ for $u$ large enough, for some $\alpha >1$.
\item $\Psi \in \nabla_1$ means that $\chi_A(u) \geq (\Psi(A)/b)\, u$ for $u$ large enough and for 
every $A$ large enough, for some $b>0$.
\end{itemize}
\item For $|a|=1$ and $0\leq r<1$, $u_{a,r}$ is the function defined on the unit disk $\D$ by:
\begin{equation}
u_{a,r}(z)=\Big(\frac{1-r}{1 - \bar{a} r z}\Big)^2, \hskip 3mm |z|<1.
\end{equation}
Note that $\| u_{a,r}\|_\infty =1$ and $\| u_{a,r}\|_1 \leq 1-r$.
\end{itemize}
\bigskip\goodbreak

\section{Composition operators on Hardy-Orlicz spaces}

\subsection{Hardy-Orlicz spaces}

It is well-known that the classical $H^p$ spaces ($1\leq p \leq \infty$) can be defined in two 
equivalent ways:\par
1) $H^p$ is the space of analytic functions $f\colon \D \to \C$ for which, setting 
$f_r(t) = f(r \e^{it})$:
\begin{displaymath}
\| f\|_{H^p} = \sup_{0 \leq r <1} \|f_r\|_p
\end{displaymath} 
is finite (recall that the numbers $\|f_r\|_p$ increase with $r$). When $f\in H^p$, the 
Fatou-Riesz Theorem asserts that the boundary limits 
$f^\ast(t)= \lim_{r\mathop{\to}\limits^{<}1} f_r(t)$ exist almost everywhere and 
$\| f\|_{H^p}=\|f^\ast\|_p$. One has 
$f^\ast \in L^p([0, 2\pi])$, and its Fourier coefficients $\widehat f^\ast (n)$ vanish for $n<0$.\par
2) Conversely, for every function $g\in L^p([0, 2\pi])$ whose Fourier coefficients $\hat g (n)$ 
vanish for $n<0$, the analytic extension $P[g]\colon \D \to \C$ defined by 
$P[g](z)= \sum_{n\geq 0} \hat g(n) z^n$ is in $H^p$ and $g$ is the boundary limit $(P[g])^\ast$ 
of $P[g]$.\par
\smallskip

Hardy-Orlicz spaces $H^\Psi$ are defined in a similar way. However, we did not find very satisfactory 
references, and, though the reasonings are essentially the same as in the classical case, the lack of homogeneity 
of $\Psi$ and the presence of the two spaces $M^\Psi$ and $L^\Psi$ gives proofs which are not so obvious and 
therefore we shall give some details.\par 
It should be noted that our definition is not exactly the same as the one given in \cite{Rao}, \S~9.1.
\smallskip

We shall begin with the following proposition. 

\begin{proposition}\label{Hardy-Orlicz}
Let $f\colon \D \to \C$ be an analytic function. For every Orlicz function $\Psi$, the 
following assertions are equivalent:\par
$1)$ $\sup_{0\leq r <1} \| f_r\|_\Psi < +\infty$, where $f_r(t)=f(r\e^{it})$;\par
$2)$ there exists $f^\ast \in L^\Psi([0, 2\pi])$ such that $\widehat f^\ast(n)=0$ for $n<0$ and 
for which $f(z)= \sum_{n\geq 0} \widehat f^\ast(n) z^n$, $z\in \D$.\par 
When these conditions are satisfied, one has 
$\|f^\ast\|_\Psi=\sup_{0\leq r <1} \| f_r\|_\Psi$.
\end{proposition}

Let us note that, since $\Psi$ is convex and increasing, $\Psi(a|f|)$ is subharmonic on $\D$, 
and hence the numbers $\int_\T \Psi(a|f_r|)\,dm$ increase with $r$, for every $a>0$.
\par\medskip

This proposition leads to the following definition.

\begin{definition}  
Given an Orlicz function $\Psi$, the \emph{Hardy-Orlicz} space $H^\Psi$ associated to $\Psi$ 
is the space of analytic functions $f\colon \D \to \C$ such that one of the equivalent 
conditions of the above proposition is satisfied. The norm of $f$ is defined by 
$\|f\|_{H^\Psi}= \|f^\ast\|_\Psi$. We shall denote by $HM^\Psi$ the 
\emph{Hardy-Morse-Transue space}, {\it i.e.} the subspace 
$\{f\in H^\Psi\,;\ f^\ast \in M^\Psi(\T)\}$.
\end{definition}

In the sequel, we shall make no distinction between $f$ and $f^\ast$, unless there may be some 
ambiguity, and shall write $f$ instead of $f^\ast$ for the boundary limit. Hence we shall allow 
ourselves to write $f(\e^{it})$ instead of $f^\ast(t)$, or even $f^\ast(\e^{it})$. Moreover, we shall 
write $\|f\|_\Psi$ instead of $\|f\|_{H^\Psi}$.\par
It follows that $H^\Psi$ becomes a subspace of $L^\Psi(\T)$ and 
$HM^\Psi = H^\Psi \cap M^\Psi(\T)$. These two spaces are closed (hence Banach spaces) since 
Proposition \ref{Hardy-Orlicz} gives:

\begin{corollaire}
$H^\Psi$ is weak-star closed in $L^\Psi=(M^\Phi)^\ast$. When $\Psi$ satisfies $\nabla_2$, it is isometrically isomorphic to the bidual of $HM^\Psi$.
\end{corollaire}

\noindent{\bf Proof.} The weak-star closure of $H^\Psi$ is obvious with Proposition \ref{Hardy-Orlicz},~2).

Suppose now that $\Phi$ satisfies $\nabla_2$, it is plainly seen that $(HM^\Psi)^\bot$ is the closed subspace of 
$L^\Phi = (M^\Psi)^\ast$ generated by all characters $e_n$ with $n<0$, where $e_n(t)= \e^{int}$ (for 
convenience, we define the duality between $f\in L^\Psi$ and $g\in L^\Phi$ by integrating the product 
$f \check{g}$, where $\check g(t) = g(-t)$). As $H^\Psi \subseteq L^\Psi = (L^\Phi)^\ast=(M^\Phi)^\ast$ is 
the orthogonal of this latter subspace. So, we have $(HM^\Psi)^{\bot\bot} = H^\Psi$.\hfill$\square$
\medskip

\noindent{\bf Proof of Proposition \ref{Hardy-Orlicz}.} Assume that $1)$ is satisfied. Since  
$\| f_r\|_1 \leq C_\Psi \| f_r\|_\Psi$, one has $f\in H^1$, and hence, by Fatou-Riesz Theorem, 
$f$ has almost everywhere a boundary limit $f^\ast \in L^1(m)$. If 
$C= \sup_{0\leq r <1} \| f_r\|_\Psi$, one has:
\begin{displaymath}
\int_\T \Psi\Big(\frac{| f_r|}{C}\Big)\,dm \leq 1
\end{displaymath}
for every $r<1$; hence Fatou's lemma implies:
\begin{displaymath}
\int_\T \Psi\Big(\frac{| f^\ast|}{C}\Big)\,dm \leq 1,
\end{displaymath}
\emph{i.e.} $f^\ast \in L^\Psi$ and $\|f^\ast\|_\Psi \leq C$.\par
Conversely, assume that $2)$ is satisfied. In particular $f^\ast \in L^1(m)$; hence 
$f\in H^1$ and $f^\ast=\lim_{r\to 1} f_r$ almost everywhere. One has $f_r = f^\ast \ast P_r$, 
where $P_r$ is the Poisson kernel at $r$. Hence, using Jensen's formula for the 
probability measure $P_r(\theta - t)\,\frac{dt}{2\pi}$, we get:
\begin{align*}
\int_{0}^{2\pi} \Psi\Big(\frac{|(f^\ast\ast P_r) (\theta)|}{\ \| f^\ast\|_\Psi}\Big)
\,\frac{d\theta}{2\pi}  
& \leq \int_{0}^{2\pi} \Psi\bigg(\int_{0}^{2\pi} 
\frac{|f^\ast(t)|}{\ \| f^\ast\|_\Psi} P_r(\theta -t)\,\frac{dt}{2\pi}\bigg)
\,\frac{d\theta}{2\pi} \\
& \leq \int_{0}^{2\pi}\bigg(\int_{0}^{2\pi}  
\Psi\Big(\frac{|f^\ast(t)|}{\ \| f^\ast\|_\Psi}\Big)
\,P_r(\theta -t)\,\frac{dt}{2\pi}\bigg)
\,\frac{d\theta}{2\pi} \\
& = \int_{0}^{2\pi}\bigg(\int_{0}^{2\pi} P_r(\theta -t)\,\frac{d\theta}{2\pi}\bigg) 
\Psi\Big(\frac{|f^\ast(t)|}{\ \| f^\ast\|_\Psi}\Big)\,\frac{dt}{2\pi} \\
& = \int_{0}^{2\pi} \Psi\Big(\frac{|f^\ast(t)|}{\ \| f^\ast\|_\Psi}\Big)\,\frac{dt}{2\pi} 
\leq 1\,,
\end{align*}
so that $\|f_r\|_\Psi \leq \|f^\ast\|_\Psi$. Hence we have $1)$, and  
$\|f\|_{H^\Psi}\leq \|f^\ast\|_\Psi$.\par 
The two parts of the proof actually give $\|f\|_{H^\Psi}= \|f^\ast\|_\Psi$.\hfill$\square$
\medskip

\begin{proposition}\label{convergence des f_r}
For every $f \in HM^\Psi$, one has 
$\| f_r - f^\ast\|_\Psi \mathop{\longrightarrow}\limits_{r\to 1} 0$. Therefore the polynomials 
on $\D$ are dense in $HM^\Psi$.
\end{proposition}
\medskip

Equivalently, on $\T=\partial \D$, the analytic trigonometric polynomial are dense in $HM^\Psi$.
\medskip

\noindent{\bf Proof.} Let $f \in HM^\Psi$ and $\varepsilon>0$. Since $M^\Psi=\overline{C(\T)}^{L^\Psi}$, there exists a continuous function $h$ on $\T$ such that $\|f-h\|_\Psi\le\varepsilon$. We have, for every $r<1$:
\begin{displaymath}
\|P_r\ast f-f\|_\Psi\le\|P_r\ast(f-h)\|_\Psi+\|P_r\ast h-h\|_\Psi+\|h-f\|_\Psi\le2\varepsilon+\|P_r\ast h-h\|_\Psi
\end{displaymath}
because $\|P_r\ast g\|_\Psi\le\|g\|_\Psi$, for every $r<1$ and every $g\in L^\Psi$.

But now, $P_r\ast h\mathop{\longrightarrow}\limits_{r\to 1} h$ uniformly. The conclusion follows. \hfill$\square$
\bigskip 

\noindent{\bf Remark.} We do not have to use a maximal function to prove the existence of 
boundary limits because we use their existence for functions in $H^1$. However, as in the 
classical case, the Marcinkiewicz interpolation Theorem, or, rather, its Orlicz space version 
ensures that the maximal non-tangential function is in $L^\Psi$. This result is undoubtedly 
known, but perhaps never stated in the following form. Recall that $N_\alpha$ is defined, 
for every $f$, say in $L^1(\T)$, as 
\begin{displaymath}
(N_\alpha f) (\e^{i\theta}) = \sup_{r\e^{it}\in S_\theta} |(f\ast P_r) (\e^{it})| = 
\sup_{z\in S_\theta} |f(z)|,
\end{displaymath}
where $S_\theta$ is the Stolz domain at $\e^{i\theta}$ with opening $\alpha$ (see \cite{Be-Sh}, page 177); 
here $f$ defines a harmonic function in $\D$.

\begin{proposition}\label{Marcinkiewicz}
Assume that the complementary function $\Phi$ of the Orlicz function $\Psi$ satisfies 
the $\Delta_2$ condition ({\it i.e.} $\Psi\in \nabla_2$). Then every linear, or sublinear, 
operator which is of weak-type $(1,1)$ and (strong) type $(\infty,\infty)$ is bounded 
from $L^\Psi$ into itself. In particular, for every $f\in L^\Psi(\T)$, the maximal 
non-tangential function $N_\alpha f$ is in $L^\Psi(\T)$ $(0<\alpha <1)$.
\end{proposition}

\noindent{\bf Proof.} If $\Psi\in \nabla_2$, then (\cite{Rao}, Theorem 3, 1~(iii), page 23), 
there exists some $\beta>1$ such that $x \Psi'(x) \geq \beta \Psi(x)$ for $x$ large enough. 
Integrating between $u$ and $v$, for $u< v$ large enough, we get 
$\frac{\Psi(u)}{\Psi(v)} \geq \big(\frac{u}{v}\big)^\beta$. Hence, for $s,t$ large 
enough $\frac{\Psi^{-1}(s)}{\Psi^{-1}(s/t)} \leq t^{1/\beta}$. This means (see 
\cite{Be-Sh}, Theorem 8.18) that the upper Boyd index of $L^\Psi$ is $\leq 1/\beta <1$. Hence 
(\cite{Be-Sh}, Theorem 5.17), $N_\alpha$ is bounded on $L^\Psi$ (it is well-known that 
$N_\alpha f$ is dominated by the Hardy-Littlewood maximal function $Mf$).
\hfill$\square$

The following, essentially well-known, criterion for compactness of operators will be very useful.

\begin{proposition}\label{critere compact} 
1) Every bounded linear operator $T\colon H^\Psi \to X$ from $H^\Psi$ into a Banach space $X$ which 
maps every bounded sequence which is uniformly convergent on compact subsets of $\D$ into a norm 
convergent sequence is compact.\par
2) Conversely, if $T\colon H^\Psi \to X$ is compact and weak-star to weak continuous, or if 
$T \colon H^\Psi \to Y^\ast$ is compact and weak-star continuous, then $T$ maps every  bounded sequence 
which is uniformly convergent on compact subsets of $\D$ into a norm convergent sequence.
\end{proposition}

Though well-known (at least for the classical case of $H^p$ spaces), the link with the weak 
(actually the weak-star) topology is usually not highlighted. Indeed, the criterion is an easy consequence 
of the following proposition. Note that Proposition \ref{critere compact} will apply to the composition 
operators on $H^\Psi$ since they are weak-star continuous.

\begin{proposition}\label{weak-star topology}
On the unit ball of $H^\Psi$, the weak-star topology is the topology of uniform convergence on every 
compact subset of $\D$.
\end{proposition}

\noindent{\bf Proof.} First we notice that the topologies are metrizable. Indeed, this is known for the topology 
of uniform convergence on every compact subset of $\D$ and, on the other hand,  $M^\Phi$ is separable, so that the weak-star topology is metrizable on the unit ball of 
its dual space $L^\Psi$, and {\it a fortiori} on that of $H^\Psi$. Now, it is sufficient to prove that the 
convergent sequences in both topologies are the same.\par
Let $(f_k)_{k \geq 1}$ be in the unit ball of $H^\Psi$ and weak-star convergent to $f\in H^\Psi$. Let us 
fix a compact subset $K$ of $\D$. We may suppose that $K$ is the closed ball of center $0$ and 
radius $r<1$. First, testing the weak-star convergence on characters, we have 
$\widehat{f_k}(n) \mathop{\longrightarrow}\limits_{k\to \infty}\hat f(n)$ for every $n\in \Z$. Then:
\begin{displaymath}
\sup_{|z|\leq r}|f_k(z) - f(z)|=\sup_{|z|=r} |f_k(z) - f(z)|
\leq \sum_{n\geq 0} r^n |\widehat{f_k}(n) -\hat f(n)|.
\end{displaymath}
The last term obviously tends to zero when $k$ tends to infinity. The result follows.\par
Conversely, let $(f_k)_{k \geq 1}$ be in the unit ball of $H^\Psi$ and converging to some holomorphic function 
$f$ uniformly on every compact subset of $\D$. We first notice that $f$ actually lies in the unit ball of 
$H^\Psi$ (by Fatou's lemma). Fix $h\in M^\Phi$  and $\eps>0$. There exists some $r<1$ such that 
$\|P_r\ast h-h\|_\Phi \leq \eps/8$, where $P_r$ is the Poisson kernel with parameter $r$. Then (see 
\cite{Rao}, page 58, inequality~$(3)$ for the presence of the coefficient $2$):
\begin{align*}
|\langle h, f_k - f \rangle |
&= |\langle P_r \ast h - h, f_k - f \rangle| + |\langle P_r\ast h, f_k -f\rangle| \\
&\leq 2\,\| P_r \ast h - h\|_\Phi \|f_k - f\|_\Psi + |\langle h , P_r\ast (f_k - f) \rangle| \\
&\leq \frac{\eps}{2}+ 2\,\|h\|_\Phi\, \|[f_k]_r - f_r\|_\Psi  \\
&\leq \frac{\eps}{2}+2\alpha\,\|h\|_\Phi \|[f_k]_r - f_r\|_\infty \\
& = \frac{\eps}{2}+2\alpha\,\|h\|_\Phi \sup_{|z|=r} |f_k(z) - f(z)| .
\end{align*}
where $\alpha$ is the norm of the injection of $L^\infty$ into $L^\Psi$.

Now, by uniform convergence on  the closed ball of center $0$ and radius $r$, there exists 
$k_\eps\geq 1$ such that for every integer $k \geq k_\eps$, one has 
\begin{displaymath}
\|h\|_\Phi \sup_{|z|=r} |f_k(z) - f(z)|\leq \eps/4.
\end{displaymath}
It follows that $(f_k)_k$ weak-star converges to $f$.\hfill$\square$
\medskip

However, we shall have to use a similar compactness criterion for Bergman-Orlicz spaces, and it is worth  
stating and proving a general criterion. We shall say that a Banach space of holomorphic functions on an 
open subset $\Omega$ of the complex plane has the \emph{Fatou property} if $X$ is continuously embedded 
(though the canonical injection) in ${\mathscr H}(\Omega)$, the space of holomorphic functions on $\Omega$, 
equipped with its natural topology of compact convergence, and if it has the following property: for every 
bounded sequence $(f_n)_n$ in $X$ which converges uniformly on compact subsets of $\Omega$ to a 
function $f$, one has $f\in X$. Then:\goodbreak

\begin{proposition} [Compactness criterion] \label{critere compacite general}
Let $X$, $Y$ be two Banach spaces of analytic functions on an open set $\Omega \subseteq \C$ which have the 
Fatou property. Let $\phi$ be an analytic self-map of $\Omega$ such that $C_\phi= f\circ \phi \in Y$ 
whenever $f\in X$. 

Then $C_\phi \colon X \to Y$ is compact if and only if for every bounded sequence 
$(f_n)_n$ in $X$ which converges to $0$ uniformly on compact subsets of $\Omega$, one has 
$\| C_\phi(f_n)\|_Y \to 0$.
\end{proposition}

Note that Hardy-Orlicz $H^\Psi$ and Bergman-Orlicz ${\mathscr B}^\Psi$ (see Section 5) spaces trivially 
have the Fatou property, because of Fatou's Lemma.
\medskip

\noindent{\bf Proof.} Assume that the above condition is fulfilled. Let $(f_n)_{n\geq 1}$ be in the unit ball 
of $X$. The assumption on $X$ implies that $(f_n)_n$ is a normal family in ${\mathscr H}(\Omega)$. 
Montel's Theorem allows us to extract a subsequence, that we still denote by $(f_n)_n$ to save notation, which 
converges to some $f \in {\mathscr H}(\Omega)$, uniformly on compact subsets of $\Omega$. 
Since $X$ has the Fatou property, one has $f\in X$. Now, since $(f_n -f )_n$ is a bounded sequence in $X$ 
which converges to $0$ uniformly on compact subsets of $\Omega$, one has 
$\| C_\phi(f_n) - C_\phi(f)\|_Y= \| C_\phi(f_n - f)\|_Y \to 0$. Hence $C_\phi$ is compact.\par
Conversely, assume that $C_\phi$ is compact. Let $(f_n)_n$ be a bounded sequence in $X$ 
which converges to $0$ uniformly on compact subsets of $\Omega$. By the compactness of $C_\phi$, 
we may assume that $C_\phi(f_n) \to g\in Y$. The space $Y$ being continuously embedded in 
${\mathscr H}(\Omega)$, $(f_n\circ \phi)_n$ converges pointwise to $g$. Since $(f_n)_n$ 
converges to $0$ uniformly on compact subsets of $\Omega$, the same is true for 
$(f_n \circ \phi)_n$. Hence $g=0$. Therefore, since $C_\phi$ is compact, we get 
$\| C_\phi(f_n)\|_Y \to 0$.\hfill $\square$

\subsection{Preliminary results}

\begin{lemme}\label{norme Psi}
Let $(\Omega, \P)$ be any probability space. For every function $g\in L^\infty (\Omega)$, one has:
\begin{displaymath}
\|g\|_\Psi \leq \frac{\|g\|_\infty}{\Psi^{-1}(\|g\|_\infty/\|g\|_1)}\cdot
\end{displaymath}
\end{lemme}

\noindent{\bf Proof.} We may suppose that $\|g\|_\infty=1$.\par
\noindent Since $\Psi(0)=0$, the convexity of $\Psi$ implies $\Psi (ax)\leq a\Psi (x)$ 
for $0\leq a\leq 1$. Hence,
for every $C>0$, one has, since $|g|\leq 1$:
\begin{displaymath}
\int_\Omega \Psi (|g|/C)\,d\P \leq \int_\Omega |g| \Psi (1/C)\,d\P =\|g\|_1\Psi (1/C).
\end{displaymath}
But $\|g\|_1 \Psi (1/C)\leq 1$ if and only if $C \geq 1/\Psi^{-1} (1/\|g\|_1)$, and
that proves the lemma.\hfill$\square$
\par\medskip\goodbreak

\begin{corollaire}\label{norme $u_{a,r}$} 
For $|a|=1$ and $0 \leq r<1$, one has:
\begin{displaymath}
\|u_{a,r}\|_\Psi \leq \frac{1}{\Psi^{-1} (\frac{1}{1-r})}\,\cdot
\end{displaymath}
\end{corollaire}

\noindent{\bf Proof.} One has $\|u_{a,r}\|_\infty=1$, and:
\begin{align*}
\|u_{a,r}\|_1 & =\int_0^{2\pi} \Big| \frac{1-r}{1- \bar{a} r \e^{it}}\Big|^2\,dm(t) \\
& = (1-r)^2\sum_{n=0}^{+\infty} r^{2n}
=\frac{(1-r)^2}{1-r^2}= \frac{1-r}{1+r} \cdot
\end{align*}
Hence $\|u_{a,r}\|_\Psi \leq 1/\Psi^{-1}(1+r/1-r)$, by using Lemma \ref{norme Psi}, giving the result since 
$(1+r)/(1 - r )\geq 1/1 - r$.\hfill $\square$

\medskip
\noindent{\bf Remark.} We hence get actually $\|u_{a,r}\|_\Psi \leq 1/\Psi^{-1}(1+r/1-r)$; the term 
$1 + r$ has no important meaning, so we omit it in the statement of Corollary \ref{norme $u_{a,r}$}, but 
sometimes, for symmetry of formulae, or in order to be in accordance with the classical case, we shall use
this more precise estimate.
\bigskip

For every $f\in L^1(\T)$ and every $z=r\,\e^{i\theta} \in \D$, one has
\begin{displaymath}
(P[f])(z)= \int_0^{2\pi} f(\e^{it})P_z(t)\,dm(t),
\end{displaymath}
where $P_z$ is the Poisson kernel:
\begin{displaymath}
P_z(t)= \frac{1-r^2}{1-2r\cos(\theta -t) +r^2} 
= \frac{1 - |z|^2}{|\e^{it} - z|^2}\,\raise 0,5mm\hbox{,}
\end{displaymath}
and $f(z)= (P[f])(z)$ when $f$ is analytic on $\D$. Since 
$P_z\in L^\infty(\T)\subseteq L^\Phi(\T)$, it follows
that the evaluation in $z\in \D$:
\begin{displaymath}
\delta_z(f)=f(z)
\end{displaymath}
is a continuous linear form on $H^\Psi$. The following lemma explicits the behaviour of 
its norm.

\begin{lemme}\label{norme evaluation}
For $|z|<1$, the norm of the evaluation functional at $z$ is:
\begin{displaymath}
\|\delta_z\|_{(HM^\Psi)^\ast} =\|\delta_z\|_{(H^\Psi)^\ast} \approx \Psi^{-1}\Big(\frac{1}{1 -|z|}\Big)\,\cdot
\end{displaymath}
More precisely:
\begin{displaymath}
\frac{1}{4} \Psi^{-1}\Big(\frac{1 + |z|}{1 -|z|}\Big)
 \leq \|\delta_z\|_{(H^\Psi)^\ast} 
 \leq 2\Psi^{-1}\Big(\frac{1 + |z|}{1 -|z|}\Big)\,\cdot
\end{displaymath}
\end{lemme}

\noindent{\bf Remark.} In particular:
\begin{displaymath}
\frac{1}{4} \Psi^{-1}\Big(\frac{1}{1 -|z|}\Big) 
\leq \|\delta_z\|_{(H^\Psi)^\ast}
\leq 4\Psi^{-1}\Big(\frac{1}{1 -|z|}\Big)\,,
\end{displaymath}
which often suffices for our purpose.
\bigskip

\noindent{\bf Proof.} The first equality $\|\delta_z\|_{(HM^\Psi)^\ast} =\|\delta_z\|_{(H^\Psi)^\ast}$ comes from the fact that $f_r\in HM^\Psi$, for every $f\in H^\Psi$ and $r<1$ (thus $f(rz)\mathop{\longrightarrow}\limits_{r\to 1}f(z)$, when $z\in\D$ and $f\in H^\Psi$).

On the one hand, we have, when $|z|=r$, using \cite{Rao}, inequality $(4)$ page 58, and 
Lemma \ref{norme Psi}, since $\|P_z\|_1=1$ and $\|P_z\|_\infty=\frac{1+r}{1-r}$:
\begin{displaymath}
\|\delta_z\|_{(H^\Psi)^\ast} \leq 2\,\|P_z\|_\Phi 
\leq 2\frac{1+r}{1-r}\frac{1}{\Phi^{-1} \big(\frac{1+r}{1-r}\big)}\,\raise0,5mm\hbox{,}
\end{displaymath}
which is less than $2\Psi^{-1}(1+r/1-r)$, by using the inequality (see 
\cite{Rao}, Proposition 1~(ii), page 14, or \cite{Krasno}, pages 12--13):
\begin{displaymath}
\Psi^{-1}(x) \Phi^{-1}(x) \geq x\,, \hskip 5mm x>0.
\end{displaymath}
On the other hand, one has, using Corollary \ref{norme $u_{a,r}$}, with $r=|z|$ and $\bar{a} z=r$:
\begin{displaymath}
\|\delta_z\|_{(H^\Psi)^\ast} \geq \frac{|u_{a,r}(z)|}{\,\,\|u_{a,r}\|_\Psi} 
\geq \frac{1 /(1+r)^2}{1/\Psi^{-1}\big(\frac{1 + r}{1 - r}\big)}
\geq \frac{1}{4} \Psi^{-1}\Big(\frac{1 + r}{1-r}\Big),
\end{displaymath}
and that ends the proof.\hfill $\square$
\par\medskip\goodbreak

\subsection{Composition operators}

We establish now some estimations for the norm of composition operators.

\begin{proposition}\label{subordination}
1) Every analytic self-map $\phi\colon \D \to \D$ induces a bounded composition operator 
$C_\phi \colon H^\Psi \to H^\Psi$ by setting $C_\phi (f) = f\circ \phi$. More precisely:
\begin{displaymath}
\| C_\phi\| \leq \frac{1 + | \phi (0)| } {1 - |\phi (0)|}\,\cdot
\end{displaymath}
In particular, $\|C_\phi\| \leq 1$ if $\phi (0)=0$.\par

2) One has:
\begin{displaymath}
\|C_\phi\| \geq \frac{1}{8 \Psi^{-1} (1)} \Psi^{-1}\bigg(\frac{1 + | \phi (0)| } {1 - |\phi (0)|}\bigg).
\end{displaymath}

3) When $\Psi \in \nabla_1$ globally: $\Psi (x) \Psi (y) \leq \Psi (b xy)$ for all $x,y \geq 0$, we also have:
\begin{displaymath}
\| C_\phi\| \leq b \Psi^{-1} \bigg(\frac{1 + | \phi (0)| } {1 - |\phi (0)|}\bigg).
\end{displaymath}
\par

4) Moreover, $C_\phi$ maps $HM^\Psi$ into $HM^\Psi$. Hence, if $\Psi\in\nabla_2$, then $C_\phi \colon H^\Psi \to H^\Psi$ is 
the bi-adjoint of the composition operator $C_\phi \colon HM^\Psi \to HM^\Psi$.
\end{proposition}

Note that when $\Psi (x) = x^p$ for $1 \leq p < \infty$, then $\Psi \in \nabla_1$ globally, with $b=1$.
\medskip

\noindent{\bf Proof.} 1) Assume first that $\phi (0)=0$. Let $f\in H^\Psi$, with $\|f\|_\Psi=1$. 
Since $\Psi$ is convex and increasing, the function $u = \Psi \circ |f|$ is subharmonic on $\D$, 
thanks to Jensen's inequality. The condition $\phi(0)=0$ allows to use Littlewood's 
subordination principle (\cite{Du}, Theorem 1.7); for $r<1$, one has:
\begin{displaymath}
\int_0^{2\pi} \Psi\big( |(f \circ \phi)(r \e^{it})|\big)\,\frac{dt}{2\pi} 
\leq \int_0^{2\pi} \Psi\big( |f(r \e^{it})|\big)\,\frac{dt}{2\pi} \leq 1.
\end{displaymath}
Hence $f\circ \phi\in H^\Psi$ and $\| f\circ \phi\|_\Psi\leq 1$.\par

Assume now that $\phi$ is an inner function, and let $a=\phi(0)$. It is known that (see 
\cite{Nord}, Theorem 1) that 
\begin{displaymath}
\phi^\ast(m) = P_a.m\,,
\end{displaymath}
where $\phi^\ast(m)$ is the image under $\phi^\ast$ (the boundary limit of $\phi$) of the 
normalized Lebesgue measure $m$, and $P_a.m$ is the measure of density $P_a$, the Poisson 
kernel at $a$. Therefore, for every $f\in H^\Psi$ with $\| f\|_\Psi=1$, one has for 
$0\leq r<1$, in setting $K_a=\|P_a\|_\infty=\frac{1 + |a|}{1-|a|}$:
\begin{align}\label{cas inner}
\hskip 1cm \int_0^{2\pi} \Psi\Big( \frac{|(f  \circ \phi) (r\e^{it})|}{K_a}\Big)\,\frac{dt}{2\pi} 
& \leq \int_\T \Psi\Big(\frac{|(f \circ \phi)^\ast|}{K_a}\Big)\,dm \notag \\
& = \int_\T \Psi\Big(\frac{|f^\ast \circ \phi^\ast|}{K_a}\Big)\,dm \hskip 3mm 
\text{(recall that $|\phi^\ast|=1$)} \notag \\
& = \int_\T \Psi\Big(\frac{|f^\ast|}{K_a}\Big)\,d\phi^\ast(m) 
= \int_\T \Psi\Big(\frac{|f^\ast|}{K_a}\Big)\,P_a\,dm \\
& \leq \int_\T \frac{1}{K_a}\Psi(|f^\ast|)\,P_a\,dm\,, \hskip 3mm 
\text{since $K_a>1$} \notag \\
& \leq \int_\T \frac{1}{K_a}\Psi(|f^\ast|)\,\|P_a\|_\infty\,dm \notag \\
& = \int_\T \Psi(|f^\ast|)\,dm \leq 1. \notag
\end{align}
Hence $\|(f\circ \phi)_r\|_\Psi \leq K_a$, and therefore $\| f\circ \phi\|_\Psi\leq K_a$.\par

Then, for an arbitrary $\phi$, let $a=\phi(0)$ again, and let $\phi_a$ be the automorphism 
$z\mapsto \frac{z -a }{1 -\bar{a}z}$\raise0,5mm\hbox{,} whose inverse is $\phi_{-a}$. Since 
$\phi=\phi_{-a}\circ (\phi_a \circ \phi)$, one has 
$C_\phi= C_{\phi_a\circ \phi}\circ C_{\phi_{-a}}$. But $\phi_{-a}$ is inner and, on the other 
hand, $(\phi_a\circ \phi)(0)=0$; hence parts a) and b) of the proof give:
\begin{displaymath}
\| C_\phi\| \leq \| C_{\phi_{-a}}\| \leq K_a =\frac{1+|a|}{1-|a|}\,\raise 0,5mm\hbox{,}
\end{displaymath} 
which gives the first part of the proof.\par
\smallskip

2) By Lemma~\ref{norme evaluation}, we have for every $f \in H^\Psi$ with $\|f\|_\Psi \leq 1$:
\begin{displaymath}
| (f \circ \phi ) (0)| \leq \|\delta_0\|_{(H^\Psi)^\ast} \| f \circ \phi \|_\Psi 
\leq 2 \Psi^{-1} (1)\,\|C_\phi\|.
\end{displaymath}
In other words:
\begin{displaymath}
| f\big( \phi (0) \big)| \leq 2 \Psi^{-1} (1) \| C_\phi \|
\end{displaymath}
for every such $f \in H^\Psi$. Hence:
\begin{displaymath}
\|\delta_{\phi (0)}\|_{(H^\Psi)^\ast} \leq 2 \Psi^{-1} (1)  \,\|C_\phi\|\,,
\end{displaymath}
giving
\begin{displaymath}
\|C_\phi\| \geq \frac{1}{8 \Psi^{-1} (1)} \Psi^{-1}\bigg(\frac{1 + | \phi (0)| } {1 - |\phi (0)|}\bigg),
\end{displaymath}
by using Lemma~\ref{norme evaluation} again, but the minoration.\par 
\smallskip

3) When $\Psi \in \nabla_1$ globally, we go back to 
the proof of $1)$. We have only to modify inequalities \eqref{cas inner} in b). Setting  
$K'_a = \Psi^{-1} (K_a)$, and writing $P_a = \Psi \big( \Psi^{-1} (P_a)\big)$, we get, 
for every $f \in H^\Psi$ with $\| f\|_\Psi =1$:
\begin{align*}
\int_0^{2\pi} \Psi\Big( \frac{|(f  \circ \phi) (r\e^{it})|}{b K'_a}\Big)\,\frac{dt}{2\pi} 
& \leq \int_\T \Psi\Big(\frac{|f^\ast|}{b K'_a}\Big)\,P_a\,dm \\
& = \int_\T \Psi\Big(\frac{|f^\ast|}{b K'_a}\Big)\, \Psi \big( \Psi^{-1} (P_a)\big)\,dm \\ 
& \leq \int_\T \Psi\Big(\frac{|f^\ast|}{K'_a} \Psi^{-1} (P_a) \Big)\,dm \\
& \leq \int_\T \Psi (|f^\ast|)\,dm \leq 1,
\end{align*}
since $\Psi^{-1} (P_a) \leq \Psi^{-1} ( \|P_a\|_\infty) = K'_a$, giving $\| f \circ \phi \|_\Psi \leq b K'_a$.
\par\smallskip

4) Suppose now that $f\in HM^\Psi$. As before, when $\phi(0) = 0$, Littlewood's subordination principle 
gives, for every $C>0$:
\begin{align*}
\int_0^{2\pi} \Psi \big(C |(f & \circ \phi)(\e^{it})|\big)\,\frac{dt}{2\pi}  
= \sup_{r<1} \int_0^{2\pi} \Psi \big(C |(f\circ \phi)(r\e^{it})|\big)\,\frac{dt}{2\pi} \\
& \leq \sup_{r<1} \int_0^{2\pi} \Psi \big(C |f (r\e^{it})|\big)\,\frac{dt}{2\pi} 
= \int_0^{2\pi} \Psi \big(C |f (\e^{it})|\big)\,\frac{dt}{2\pi} < +\infty\,;
\end{align*} 
hence $f\circ \phi \in HM^\Psi$.  When $\phi$ is inner, the same computations as in 1)~b) above, using 
that $\phi^\ast(m)= P_a.m$, where $a=\phi(0)$, give, for every $C>0$:
\begin{align*}
\int_0^{2\pi} \Psi \big(C |(f\circ \phi)(r\e^{it})|\big)\,\frac{dt}{2\pi} 
& \leq \int_\T \Psi( C |f^\ast|) \|P_a\|_\infty\,dm \\
& = K_a \int_\T \Psi( C |f^\ast|)\,dm <+\infty,
\end{align*}
and $f\circ \phi \in HM^\Psi$ again. The general case follows, as in 1)~c) above, since $f\in HM^\Psi$ 
implies $f \circ \phi_{-a} \in HM^\Psi$, because $\phi_{-a}$ is inner, and then 
$f\circ \phi = (f\circ \phi_{-a})\circ (\phi_a \circ \phi) \in HM^\Psi$ since 
$(\phi_a \circ \phi)(0)=0$.
\hfill$\square$

\subsection{Order bounded composition operators}

Recall that an operator $T\colon X \to Z$ from a Banach space $X$ into a Banach subspace 
$Z$ of a Banach lattice $Y$ is order bounded if there is some positive $y\in Y$ such that 
$| Tx | \leq y$ for every $x$ in the unit ball of $X$.
\par\medskip

Before studying order bounded composition operators, we shall recall the following, certainly well-known,  
result, which says that order boundedness can be seen as stronger than compactness.

\begin{proposition}
Let $T\colon L^2(\mu) \to L^2(\mu)$ be a continuous linear operator. Then $T$ is order bounded 
if and only if it is a Hilbert-Schmidt operator. 
\end{proposition}

The proof is straightforward: if $B$ is the unit ball of $L^2(\mu)$, and $(e_i)_i$ is an 
orthonormal basis, one has $\sup_{f\in B} |Tf| = \big(\sum_i |Te_i|^2\big)^{1/2}$. Hence 
$\sup_{f\in B} |Tf|\in L^2(\mu)$ if and only if 
$\int \big(\sum_i |Te_i|^2\big)\,d\mu = \sum_i \| Te_i\|^2 < +\infty$, \emph{i.e.} if and only if 
$T$ is Hilbert-Schmidt.\par
J. H. Shapiro and P. D. Taylor proved in \cite{Shap-T} that there exist composition operators 
on $H^2$ which are compact but not Hilbert-Schmidt. We are going to see that, when the 
Orlicz function $\Psi$ grows fast enough, the compactness of composition operators on $H^\Psi$ is 
equivalent to their order boundedness.

\begin{proposition}\label{equiv order bounded}
The composition operator $C_\phi \colon H^\Psi \to H^\Psi$ is order bounded  (resp. order bounded into
$M^\Psi(\T)$) if and only if $\Psi^{-1}\big(\frac{1}{1-|\phi|}\big)\in L^\Psi(\T)$ (resp. 
$\Psi^{-1}\big(\frac{1}{1-|\phi|}\big)\in M^\Psi(\T)$). Equivalently, if and only if
\begin{equation}\label{integrabilite-1}
\chi_A\Big(\frac{1}{1-|\phi|}\Big)\in L^1(\T) \hskip 3mm \text{for some } A>0,
\tag*{(OB1)}
\end{equation}
respectively:
\begin{equation}\label{integrabilite}
\chi_A\Big(\frac{1}{1-|\phi|}\Big)\in L^1(\T) \hskip 3mm \text{for every } A>0,
\tag*{(OB2)}
\end{equation}
\end{proposition}

\noindent In other words (recall that $\chi_A(x)=\Psi\big(A \Psi^{-1}(x)\big)$), if and 
only if:
\begin{displaymath}
\int_\T \Psi\bigg[A\Psi^{-1}\Big(\frac{1}{1-|\phi|}\Big)\bigg]\,dm < +\infty 
\hskip 3mm \text{for some
(resp. every) $A>0$.}
\end{displaymath}
\medskip

\noindent{\bf Remark} For a sequence $(g_n)_n$ of elements of $L^\Psi(\Omega)$, one has
$\|g_n\|_\Psi \mathop{\longrightarrow}\limits_{n\to+\infty} 0$
if and only if
\begin{equation}\label{norme psi tend vers 0}
\hskip 1,5cm \int_\Omega \Psi\Big(\frac{|g_n|}{\eps}\Big)\,d\P \mathop{\longrightarrow}_{n\to+\infty} 0 
\hskip 3mm \text{for every $\eps >0$}.
\end{equation}
In fact, if \eqref{norme psi tend vers 0} holds, then for $n\geq n_\eps$, the above integrals are $\leq 1$, and 
hence $\|g_n\|_\Psi \leq \eps$. Conversely, assume that 
$\|g_n\|_\Psi \mathop{\longrightarrow}\limits_{n\to+\infty} 0$, and let $\eps>0$ be given. Let 
$0< \delta \leq 1$. Since 
$\|g_n/(\eps \delta)\|_\Psi \mathop{\longrightarrow}\limits_{n\to+\infty} 0$, one has   
$\|g_n/(\eps \delta)\|_\Psi \leq 1$, and hence 
$\int_\Omega \Psi\big(\frac{|g_n |}{\eps \delta}\big)\,d\P \leq 1$, for $n$ large enough. Then, using the 
convexity of $\Psi$:
\begin{displaymath}
\int_\Omega \Psi\Big( \frac{|g_n|}{\eps}\Big)\,d\P 
= \int_\Omega \Psi\Big( \delta\,\frac{|g_n|}{\eps \delta }\Big)\,d\P 
\leq \delta \int_\Omega \Psi\Big( \frac{|g_n|}{\eps\delta }\Big)\,d\P \leq \delta\,,
\end{displaymath}
for $n$ large enough.
\par\smallskip

Therefore, using Lebesgue's dominated convergence Theorem, it follows that if 
$C_\phi \colon H^\Psi \to H^\Psi$ is order bounded \emph{into $M^\Psi(\T)$}, then the 
composition operator $C_\phi  \colon H^\Psi \to H^\Psi$ is compact.
\medskip

\noindent{\bf Proof of Proposition \ref{equiv order bounded}.} 

As $HM^\Psi$ is separable, there exists a sequence $(f_n)_{n\ge1}$ in the unit ball of $HM^\Psi$ such that 
$$\sup_n\bigl|f_n\circ\phi(r\e^{i\theta})\bigr|=\|\delta_{\phi(r \e^{i\theta})}\|_{(HM^\Psi)^\ast}.$$

Now, suppose that $C_\phi$ is order bounded into $L^\Psi(\T)$ (resp. into $M^\Psi(\T)$). Then there exists some $g$ in  $L^\Psi(\T)$ (resp. in $M^\Psi(\T)$) such that $g\ge|C_\phi(f)|$ {\it a.e.}, for every $f$ in the unit ball of $H^\Psi$. Using Lemma \ref{norme evaluation}, we have {\it a.e.}
$$\Psi^{-1}\Bigl(\frac{1}{1-|\phi(r\e^{i\theta})|}\Bigr)\le 4\|\delta_{\phi(r \e^{i\theta})}\|_{(H^\Psi)^\ast}=4\sup_n\bigl|f_n\circ\phi(r\e^{i\theta})\bigr|\le4g.$$
The result hence follows letting $r\uparrow 1$.

The converse is obvious.
\hfill$\square$
\par\medskip

\begin{theoreme}\label{impli-order bounded}
If the composition operator $C_\phi \colon H^\Psi \to H^\Psi$ is order bounded, then:
\begin{equation}\label{grand O-1}
m(1-|\phi| < \lambda )= O\,\Big(\frac{1}{\chi_A(1/\lambda)}\Big),
\hskip 2mm \text{as}\ \lambda\to 0,\hskip 3mm \text{for some } A>0.
\tag*{(OB3)}
\end{equation}
and if it is order bounded into $M^\Psi(\T)$, then:
\begin{equation}\label{grand O}
m(1-|\phi| < \lambda )= O\,\Big(\frac{1}{\chi_A(1/\lambda)}\Big),
\hskip 2mm \text{as}\ \lambda\to 0,\hskip 3mm \text{for every } A>0.
\tag*{(OB4)}
\end{equation}
Under the hypothesis $\Psi\in \Delta^1$, the converse holds.
\end{theoreme}

\par\medskip

\noindent{\bf Proof.} The necessary condition follows from Proposition~\ref{equiv order bounded} 
and Markov's inequality:
\begin{displaymath}
m\Big(\frac{1}{1- |\phi|} >t \Big) \leq \frac{1}{\chi_A (t)} 
\int_\T \chi_A\Big(\frac{1}{1- |\phi|}\Big)\,dm.
\end{displaymath}
\par\medskip

For the converse, we shall prove a stronger result, and for that, we define the 
\emph{weak-$L^\Psi$ space}
as follows:

\begin{definition}
The weak-$L^\Psi$ space $L^{\Psi,\infty}(\Omega)$ is the space of measurable functions 
$f\colon \Omega\to \C$
such that, for some constant $c>0$, one has, for every $t>0$:
\begin{displaymath}
\P(|f| >t) \leq \frac{1}{\Psi(ct)}\,\cdot
\end{displaymath}
\end{definition}
\medskip

For subsequent references, we shall state separately the following elementary result.

\begin{lemme}\label{Markov} 
For every $f \in L^\Psi (\Omega)$, one has, for every $t>0$:
\begin{displaymath}
\| f \|_\Psi \geq \frac{t\,}{\Psi^{-1}\big(\frac{1}{\P (|f| > t)}\big)}\,\cdot
\end{displaymath}
\end{lemme}

\noindent{\bf Proof.}  By  Markov's inequality, one has, for $t>0$:
\begin{displaymath}
\Psi \Big(\frac{t}{\ \|f \|_\Psi}\Big) \P (|f| >  t) 
\leq \int_\Omega \Psi\Big(\frac{|f|\ }{\ \|f \|_\Psi}\Big)\,d\P \leq 1;
\end{displaymath}
and that gives the lemma.\hfill$\square$
\bigskip

Since Lemma~\ref{Markov} can be read:
\begin{displaymath}
\P(|f|>t) \leq \frac{1}{\Psi(t/\|f\|_\Psi)}\,\raise0,5mm\hbox{,}
\end{displaymath}
we get that $L^\Psi(\Omega)\subseteq L^{\Psi,\infty}(\Omega)$.\par\medskip

The converse of Theorem \ref{impli-order bounded} is now an immediate consequence 
of the following
proposition.

\begin{proposition}\label{egalite L-Psi}
\vbox{\hfill}\par
$1)\, a)$ If $\Psi\in \Delta^1$, then $L^\Psi(\Omega) = L^{\Psi,\infty}(\Omega)$.\par
$\phantom{1)}\, b)$ If  $L^\Psi(\Omega) = L^{\Psi,\infty}(\Omega)$, then $\Psi\in\Delta^0$.\par
$2)$ If $L^\Psi(\T) = L^{\Psi,\infty}(\T)$, then condition \ref{grand O-1} implies that 
$C_\phi\colon H^\Psi \to H^\Psi$ is order bounded, and condition \ref{grand O} implies that it 
is order bounded into $M^\Psi(\T)$.
\end{proposition}

\begin{lemme}\label{chi-integrable}

The following assertions are equivalent

i) $L^\Psi(\Omega) = L^{\Psi,\infty}(\Omega)$.\par

ii)  $\displaystyle\int_1^\infty \frac{\Psi' (u)}{\Psi (B u)}\,du\equiv \int_{\Psi(1)}^{+\infty} \frac{1}{\chi_B(x)}\,dx<+\infty$, 
for some $B>1$.\par
\end{lemme}

\noindent{\bf Proof of the Lemma.} Assume first that $1/\chi_B$ is integrable on 
$(\Psi(1),\infty)$. For
every $f\in L^{\Psi,\infty}(\Omega)$, there is a $c>0$ such that 
$\P(|f|>t)\leq 1/\Psi(c t)$. Then
\begin{align*}
\int_\Omega \Psi\Big(c\frac{|f|}{B}\Big)\,d\P 
&=\int_0^{+\infty} \P(|f|>Bt/c)\,\Psi'(t)\,dt \\
&\leq\Psi(1)+\int_1^{+\infty} \P(|f|>Bt/c)\,\Psi'(t)\,dt \\
& \leq\Psi(1)+\int_1^{+\infty} \frac{\Psi'(t)}{\Psi(Bt)}\,dt 
=\Psi(1)+\int_{\Psi(1)}^{+\infty} \frac{du}{\chi_B(u)} <+\infty,
\end{align*}
so that $f\in L^\Psi(\Omega)$.\par
Conversely, assume that $L^\Psi(\Omega) = L^{\Psi,\infty}(\Omega)$. Since $1/\Psi$ 
is decreasing, there is a measurable
function $f\colon \Omega \to \C$ such that $\P(|f|>t)=1/\Psi(t)$, when $t\ge\Psi^{-1}(1)$. Such a function is 
in $L^{\Psi,\infty}(\Omega)$;
hence it is in $L^\Psi(\Omega)$, by our hypothesis. Therefore, there is a $B>1$ such that
\begin{displaymath}
\int_\Omega \Psi(|f|/B)\,d\P <+\infty.
\end{displaymath}
But
\begin{align*}
\int_\Omega \Psi(|f|/B)\,d\P & = \int_0^{+\infty} \P(|f|>Bt)\,\Psi'(t)\,dt \ge \int_{\frac{1}{B}\Psi^{-1}(1)}^{+\infty} \frac{\Psi'(t)}{\Psi(Bt)}\,dt  \\
& \ge \int_{\Psi^{-1}(1)}^{+\infty} \frac{\Psi'(t)}{\Psi(Bt)}\,dt 
= \int_1^{+\infty} \frac{du}{\chi_B(u)}\,\raise0,5mm\hbox{,}
\end{align*}
and hence $1/\chi_B$ is integrable on $(1,\infty)$.\hfill$\square$
\par\medskip

\noindent{\bf Proof of Proposition \ref{egalite L-Psi}.} 1)a) We first remark that for 
every Orlicz function
$\Psi$, one has $\Psi(2x)\geq x\Psi'(x)$ for every $x>0$, because, since $\Psi'$ is 
positive and increasing:
\begin{displaymath}
\Psi(2x) =\int_0^{2x} \Psi'(t)\,dt \geq \int_x^{2x} \Psi'(t)\,dt \geq x\Psi'(x).
\end{displaymath}
Assume now that $\Psi\in \Delta^1$: $x\Psi(x)\leq \Psi(\beta x)$
for some $\beta>0$ and for $x\geq x_0>0$. By Lemma \ref{chi-integrable}, it suffices to 
show that $1/\chi_{2\beta}$
is integrable on $(\Psi(x_0),+\infty)$. But
\begin{displaymath}
\int_{\Psi(x_0)}^{+\infty} \frac{dx}{\chi_{2\beta}(x)}
= \int_{x_0}^{+\infty} \frac{\Psi'(u)}{\Psi(2\beta u)}\,du
\leq \int_{x_0}^{+\infty} \frac{\Psi(2u)/u}{2u\Psi(2u)}\,du =
\int_{x_0}^{+\infty} \frac{du}{2u^2} <+\infty.
\end{displaymath}
\par
1)b) Suppose now that $L^\Psi(\Omega) = L^{\Psi,\infty}(\Omega)$. By the preceding lemma, there exists some $B>1$ such that 

\begin{displaymath}
\mathop{\lim}\limits_{x\to +\infty}\int_x^{2x}\frac{\Psi' (u)}{\Psi (B u)}\,du=0.
\end{displaymath}
By convexity, $\Psi(2x)\ge2\Psi(x)$ and $\Psi'$ is nonnegative, so that 
$$\int_x^{2x}\frac{\Psi' (u)}{\Psi (B u)}\,du\ge\frac{1}{\Psi (2Bx)}\int_x^{2x}\Psi'(u)\,du\ge\frac{\Psi (2x)-\Psi (x)}{\Psi (2Bx)}\ge\frac{\Psi (x)}{\Psi (2Bx)}\cdot$$

Therefore $\Psi$ satisfies $\Delta^0$.
\par
2) Assume that $L^\Psi(\T) = L^{\Psi,\infty}(\T)$ and that condition \ref{grand O-1} 
(resp. \ref{grand O}) holds. 
By Lemma \ref{chi-integrable}, there is a $B>0$ such that $1/\chi_B$ is integrable on $(1,+\infty)$. We get, 
using condition \ref{grand O-1} (resp. \ref{grand O}) and setting $C=A/B$:
\begin{align*}
\int_\T \chi_C \Big(\frac{1}{1-|\phi|}\Big)\,dm
& = \int_0^{+\infty} m\Big(\frac{1}{1-|\phi|} > t\Big) \chi'_C(t)\,dt \\
& = \int_0^{+\infty} m(1-|\phi| <1/ t)\, \chi'_C (t)\,dt\\
&\leq \chi_C (1)+K \int_1^{+\infty} \frac{\chi'_C (t)}{\chi_{A} (t)}\,dt.
\end{align*}
But, if we set $u= \chi_C(t)$, \emph{i.e.} $u=\Psi\big(C\Psi^{-1}(t)\big)$, then 
$\Psi^{-1}(u)= C \Psi^{-1}(t)$,
and hence
\begin{displaymath}
\chi_{A}(t)= \Psi\big( A \Psi^{-1}(t)\big) = \Psi\big( B \Psi^{-1}(u)\big) = \chi_B(u).
\end{displaymath}
Therefore:
\begin{displaymath}
\int_\T \chi_C\Big(\frac{1}{1-|\phi|}\Big)\,dm 
\leq \chi_C (1)+K\int_{\chi_C (1)}^{+\infty} \frac{du}{\chi_B(u)}\,du <+\infty.
\end{displaymath}
It follows from Proposition~\ref{equiv order bounded} that $C_\phi$ is order bounded (resp. into
$M^\Psi(\T)$).\par\hfill$\square$
\par\medskip

\noindent{\bf Remark.} The condition $\Psi\in \Delta^1$ is not equivalent to 
$L^\Psi(\Omega) =L^{\Psi,\infty}(\Omega)$.
For example, we can take:
\begin{displaymath}
\Psi(x)= \exp\big[\big(\log (x+1)\big)^{3/2}\big] -1.
\end{displaymath}
Then, as $x$ tends to infinity,
\begin{displaymath}
\Psi(K x)\sim \Psi(x)\,\exp\big(\frac{3}{2}(\log K )(\log x)^{1/2}\big),
\end{displaymath}
and hence $\Psi\not \in \Delta^1$. On the other hand, 
$\Psi'(x)\sim \frac{3}{2}\frac{(\log x)^{1/2}}{x}\Psi(x)$; hence:
\begin{align*}
\int_{\Psi(1)}^{+\infty} \frac{du}{\chi_K(u)} & = \int_1^{+\infty} \frac{\Psi'(t)}{\Psi(Kt)}\,dt \\
& \sim \int_0^{+\infty} \frac{3u^2}{\exp(\frac{3}{2} (\log K) u)}\,du
= \int_0^{+\infty} \frac{3 u^2}{K^{3u/2}}\,du <+\infty
\end{align*}
for $K>1$. Therefore $L^\Psi(\Omega) =L^{\Psi,\infty}(\Omega)$ by Lemma \ref{chi-integrable}.
\bigskip

\subsection{Weakly compact composition operators}

We saw in Lemma \ref{norme $u_{a,r}$} that
\begin{displaymath}
\|u_{a,r}\|_\Psi \leq \frac{1}{\Psi^{-1} (\frac{1}{1-r})}\,\cdot
\end{displaymath}
The next result shows that the weak compactness of $C_\phi$ transforms the ``\emph{big-oh}'' into a
``\emph{little-oh}'', when $\Psi$ grows fast enough.  

\begin{theoreme}\label{faiblement compact}
Assume that $\Psi\in \Delta^0$. Then the weak compactness of the composition operator
$C_\phi \colon H^\Psi \to H^\Psi$ implies that:
\begin{equation}\label{petit o}
\sup_{a\in \T} \|C_\phi (u_{a,r})\|_\Psi
=o\,\bigg(\frac{1}{\Psi^{-1}\big(\frac{1}{1-r}\big)}\bigg),
\hskip 3mm \text{as}\ r\to 1.
\tag*{(W)}
\end{equation}
\end{theoreme}

\noindent{\bf Proof.} Actually, we only need to use that the 
restriction of $C_\phi\colon H^\Psi \to H^\Psi$ to $HM^\Psi$ is weakly compact. We proved in 
\cite{LLQR}, Theorem 4, that, under the hypothesis $\Psi\in \Delta^0$, the 
operator $C_\phi \colon HM^\Psi \to HM^\Psi$ is weakly compact if and only if for every $\eps>0$, we 
can find $K_\eps>0$ such that, for every $f\in HM^\Psi$, one has:
\begin{displaymath}
\|C_\phi(f)\|_\Psi \leq K_\eps \|f\|_1 +\eps\,\|f\|_\Psi.
\end{displaymath}
Using Corollary \ref{norme $u_{a,r}$}, we get, for every $\eps>0$, since $\|u_{a,r}\|_1\leq 1-r$:
\begin{displaymath}
\|C_\phi(u_{a,r})\|_\Psi 
\leq K_\eps (1-r) +\eps\,\frac{1}{\Psi^{-1}(\frac{1}{1-r})}\,\cdot
\end{displaymath}
But $\frac{\Psi(x)}{x}\mathop{\longrightarrow}\limits_{x\to +\infty} +\infty$. Hence
\begin{displaymath}
1-r = o\Big(\frac{1}{\Psi^{-1}(\frac{1}{1-r})}\Big) \hskip 3mm \text{as } r\to 1,
\end{displaymath}
and that proves the theorem.\hfill$\square$
\medskip

We shall see in Section 4 that the converse holds when $\Psi$ satisfies $\Delta^0$ and moreover that $C_\phi$ is then compact. That will use some other techniques. Nevertheless, we can prove, from now on, the following result.

\begin{theoreme}\label{Delta2}
If $\Psi\in \Delta^2$, then condition:
\begin{equation}
\sup_{a\in \T} \|C_\phi (u_{a,r})\|_\Psi
=o\,\bigg(\frac{1}{\Psi^{-1}\big(\frac{1}{1-r}\big)}\bigg),
\hskip 3mm \text{as}\ r\to 1.\tag*{\ref{petit o} }
\end{equation}
in Theorem~\ref{faiblement compact} implies condition:
\begin{equation}
m(1-|\phi| < \lambda )= O\,\Big(\frac{1}{\chi_A(1/\lambda)}\Big),
\hskip 2mm \text{as}\ \lambda\to 0,\hskip 3mm \text{for every } A>0.\tag*{\ref{grand O} }
\end{equation}
 in Theorem~\ref{impli-order bounded}.
\end{theoreme}
\medskip

Remark that when $\Psi\in \Delta^0$ (and in particular when $\Psi\in \Delta^2$), one has, 
for any $B> \beta A$ (where $\beta$ is given by the definition of $\Delta^0$):
$\frac{\Psi(Bx)}{\Psi(Ax)}\to +\infty$ as $x\to +\infty$; hence
$1/\chi_A(x) =o\big(1/\chi_B(x)\big)$ and the ``\emph{big-oh}'' condition in \ref{grand O} 
may be replaced by a ``\emph{little-oh}'' condition.
\par\medskip

Before proving this theorem, let us note that:

\begin{proposition}\label{mesure nulle}
Condition
\begin{equation}
\sup_{a\in \T} \|C_\phi (u_{a,r})\|_\Psi =
o\,\Big(\frac{1}{\Psi^{-1}\big(\frac{1}{1-r}\big)}\Big) \hskip 3mm \text{as}\ r\to 1 
\tag*{\ref{petit o} }
\end{equation}
implies that $m(|\phi|=1)=0$.
\end{proposition}

\noindent{\bf Proof.} Otherwise, one has $m(|\phi|=1)=\delta>0$. Splitting the unit circle $\T$ into 
$N$ parts, we get some $a\in \T$ such that:
\begin{displaymath}
m(|\phi - a|\leq \pi/N)\geq \delta/N.
\end{displaymath}
But, the inequality $|\phi - a |\leq \pi/N$ implies, with $r=1 - 1/N$ (since $|\phi|\leq 1$):
\begin{displaymath}
|1 - \bar{a} r\phi| \leq |1 - \bar{a} \phi| + (1-r)| \bar{a} \phi| = |a -\phi| + (1-r) |\phi| < 5/N = 5(1-r).
\end{displaymath}
Hence:
\begin{align*}
m\big(|C_\phi (u_{a,r})|> 1/25) \big) & = m\big(|1 - \bar{a} r\phi|< 5(1-r)\big) \\
& \geq m\big(| \phi - a| \leq  \pi(1-r)\big) \geq \delta (1-r) \,, 
\end{align*}
and therefore, by Lemma \ref{Markov}:
\begin{displaymath}
\| C_\phi(u_{a,r}) \|_\Psi \geq \frac{1/25}{\Psi^{-1}\big(1/\delta (1-r)\big)}
\geq \frac{\delta/25}{\Psi^{-1}\big(1/(1-r)\big)}\,\cdot
\end{displaymath}
Since $r$ can be taken arbitrarily close to $1$, that proves the proposition. \hfill$\square$
\bigskip

\noindent{\bf Proof of Theorem \ref{Delta2}.} Assume that condition~\ref{petit o} is 
satisfied, and fix $A>1$. Let $\eps >0$ to be adjusted later. We can find $r_\eps < 1$ such that 
$r_\eps \leq r < 1$ implies:
\begin{displaymath}
\| C_\phi( u_{a,r}) \|_\Psi \leq \frac{\eps}{\Psi^{-1} \big(\frac{1}{1-r}\big)}\,
\raise 0,5mm \hbox{,}\quad \forall a\in \T.
\end{displaymath}
Now, Lemma~\ref{Markov} also reads:
\begin{displaymath}
m ( |f| > t ) \leq \frac{1}{\Psi\big( \frac{t}{ \|f\|_\Psi}\big)}\,\raise 0,5mm\hbox{,}
\end{displaymath}
so that, if one sets $B= 1/9\eps$:
\begin{displaymath}
m \big( | 1 - \bar{a} r \phi| < 3 (1-r) \big) = m \big(|C_\phi (u_{a,r}) | > 1/9\big)  \leq 
\frac{1} {\Psi\big[B \Psi^{-1} \big( \frac{1}{1-r}\big)\big]}\,\cdot
\end{displaymath}
\par
We claim that this implies a good upper bound on $m ( |\phi | > r)$, even if we loose a factor $1/(1-r)$, 
due to the effect of a rotation on $\phi$. For that, we shall use the following lemma.
\medskip

\begin{lemme}\label{rotation} 
Let $\phi \colon \D \to \D$ be an analytic function. Then, for every $r$ with $0 < r <1$, 
there exists $a\in \T$ such that:
\begin{displaymath}
m\big(|1 - \bar{a} r\, \phi| < 3(1-r)\big) \geq \frac{1-r}{8}\, m(|\phi| > r).
\end{displaymath}
\end{lemme}
\medskip

Admitting this for a while, we are going to finish the proof. \par
Fix an $r$ such that $ r_\eps \leq r < 1$, and take an $a\in \T$ as in Lemma~\ref{rotation}. We get, 
from the preceding, in setting $\lambda = 1-r$:
\begin{align*}
m ( 1 -|\phi| < \lambda ) 
&= m ( |\phi | > r) \\
& \leq \frac{8}{1 -r} m \big( | 1- \bar{a} r \phi | < 3 (1 - r)\big) \\
& \leq \frac{8}{1 -r} \frac{1}{\Psi\big[ B \Psi^{-1} \big(\frac{1}{1-r}\big)\big]}\, 
\raise 0,5mm \hbox{,}
\end{align*}
{\it i.e.}, setting $x = \Psi^{-1} (1/1-r) = \Psi^{-1} (1/\lambda)$:
\begin{displaymath}
m ( 1 -|\phi| < \lambda ) \leq 8\, \frac{\Psi (x)}{ \Psi (Bx)}\,\cdot
\end{displaymath}
But $\Psi$ satifies the $\Delta^2$-condition: $\big[\Psi(y)\big]^2\leq \Psi(\alpha y)$ for 
some $\alpha >1$ and $y$ large enough. Then, adjusting now $\eps >0$ as $\eps = 1/ 9\alpha A$, 
in order that $B =\alpha A$, we get, for $x$ large enough, since $A>1$:
\begin{displaymath}
\Psi (x) \Psi (A x) \leq [\Psi (Ax)]^2 \leq \Psi (Bx).
\end{displaymath}
Therefore, for $r$ close enough to $1$:
\begin{displaymath}
m ( 1 -|\phi| < \lambda ) \leq \frac{8}{\Psi (A x)} = \frac{8}{\chi_A (1/\lambda)}\,\cdot
\end{displaymath}
We hence get condition~\ref{grand O}, and that proves Theorem~\ref{Delta2}.\hfill $\square$
\bigskip

\noindent{\bf Proof of Lemma \ref{rotation}.} Let $\lambda=1-r$, and let $\delta>0$ be a 
number which we shall specify later. Consider the set:
\begin{displaymath}
C_\delta=\{z\in \D\,;\ |z|\geq 1-\lambda \hskip 2mm \text{and} 
\hskip 2mm |{\rm arg}\,z| \leq \delta\}
\end{displaymath}
(for $\delta=\lambda$, $C_\delta$ is a closed Carleson window). It is geometrically clear that 
$(1-\lambda) C_\delta$ is contained in the closed disk of center $1$ and whose edge contains 
$(1-\lambda)^2\e^{i\delta}$; hence, for every $z\in C_\delta$, one has:
\begin{align*}
|1 - (1-\lambda) z|^2 & \leq |1 - (1-\lambda)^2 \e^{i\delta} |^2
= 2(1-\lambda)^2(1-\cos \delta) + \lambda^2 (2 - \lambda)^2 \\
& \leq (1 - \lambda)^2 \delta^2 + \lambda^2 (2 - \lambda)^2 \leq 9\lambda^2
\end{align*}
if $\delta \leq \lambda$.\par
By rotation, one has, for every $a\in \T$:
\begin{displaymath}
|z|\geq 1-\lambda \hskip 3mm \text{and }\hskip 3mm  |{\rm arg}\,(\bar{a} z)| \leq \delta \hskip 3mm
\Rightarrow \hskip 3mm |1 - (1-\lambda) \bar{a} z|\leq 3\lambda.
\end{displaymath}
Let now $N\geq 2$ be the integer such that:
\begin{displaymath}
\frac{\pi}{N} \leq \lambda < \frac{\pi}{N-1}\,\raise 0,5mm\hbox{,}
\end{displaymath}
and take $\delta = \pi/N$.\par
One has, by the previous inequalities, setting $a_k=\e^{2ik\delta}$:
\begin{displaymath}
\{z\in \D\,;\ |z|\geq 1-\lambda \} = \bigcup_{1 \leq k\leq N} \bar{a}_k C_\delta
\subseteq \bigcup_{1\leq k\leq N}
\{z\in\D\,;\ |1 - (1-\lambda) \bar{a}_k z| \leq 3\lambda\}.
\end{displaymath}
Hence, with $z=\phi(\e^{it})$ (remark that, by Proposition \ref{mesure nulle}, we have only 
to consider the values of $\e^{it}$ for which $|\phi(\e^{it})|<1$; however, in this lemma, 
we may replace $\D$ by $\overline{\D}$), and get:
\begin{displaymath}
m(|\phi|\geq 1-\lambda) \leq N \sup_{1\leq k \leq N}
m(|1 - (1-\lambda) \bar{a}_k \phi| \leq 3\lambda).
\end{displaymath}
Therefore, we can find some $k$ such that:
\begin{displaymath}
m(|1 - (1-\lambda) \bar{a}_k \phi| \leq 3\lambda)
\geq \frac{1}{N} m(|\phi|\geq 1-\lambda) 
\geq \frac{\lambda}{8} m(|\phi|\geq 1-\lambda),
\end{displaymath}
since $\lambda \leq 2\pi/N\leq 8/N$. That proves Lemma \ref{rotation}.\hfill$\square$
\par\bigskip

Since the $\Delta^2$-condition implies the $\Delta^1$-condition, which, in its turn, implies  
the $\Delta^0$-con\-dition, we get, from Theorem \ref{impli-order bounded}, 
Theorem \ref{faiblement compact} and Theorem \ref{Delta2} that the
weak compactness of $C_\phi$ implies its order boundedness into $M^\Psi(\T)$, and thanks to the Remark after Proposition \ref{equiv order bounded}, its compactness. We get

\begin{theoreme}\label{corollaire Delta2}
If $\Psi$ satisfies the $\Delta^2$-condition, then the following assertions for 
composition operator $C_\phi \colon H^\Psi \to H^\Psi$ are equivalent:
\begin{itemize}
\item[$1)$] $C_\phi$ is order bounded into $M^\Psi(\T)$;
\item[$2)$] $C_\phi$ is compact;
\item[$3)$] $C_\phi$ is weakly compact;
\item[$4)$] $\Psi^{-1}\big(\frac{1}{1-|\phi|}\big)\in M^\Psi(\T)$ ({\it i.e.}: 
$ \chi_B (1/ 1 -|\phi|) \in L^1 (\T)$ for every $B>0$);
\item[$5)$] $m(1-|\phi| < \lambda )= O\,\Big(\frac{1}{\chi_A(1/\lambda)}\Big)$ as 
$\lambda\to 0$, for every $A>0$;
\item[$6)$] $\sup_{a\in \T} \|C_\phi (u_{a,r})\|_\Psi
=o\,\Big(\frac{1}{\Psi^{-1}\big(\frac{1}{1-r}\big)}\Big)$ as $r\to 1$\hskip15pt $(W)$
\end{itemize}
\end{theoreme}

\noindent{\bf Remark.} We shall see in the next section (Theorem \ref{suite exemple}) that the 
assumption $\Psi\in \Delta^2$ cannot be removed in general: Theorem \ref{corollaire Delta2} is not 
true for the Orlicz function $\Psi(x)= \exp\big[\big((\log (x+1)\big)^2\big] - 1$ (which 
nevertheless satisfies $\Delta^1$).
\medskip

If one specializes this corollary to the case where $\Psi(x)=\Psi_2(x)=\e^{x^2}-1$, 
which verifies the $\Delta^2$-condition, we get, using Stirling's formula:
\goodbreak

\begin{corollaire}
The following assertions are equivalent:
\begin{itemize}
\item[$1)$] $C_\phi\colon H^{\Psi_2} \to H^{\Psi_2}$ is order bounded into $M^{\Psi_2}(\T)$;
\item[$2)$] $C_\phi\colon H^{\Psi_2} \to H^{\Psi_2}$ is compact;
\item[$3)$] $C_\phi\colon H^{\Psi_2} \to H^{\Psi_2}$ is weakly compact;
\item[$4)$] $\frac{1}{1-|\phi|}\in L^p(\T), \forall p\geq 1$;
\item[$5)$] $\forall q\geq 1$ $\exists C_q>0$:\quad $m(1-|\phi| < \lambda)\leq C_q \lambda^q$;
\item[$6)$] $\forall q\geq 1$ $\|\phi^n\|_1 = o\,(n^{-q})$;
\item[$7)$] $\|\phi^n\|_{\Psi_2} = o\,(1/\sqrt{\log n})$.
\end{itemize}
\end{corollaire}

As a consequence of Theorem \ref{corollaire Delta2}, we obtain the following:

\begin{corollaire}\label{exemple non compact}
Assume that $\Psi \in \Delta^2$. Then there exist compact composition operators 
$C_\phi \colon H^p \to H^p$ for $1\leq p <\infty$ which are not compact as operators 
$C_\phi \colon H^\Psi \to H^\Psi$.
\end{corollaire}

\noindent{\bf Remark.} We shall see in Theorem~\ref{implique Carleson} that compactness on $H^\Psi$ 
implies compactness on $H^p$ for $p < \infty$. Note that this shows that, though $H^\Psi$ is an interpolation 
space between $H^1$ and $H^\infty$ (see \cite{Be-Sh}, Theorem V.10.8), the compactness of 
$C_\phi \colon H^1 \to H^1$ with the continuity of $C_\phi\colon H^\infty \to H^\infty$ does not suffice to 
have compactness for $C_\phi \colon H^\Psi \to H^\Psi$.
\medskip

\noindent{\bf Proof.} 
As the $\Delta^2$ condition implies the $\Delta^0$ condition, we have $x = o\,\big(\chi_\beta (x)\big)$ 
as $x \to \infty$, for some $\beta >1$. It follows that we can find a positive 
function $a \colon \T \to \R_+$ such that $a\geq 2$, $a \in L^1$ but 
$\chi_\beta (a) \notin L^1$. Set $h = 1 - \frac{1}{a}$. One has $1/2 \leq h \leq 1$ and in particular 
$\log h \in L^1(\T)$. Then the outer function $\phi \colon \D \to \C$ defined for $z\in\D$ by:
\begin{displaymath}
\phi(z)= \exp\bigg[\int_\T \frac{u +z}{u-z}\,\log h(u)\,dm(u)\bigg]
\end{displaymath}
is analytic on $\D$ and its boundary limit verifies $|\phi| = h \leq 1$ on $\T$. 
By \cite{Shap-T}, Theorem 6.2,  $C_\phi \colon H^2 \to H^2$ is Hilbert-Schmidt, and hence compact, 
since $\int_\T \frac{1}{1 - |\phi|}\,dm = \int_\T a\,dm <+\infty$. It is then compact from $H^p$ to 
$H^p$ for every $p<\infty$ (\cite{Shap-T}, Theorem 6.1). However, 
$\int_\T \chi_\beta\big(\frac{1}{1 - |\phi|}\big)\,dm = \int_\T \chi_\beta(a)\,dm =+\infty$, and hence, 
by our Theorem \ref{corollaire Delta2}, $C_\phi$ is not compact on $H^\Psi$.\hfill$\square$

\bigskip\goodbreak

\subsection{$\boldsymbol{ p}$-summing operators.}

Recall that an operator $T\colon X \to Y$ between two Banach spaces is said to be 
$p$-summing ($1\leq p<+\infty$) if there is a constant $C>0$ such that, for every choice 
of $x_1,\ldots, x_n\in X$, one has:
\begin{displaymath}
\sum_{k=1}^n \| Tx_k\|^p \leq C^p \sup_{\stackrel{x^\ast\in X^\ast}{\scriptscriptstyle \|x^\ast\|\leq 1}}
\sum_{k=1}^n | x^\ast(x_k)|^p.
\end{displaymath}
In other terms, $T$ maps weakly unconditionaly $p$-summable sequences into norm 
$p$-summable sequences. When $X\subseteq Y=L^\Psi$, this implies that whenever 
$\sum_{n\geq 1} |g_n| \in L^\Psi$, then $\sum_{n\geq 1} \|Tg_n\|_\Psi^p <+\infty$.
\medskip

For $1\leq p <+\infty$, J. H. Shapiro and P. D. Taylor proved in \cite{Shap-T}, 
Theorem 6.2, that the condition:
\begin{equation}\label{ST}
\int_\T \frac{dm}{1 -|\phi|} <+\infty
\end{equation}
implies that the composition operator $C_\phi\colon H^p \to H^p$ is $p$-summing (and 
condition (\ref{ST}) is necessary  for $1\leq p\leq 2$; in particular, for $p=2$, it is equivalent to 
say that $C_\phi$ is Hilbert-Schmidt). Actually, they proved that (\ref{ST}) is equivalent to 
the fact that $C_\phi$ is order bounded on $H^p$, and 
(acknowledging to A. Shields, L. Wallen, and J. Williams) every order bounded operator into an $L^p$-space 
is $p$-summing. The counterpart of (\ref{ST}) in our setting, are conditions \ref{integrabilite-1} and 
\ref{integrabilite}
\begin{displaymath}
\int_\T \chi_A\Big(\frac{1}{1-|\phi|}\Big)\,dm <+\infty
\end{displaymath}
in Proposition~\ref{equiv order bounded}. We are going to see that, if $\Psi$ grows fast 
enough, order boundedeness does not imply that $C_\phi$ is $p$-summing. Note that, for 
composition operators on $H^\infty$, being $p$-summing is equivalent to being compact 
(\cite{Pascal}, Theorem 2.6), but $H^\infty$ corresponds to the very degenerate Orlicz 
function $\Psi(x)=0$ for $0\leq x \leq 1$ and $\Psi(x)=+\infty$ for $x>1$, which does not 
match in the proof below. 

\begin{theoreme}\label{exemple}
If $\Psi\in \Delta^2$, then there exists a composition operator 
$C_\phi\colon H^\Psi \to H^\Psi$ which is 
order bounded into $M^\Psi(\T)$, and hence compact, but which is $p$-summing for no 
$p\geq 1$.
\end{theoreme}

Note that every $p$-summing operator is Dunford-Pettis (it maps the weakly convergent sequences into norm 
convergent sequences); therefore, when it starts from a reflexive space, it is compact. However, when 
$\Psi\in \Delta^2$, being Dunford-Pettis implies compactness for composition operators on $H^\Psi$, 
though $H^\Psi$ is not reflexive, thanks to the next proposition and Theorem \ref{corollaire Delta2}. Later (see Theorem \ref{equiv weak-comp}), we shall see that, under condition $\Delta^0$, every Dunford-Pettis composition operator is compact.

\begin{proposition}\label{DP}
When $\Psi\in\nabla_2$, every Dunford-Pettis composition operator satisfies condition $(W)$.
\end{proposition}

\noindent{\bf Proof.} Let $g_{a,r} = \Psi^{-1}\big(1/(1-r)\big) u_{a,r}$. If condition \ref{petit o} were not satisfied, , we could find a sequence $(a_n)_{n\geq 1}$ in $\T$ and a sequence of numbers $(r_n)_{n\geq 1}$ tending to $1$ such 
that  $\| C_\phi (g_{a_n,r_n})\|_\Psi \geq \delta >0$ for all $n \geq 1$. But 
$(1-r)^2 \Psi^{-1}\big(1/(1-r)\big)\mathop{\longrightarrow}\limits_{r\to 1} 0$.
Therefore $g_{a_n,r_n} (z) = (1-r_n)^2 \Psi^{-1}\big(1/(1-r_n)\big) / (1 - \bar{a}_n r_n z) $ tends to $0$ 
uniformly on compact sets of $\D$. Hence, by Proposition~\ref{weak-star topology}, 
$(g_{a_n, r_n})_{n\geq 1}$ tends weakly to $0$ (because $g_{a_n,r_n} \in HM^\Psi$ and, on $HM^\Psi$, 
the weak-star topology of $H^\Psi$ is the weak topology). Since $C_\phi$ is Dunford-Pettis, 
$\big(C_\phi (g_{a_n,r_n})\big)_{n\geq 1}$ tends in norm to $0$, and we get a contradiction, proving 
the proposition.\hfill $\square$
\medskip 

\noindent{\bf Proof of Theorem \ref{exemple}.} 
We shall begin with some preliminaries. First, since $\Psi\in \Delta^2$, there exists 
$\alpha >1$ such that $\big[\Psi(x)\big]^2 \leq \Psi(\alpha x)$ for $x$ large enough. Hence:
\begin{displaymath}
\frac{\Psi (x)}{\Psi (x/\alpha)} \geq \Psi \Big( \frac{x}{\alpha}\Big) 
\mathop{\longrightarrow}_{x\to +\infty} +\infty.
\end{displaymath} 
Therefore, there exists, for every $n\geq 1$, some $x_n >0$ such that:
\begin{displaymath}
\frac{\Psi (x)}{\Psi (x/\alpha)} \geq 2^n\hskip 5mm \forall x\geq x_n.
\end{displaymath}
Then
\begin{displaymath}
\Psi\Big(\frac{x}{\alpha}\Big) \leq \frac{1}{2^n} \Psi(x) +\Psi \Big(\frac{x_n}{\alpha}\Big) 
\leq \frac{1}{2^n} \Psi(x) +\Psi (x_n)  \hskip 3mm \forall x>0.
\end{displaymath}
For convenience, we shall assume, as we may, that $\Psi(x_n)\geq 1$.\par\noindent
Remark also that, setting $a=\Psi^{-1}(1)$, one has, for every $f\in L^\infty$:
\begin{displaymath}
\int_\T \Psi\Big(a\,\frac{|f|\,\,}{\|f\|_\infty}\Big)\,dm \leq 1,
\end{displaymath}
so that:
\begin{displaymath}
\| f\|_\Psi\leq \frac{1}{a}\|f\|_\infty.
\end{displaymath}
\medskip
We are now going to start the construction.\par
For $n\geq 1$, let $M_n=\log (n+1)$. Choose positive numbers $\beta_n$ which tend to 
$0$ fast enough to have: 
\begin{displaymath}
\sum_{k>n} \beta_k \leq \beta_n\,,\hskip 3mm  \forall n\geq 1,
\end{displaymath}
and
\begin{displaymath}
t_n= \frac{\Psi^{-1}(8/\beta_n)}{M_n} \mathop{\longrightarrow}_{n\to +\infty} +\infty.
\end{displaymath}
Set:
\begin{displaymath}
r_n= 1- \frac{1}{\Psi(t_n)}\,\cdot
\end{displaymath}
One has $r_n\mathop{\longrightarrow}\limits_{n\to +\infty} 1$ and:
\begin{displaymath}
\chi_{M_n}\Big(\frac{1}{1-r_n}\Big) = \frac{8}{\beta_n}\,\cdot
\end{displaymath}
\smallskip
Actually, for the end of the proof, we shall have to choose the $\beta_n$'s 
decreasing so fast that:
\begin{displaymath} 
\bigg[\Psi \Big(\frac{t_1 +\cdots + t_{n-1}}{\alpha}\Big) +\Psi (x_n)\bigg]
\frac{2^n t_n}{\Psi (t_n)} \leq \frac{1}{2^n} \,\cdot
\end{displaymath}
This is possible, by induction, since 
$t/\Psi (t)\mathop{\longrightarrow}\limits_{t\to +\infty} 0$. Note that, since $\Psi(x_n)\geq 1$, 
one has, in particular: 
$\sum_{n=1}^{+\infty} \frac{2^n t_n}{\Psi(t_n)} < +\infty$
\medskip

Let $B_n$ be disjoint measurable subsets of $\T$ with measure $m(B_n)=c \beta_n$ 
(where $c \geq 1$ is such that $\sum_{n\geq 1} \beta_n=1/c$), and whose union is $\T$. 
Define $h \colon \T \to \C$ by:
\begin{displaymath}
h =\sum_{n\geq 1} r_n \ind_{B_n}.
\end{displaymath}
One has $\log h \in L^1(\T)$, since $h$ does not vanish, $r_n\geq 1/2$ for $n$ large enough, 
and $\sum_n m(B_n) = 1 <+\infty$. We can define the outer function:
\begin{displaymath}
\phi(z)= \exp\bigg[\int_\T \frac{u +z}{u-z}\,\log h(u)\,dm(u)\bigg]\,,\hskip 3mm |z|<1.
\end{displaymath}
$\phi$ is analytic on $\D$ and its boundary limit verifies $|\phi| = h \leq 1$ on $\T$. 
Hence $\phi$ defines a composition operator on $H^\Psi$.\par
For any $A>0$, one has, when $n$ is large enough to ensure $M_n\geq A$, and when 
$r_n \leq r <r_{n+1}$:
\begin{align*}
m(|\phi| > r) &= \sum_{k>n} m(|\phi|=r_k) =\sum_{k>n} c\,\beta_k \leq c\,\beta_n 
= \frac{8\,c}{\chi_{M_n} \big(1/(1-r_n)\big)} \\
& \leq \frac{8\,c}{\chi_A \big(1/(1-r_n)\big)} 
\leq \frac{8\,c}{\chi_A \big(1/(1-r)\big)}\,\cdot 
\end{align*}
Since $r_n\mathop{\longrightarrow}\limits_{n\to +\infty} 1$, it follows from 
Theorem~\ref{impli-order bounded} that $C_\phi \colon H^\Psi \to H^\Psi$ is order-bounded 
into $M^\Psi(\T)$ (and hence is compact).\par\smallskip

We are now going to construct a sequence of functions $g_n\in H^\Psi$ such that 
$\sum_n |g_n| \in L^\Psi$, but $\sum_n \|C_\phi(g_n)\|_\Psi^p =+\infty$ for all $p\geq 1$. 
That will prove that $C_\phi$ is $p$-summing for no $p\geq 1$.\par\smallskip

Since
\begin{displaymath}
m(|\phi| \geq r_n) \geq m(|\phi | = r_n) = c\,\beta_n \geq \beta_n, 
\end{displaymath}
we can apply Lemma \ref{Markov} and Lemma \ref{rotation} (which remain valid with non-strict inequalities 
instead of strict ones), and we are able to find, for every $n\geq 1$, some $a_n\in \T$ such that:
\begin{displaymath}
\|C_\phi(u_{a_n,r_n})\|_\Psi \geq 
\frac{1/9}{\displaystyle \Psi^{-1} \Big(\frac{8}{(1 - r_n) \beta_n}\Big)}\,\cdot
\end{displaymath}
But:
\begin{displaymath}
\frac{8}{(1 - r_n) \beta_n} = \Psi (t_n) \Psi (M_n t_n).
\end{displaymath}
Since now $\Psi$ satisfies $\Delta^2$: $[\Psi(x)]^2\leq \Psi(\alpha x)$, for $x$ large enough and one has, for $n$ large enough, 
since $M_n \geq 1$ (for $n \geq 2$):
\begin{displaymath}
\Psi (t_n) \Psi (M_n t_n) \leq \big[\Psi (M_n t_n)\big]^2  \leq \Psi (\alpha M_n t_n).
\end{displaymath}
Therefore:
\begin{displaymath}
\|C_\phi(u_{a_n,r_n})\|_\Psi \geq \frac{1/9}{\alpha M_n t_n}\,\cdot
\end{displaymath}
Taking now:
\begin{displaymath}
g_n = \Psi^{-1}\Big(\frac{1}{1 - r_n}\Big) u_{a_n, r_n}= t_n u_{a_n, r_n} \,,
\end{displaymath}
one has $\|g_n\|_\Psi \leq 1$ (by Corollary \ref{norme $u_{a,r}$}), and 
\begin{displaymath}
\| C_\phi (g_n)\|_\Psi \geq \frac{1/9}{\alpha M_n} = \frac{1/9}{\alpha \log (n+1)}\,\cdot
\end{displaymath} 
Therefore
\begin{displaymath}
\sum_{n=1}^{+\infty} \| C_\phi(g_n)\|_\Psi^p = +\infty
\end{displaymath}
for every $p\geq 1$.
\par\medskip

It remains to show that $g= \sum_n |g_n| \in L^\Psi$.\par
We shall follow the lines of proof of Theorem II.1 in \cite{LLQR-Studia}.\par
By Markov's inequality, one has:
\begin{displaymath}
m (|g_n| > 2^{-n}) \leq 2^n t_n \| u_{a_n, r_n}\|_1 \leq 
2^n t_n (1-r_n) = \frac{2^n t_n}{\Psi(t_n)}\,\cdot
\end{displaymath}
Set:
\begin{displaymath}
A_n= \{ |g_n| > 2^{-n}\} \hskip 1mm ; \hskip 5mm 
\tilde A_n =A_n \setminus \bigcup_{j>n} A_j,
\end{displaymath}
and
\begin{displaymath}
\breve g_n = g_n \ind_{\{ |g_n| > 2^{-n}\}}.
\end{displaymath}
Since: 
\begin{displaymath}
\sum_{n=1}^{+\infty} \| g_n - \breve g_n \|_\Psi \leq 
\frac{1}{a} \sum_{n=1}^{+\infty} \| g_n - \breve g_n \|_\infty \leq 
\frac{1}{a} \sum_{n=1}^{+\infty} \frac{1}{2^n} = \frac{1}{a} < +\infty,
\end{displaymath}
it suffices to show that $\breve g = \sum_n |\breve g_n | \in L^\Psi$.\par
But $\breve g$ vanishes out of 
$\bigcup_{n\geq 1} \tilde A_n \cup\big(\limsup_n A_n\big)$, and 
$m\big(\limsup_n A_n\big)=0$, since
\begin{displaymath}
\sum_{n=1}^{+\infty} m(A_n) \leq \sum_{n=1}^{+\infty} \frac{2^n t_n}{\Psi(t_n)} < +\infty\,. 
\end{displaymath}
Therefore:
\begin{displaymath}
\int_\T \Psi\Big(\frac{| \breve g|}{2\alpha}\Big)\,dm 
= \sum_{n=1}^{+\infty} \int_{\tilde A_n} \Psi\Big(\frac{| \breve g|}{2\alpha}\Big)\,dm.
\end{displaymath}
Since $\breve g_j=0$ on $\tilde A_n$ for $j>n$, we get:
\begin{displaymath}
\int_\T \Psi\Big(\frac{| \breve g|}{2\alpha}\Big)\,dm 
= \sum_{n=1}^{+\infty} \int_{\tilde A_n} 
\Psi\Big(\frac{| \breve g_1| + \cdots + | \breve g_n|}{2\alpha}\Big)\,dm.
\end{displaymath}
Now, by the convexity of $\Psi$:
\begin{displaymath}
\Psi\Big(\frac{| \breve g_1| + \cdots + | \breve g_n|}{2\alpha}\Big) 
\leq \frac{1}{2} \bigg[ \Psi\Big(\frac{| \breve g_1| + \cdots + | \breve g_{n-1}|}{\alpha}\Big) 
+ \Psi\Big(\frac{| \breve g_n |}{\alpha}\Big)\bigg]\,\cdot
\end{displaymath}
But:
\begin{displaymath}
\Psi\Big(\frac{| \breve g_n |}{\alpha}\Big) 
\leq \frac{1}{2^n} \Psi ( |\breve g_n |) + \Psi(x_n)\,,
\end{displaymath}
and
\begin{displaymath}
\Psi\Big(\frac{| \breve g_1| + \cdots + | \breve g_{n-1}|}{\alpha}\Big) 
\leq \Psi\Big(\frac{t_1 + \cdots + t_{n-1}}{\alpha}\Big)\,;
\end{displaymath}
therefore, using that $\int_\T \Psi(|\breve g_n|)\,dm \leq \int_\T \Psi(|g_n|)\,dm \leq 1$:
\begin{align*}
\int_\T \Psi\Big(\frac{| \breve g|}{2\alpha}\Big)\,dm 
& \leq \sum_{n=1}^{+\infty} \frac{1}{2} \bigg[
\Psi\Big(\frac{t_1 + \cdots + t_{n-1}}{\alpha}\Big) m(\tilde A_n) \\
& \hskip 2,5cm + \frac{1}{2^n} \int_\T \Psi(|\breve g_n|)\,dm + \Psi(x_n) m(\tilde A_n)\bigg] \\
& \leq \sum_{n=1}^{+\infty} \frac{1}{2} \bigg[
\Big[\Psi\Big(\frac{t_1 + \cdots + t_{n-1}}{\alpha}\Big) + \Psi(x_n)\Big]
\frac{2^n t_n}{\Psi(t_n)} + \frac{1}{2^n}\bigg] \\
& \leq \sum_{n=1}^{+\infty} \frac{1}{2} \Big[\frac{1}{2^n} + \frac{1}{2^n}\Big] =1\,,
\end{align*}
which proves that $\breve g \in L^\Psi$, and $\|\breve g \|_\Psi\leq 2\alpha$.\par
The proof is fully achieved.\hfill$\square$
\medskip

\noindent{\bf Remark.} In the above proof, we chose $M_n= \log (n+1)$. This choice was only 
used to conclude that $\sum_n \| C_\phi (g_n) \|_\Psi^p =+\infty$ for every $p< \infty$. Therefore, 
the above proof shows that, given any increasing function $\Upsilon \colon (0,\infty) \to (0,\infty)$ 
tending to $\infty$, we can find, with a suitable choice of a slowly increasing sequence $(M_n)_{n\geq 1}$, a symbol $\phi$ and a sequence $(g_n)_{n\geq 1}$ in $H^\Psi$ such that $\sum_n |g_n| \in L^\Psi$, although 
$\sum_n \Upsilon (\| C_\phi (g_n)\|_\Psi) =+\infty$.
\bigskip\goodbreak

%%%%%%%%%%%%%%%%%%%%%%%%%%%%%%%%%%%%%%%%%%%%%%%%%%%%%%%%%%%%%%%%%%%%%%%%%%%%%%%%%%%%%%%%%%%%%

\section{Carleson measures}

\subsection{Introduction}

B. MacCluer (\cite{McC}; see also \cite{Co-McC}, Theorem 3.12) has characterized compact 
composition operators on Hardy spaces $H^p$ ($p<\infty$) in term of Carleson measures. 
In this section, we shall give an analogue of this result for Hardy-Orlicz spaces $H^\Psi$, but in terms 
of ``$\Psi$-Carleson measures''. Indeed, Carleson measures do not characterize the compactness 
of composition operators when $\Psi$ grows too quickly, as it follows from 
Corollary \ref{exemple non compact}.\par
\medskip

Before that, we shall recall some definitions (see for example \cite{Co-McC}, pages 37--38, or 
\cite{Du}, page 157).\par
\smallskip

Let $\xi\in{\T}$ and $h\in(0,1)$. Define
\begin{equation}\label{def S}
S(\xi,h)=\{z\in\overline{\D}\, ; \ |\xi-z|<h \}.
\end{equation}
The \emph{Carleson window} $W(\xi,h)$ is the following subset of $\overline{\D}$:
\begin{equation}\label{def W}
W(\xi,h)=\{z \in\overline{\D}\,;\ 1-h < |z|\le1\hskip 15pt \hbox{and} \hskip 15pt |{\rm arg}(z\bar\xi)|<h\}.
\end{equation}
It is easy to show that we have for every $\xi\in \T$ and $h\in(0,1)$:
\begin{displaymath}
S(\xi,h/2) \subseteq W(\xi,h)\hskip1cm\hbox{and}\hskip1cm W(\xi,h/2)\subseteq S(\xi,h)\,,
\end{displaymath}
so that, in the sequel, we may work equivalently with either $S(\xi,h)$ or $W(\xi,h)$.
Recall that a positive Borel measure $\mu$ on $\D$ (or $\overline \D$) is called a 
\emph{Carleson measure} if there exists some constant $K>0$ such that:
\begin{displaymath}
\hskip 1,5cm \mu \big( S(\xi, h)\big) \leq K h\,,\hskip 5mm \forall \xi\in \T\,,\ \forall h\in (0,1).
\end{displaymath}
Carleson's Theorem (see \cite{Co-McC}, Theorem 2.33, or \cite{Du}, Theorem 9.3) asserts that, for 
$0< p < \infty$, the Hardy space  $H^p$ is continuously embedded into $L^p(\mu)$ if and 
only if $\mu$ is a Carleson measure.
\smallskip

Given an analytic self-map $\phi \colon \D \to \D$, we define the 
\emph{pullback measure} $\mu_\phi$ on the closed unit disk $\overline{\D}$ (which we shall
denote simply $\mu$ when this is unambiguous) as the image of the Haar measure $m$ of 
$\T = \partial \D$ under the map $\phi^\ast$ (the boundary limit of $\phi$):
\begin{equation}\label{pull-back}
\mu_\phi(E) = m \big({\phi^\ast}^{-1}(E)\big),
\end{equation}
for every Borel subset $E$ of $\overline{\D}$.\par
\medskip

The automatic continuity of composition operators $C_\phi$ on the Hardy space $H^p$, 
combined with Carleson's Theorem means that $\mu_\phi$ is always a Carleson measure.
\smallskip

B. MacCluer (\cite{McC}, \cite{Co-McC}, Theorem 3.12) showed that:
\begin{equation}\label{MC}
\begin{split}
& \textit{The composition operator $C_\phi$ is compact on $H^2$ if and only if:} \\
& \hskip 10mm \mu_\phi\big( S(\xi,h)\big) = o\,(h) \hskip 2mm \textit{as}\  h\rightarrow 0,
\ \textit{uniformly for}\ \xi\in \T. 
\end{split} 
\tag*{(MC)}
\end{equation}

While the Shapiro's compactness criterion, {\it via} the Nevanlinna counting function 
(\cite{Shap}), deals with the behavior of $\phi$ \emph{inside} the open unit disk, the characterization $(MC)$ deals with its boundary values $\phi^\ast$. It is natural to wonder whether the modulus of $\phi^\ast$ on $\T =\partial \D$ suffices to characterize the compactness of $C_\phi$. This leads to the following question: if two functions $\phi_1$ and $\phi_2$ have the same modulus on $\T$,  are the compactness of the two associated composition operators equivalent? We have seen 
in Theorem \ref{corollaire Delta2} that the answer is positive on $H^\Psi$ when $\Psi \in \Delta^2$. However, on $H^2$ it turns out to be  negative. We give the following counterexample.

\begin{theoreme}\label{meme module}
There exist two analytic functions $\phi_1$ and $\phi_2$ from $\D$ into itself  such that 
$|\phi_1^\ast | = |\phi_2^\ast|$ on $\T$ but for which the composition operator  
$C_{\phi_2}\colon H^2 \to H^2$ is compact, though $C_{\phi_1}\colon H^2 \to H^2$ is not compact.
\end{theoreme}

\noindent{\bf Remark.} Let $\Psi$ be an Orlicz function which satisfies $\Delta^2$. We shall see in 
Theorem \ref{implique Carleson} that every composition operator $C_\phi\colon H^2 \to H^2$ is 
compact as soon as $C_\phi\colon H^\Psi \to H^\Psi$ is compact. Hence, in the above theorem, 
$C_{\phi_1}\colon H^\Psi \to H^\Psi$ is not compact. It follows hence from 
Theorem \ref{corollaire Delta2}, since $\phi_1^\ast$ and $\phi_2^\ast$  have the same modulus, that 
$C_{\phi_2}\colon H^\Psi \to H^\Psi$ is not compact (and, even, not weakly compact), though 
$C_{\phi_2}\colon H^2 \to H^2$ is compact. We have already seen such a phenomenon in 
Corollary \ref{exemple non compact}. However, the results of the next subsection will allow us to 
conclude (Theorem \ref{suite exemple}) that, when $\Psi(x)= \exp\big[\big(\log(x+1)^2\big]- 1$, 
which does not satisfy condition $\Delta^2$, but satisfies conditions $\Delta^1$ and  $\nabla_1$, the 
composition operator $C_{\phi_2}\colon H^\Psi \to H^\Psi$ \emph{is} compact, but not order bounded 
into $M^\Psi(\T)$. That will show that our assumption that $\Psi\in \Delta^2$ in 
Theorem \ref{corollaire Delta2} is not only a technical one.
\medskip

\noindent{\bf Proof.} We start simply with $\phi_1 (z)=\frac{1+z}{2}\,$. It is well 
known that $C_{\phi_1}$ is not compact on $H^2$ (this was first observed in 
H. J. Schwartz's thesis~\cite{Schw}: see \cite{Shap-T}, page 471). Now, let:
\begin{displaymath}
M(z)=\exp \Big(- \frac{1 + z}{1 - z}\Big)
\end{displaymath}
and
\begin{displaymath}
\phi_2 (z) = \phi_1(z)\,M(z). 
\end{displaymath}
For simplicity, we shall write $\phi=\phi_2$, and we are going to show that $C_\phi$ is a 
compact operator on $H^2$, using the criterion \ref{MC}.\par 
Let $\xi=\e^{i\alpha}\in \T$, with $|\alpha|\leq\pi$. We are going to prove that:
\begin{displaymath}
\mu_\phi \big(S(\xi,h)\big)=O\,(h^{3/2}).
\end{displaymath}
First, observe that, for $h\in(0,1)$:
\begin{displaymath}
S(\xi,h) \subseteq 
\{z\in\overline{\D}\, ;\ 1-h < |z|\le1 \hbox{ and }  |\arg (\bar z\xi)| \leq 2h \}.
\end{displaymath}
Hence for $h$ small enough,
\begin{align*}
\mu_\phi \big(S(\xi, & h) \big) 
\leq \\ 
& m\big( \{\theta\in (-\pi,\pi) \,;\ 1-h < |\phi (\e^{i\theta})|\le1 \hbox{ and } 
|\arg (\e^{-i\alpha}\phi(\e^{i\theta})) |\leq 2h \}).
\end{align*}
For $\theta\in(-\pi,\pi)$, one has  
$|\phi (\e^{i\theta})| =|\phi_1 (\e^{i\theta}) |=\cos(\theta/2)$, and so the condition  
$1-h < |\phi (\e^{i\theta})|\le1$ is equivalent to $ 1-h < \cos(\theta/2) <1$, which implies, since 
$\cos t = 1 - 2\sin^2 (t/2) \leq 1 -  2 t^2/\pi^2 \leq 1 - t^2/5$ for $0 \leq t \leq \pi/2$  
(because $\sin t \geq \frac{2}{\pi} t$), that $1-h < 1 - (\theta/2)^2/5$, {\it i.e.} 
$\theta^2 \leq 20 h$ and so  $|\theta | \leq 6 \sqrt h$. 
On the other hand, $M(\e^{i\theta})=\exp\big(-i \cot (\theta/2)\big)$; hence 
$\arg \phi (\e^{i\theta}) = \theta/2 - \cot (\theta/2)$, modulo        
$2\pi$. Therefore, for $h$ small enough:
\begin{align*}
\mu_\phi \big(S(\xi,h)\big) 
&\leq m (\{|\theta|\leq 6 \sqrt{h} \,;\ |-\alpha+\theta/2  - \cot (\theta/2)|\leq 2h, \bmod\,{2\pi} \}) \\
& \leq2 \sum_{n\in\Z} m( \{ |t|\leq 3 \sqrt h\,;\ |-\alpha + t - \cot t + 2\pi n|\leq 2h\}).
\end{align*}
We have to majorize both $\displaystyle\sum_{n\in\Z}m( \{0<t\leq 3 \sqrt h\,;\ |-\alpha + t - \cot t + 2\pi n|\leq 2h\})$ and $\displaystyle\sum_{n\in\Z}m( \{0<t\leq 3 \sqrt h\,;\ |\alpha + t - \cot t + 2\pi n|\leq 2h\})$.

The function 
\begin{displaymath}
F(t) = t - \cot (t)
\end{displaymath}
is increasing , and we define $a_n, b_n \in (0,\pi)$ by:
\begin{displaymath}
F(a_n)=\alpha - 2\pi  (n_h +n) - 2h \quad \text{and} \quad F(b_n)=\alpha - 2\pi (n_h +n) + 2h,
\end{displaymath} 
where the integer $n_h$ is given large enough to ensure that $a_0 \leq 3\sqrt h$. Of course, 
$a_n < b_n < a_{n-1}$. Observe that $4h=F(b_0)-F(a_0)\ge b_0-a_0$, and then $b_0\le 3\sqrt h + 4h \le 4\sqrt h$ for $h$ small enough. One has:
\begin{displaymath}
\displaystyle\sum_{n\in\Z} m( \{0<t\leq 3 \sqrt h\,;\ |-\alpha + t - \cot t + 2\pi n|\leq 2h\})\leq \sum_{n=0}^\infty (b_n-a_n).
\end{displaymath}
Since $\displaystyle F'(t) = 1 +\frac{1}{\sin^2 t} \geq \frac{1}{t^2}$, one has, on the one hand:
\begin{equation}\label{minoration de F}
4h =F(b_n) - F(a_n) = \int_{a_n}^{b_n} F'(t)\,dt \geq\frac{b_n - a_n}{a_n b_n}; 
\end{equation}
hence:
\begin{equation}\label{somme F}
b_n - a_n\le4h a_n b_n \le 4hb_n^2\quad,\,\hbox{for all }n\ge0. 
\end{equation}
On the other hand, let us first point out that, for $0\leq t \leq 1$,
\begin{displaymath}
\displaystyle F'(t) = 1 +\frac{1}{\sin^2 t} \leq 1 +\frac{\pi^2/4}{t^2} \leq \frac{4}{t^2}\cdot
\end{displaymath} 
Hence, for $h$ small enough: 
\begin{displaymath}
2\pi = F(b_n)- F(b_{n+1}) = \int_{b_{n+1}}^{b_n} F'(t)\,dt \leq 
4 \frac {b_n-b_{n+1}}{b_{n+1} b_n}\,\raise 0,5mm \hbox{,}
\end{displaymath}
and we get:
\begin{displaymath}
b_{n+1}^2 \leq b_{n+1} b_n \leq \frac{2}{\pi} (b_n - b_{n+1}).
\end{displaymath}
Hence, using the fact that \eqref{minoration de F} gives $b_0 - a_0 \leq 4h a_0b_0\leq 48\,h^2$, we get, 
from \eqref{somme F}:
\begin{align*}
\sum_{n=0}^\infty (b_n - a_n) 
& \leq (b_0 - a_0) + \frac{8 h}{\pi} \sum_{n=0}^\infty (b_n - b_{n+1}) \\ 
& \leq 48\,h^2 + \frac{8 h}{\pi} b_0 \leq 48\,h^2+ \frac{8 h}{\pi} 4 \sqrt{h} 
\leq \Big (48\sqrt{h} + \frac{32}{\pi} \Big) h^{3/2} \leq 11\, h^{3/2},
\end{align*}
for $h$ small enough.
\par
In the same way, we have 
\begin{displaymath}
\displaystyle\sum_{n=0}^\infty m( \{0< t\leq 3 \sqrt h\,;\ |\alpha + t - \cot t + 2\pi n|\leq 2h\})\le 11\, h^{3/2}.
\end{displaymath}

We can hence conclude that $\mu_\phi \big( S(\xi,h)\big) \leq C h^{3/2}$, where $C$ is a numerical 
constant.\hfill $\square$
\medskip

\noindent{\bf Remark.} $C_\phi$ actually maps continuously $H^2$ into $H^3$, and compactly 
$H^2$ into $H^p$, for any $p<3$ (see \cite{HJ} or Theorem \ref{continuite Carleson} and \ref{compacite Carleson}).
\medskip

However, in some cases, the behaviour of $|\phi^\ast|$ on the boundary $\partial \D$ suffices. 

\begin{proposition}
Let $\phi_1$ and $\phi_2$ be two analytic self-maps of $\D$ such that 
$|\phi_1^\ast| \leq |\phi_2^\ast |$ on $\partial \D$. Assume that they are both one-to-one on $\D$, and 
that there exists $a\in \D$ such that $\phi_2(a)=0$. Then the compactness of 
$C_{\phi_2} \colon H^2 \to H^2$ implies that of $C_{\phi_1} \colon H^2 \to H^2$.
\end{proposition}

\noindent{\bf Proof.} By composing $\phi_1$ and $\phi_2$ with the automorphism of $\D$ which maps $0$ into $a$, we 
may assume that $a=0$. We can hence write $\phi_2(z) = z \phi(z)$, with $\phi\colon \D \to \C$ analytic 
in $\D$. $\phi$ does not vanish in $\D$ because of the injectivity of $\phi_2$ (this is obvious for 
$z\not=0$, and for $z=0$, follows from the fact that the injectivity of $\phi_2$ implies $\phi_2'(0)\not= 0$).
\par 
Then there is some $\delta >0$ such that 
$|\phi(z)| \geq \delta$ for every $z\in \D$. In fact, by continuity, there is some $\alpha >0$ and some 
$0 < r <1$ such that $|\phi(z)| \geq \alpha$ and $|\phi_2(z)|\leq \alpha$ for $|z| \leq r$. But being 
analytic and non constant, $\phi_2$ is an open map, so there is some $\rho >0$ such that 
$\rho \D \subseteq \phi_2(r\overline{\D})$. Injectivity of $\phi_2$ shows that 
$\phi_2(\D \setminus r\overline{\D}) \cap \rho \D =\emptyset$, that is to say that 
$|\phi_2(z)| \geq \rho$ for $|z| >r$. {\it A fortiori} $|\phi (z)| \geq \rho$ for $|z| >r$. The claim is proved 
with $\delta=\min(\alpha, \rho)$.\par
Then $\frac{1}{\phi }\in H^\infty$, as well as $\frac{\phi_1}{\phi}$. Since 
$\big|\frac{\phi_1^\ast}{\phi^\ast}\big| = \big|\frac{\phi_1^\ast}{\phi_2^\ast}\big| \leq 1$, one has 
$\big|\frac{\phi_1(z)}{\phi(z)}\big| \leq 1$ for every $z \in \D$. Hence:
\begin{displaymath}
 \frac{ 1 - |\phi_1 (z)|}{ 1 - |z|} \geq  \frac{ 1 - |\phi (z) |}{ 1 - |z|} = 
\frac{1 - |\phi_2 (z)|}{1 -|z|} - |\phi (z)|.
\end{displaymath}
Now (\cite{Shap}, Theorem 3.5), the compactness of $C_{\phi_2} \colon H^2 \to H^2$ implies that
\begin{displaymath}
\lim_{z \stackrel{<}{\to} 1} \frac{1 - |\phi_2(z)|}{1 - |z|} = +\infty;
\end{displaymath}
we get:
\begin{displaymath}
\lim_{z \stackrel{<}{\to} 1} \frac{1 - |\phi_1(z)|}{1 - |z|} = +\infty,
\end{displaymath}
which implies the compactness of $C_{\phi_1} \colon H^2 \to H^2$, thanks to the injectivity of $\phi_1$ 
(\cite{Shap}, Theorem 3.2).\hfill $\square$
\bigskip

\subsection{Compactness on $H^\Psi$ versus compactness on $H^2$}
 
 The equivalence \ref{MC} holds actually for every $H^p$ space (with $p<\infty$) instead of $H^2$. 
 We are going to see in this section that for Hardy-Orlicz spaces $H^\Psi$, one needs a new notion of 
Carleson measures, which one may call \emph{$\Psi$-Carleson measures}. Before that, we are going to 
see that condition \ref{MC} allows to get that the compactness of composition operators on $H^\Psi$ 
always implies that on $H^p$ for $p<\infty$. Recall that, when $\Psi\in \Delta^2$, we have seen in 
Corollary \ref{exemple non compact} that the converse is not true.
\goodbreak

\begin{theoreme}\label{implique Carleson}
Let  $\phi \colon \D \to \D$ be an analytic function. If one of the following conditions: 
\begin{itemize}
\item [i)] $C_\phi$ is a compact operator on $H^\Psi$ 
\item [ii)] $\Psi \in \Delta^0$ and $C_\phi$ is a weakly compact operator on $H^\Psi$
\end{itemize}
is satisfied, the composition operator $C_\phi$ is compact on $H^2$.
\end{theoreme}

Note that we have proved in Theorem \ref{corollaire Delta2} that the weak compactness of 
$C_\phi \colon H^\Psi \to H^\Psi$ implies its compactness only when $\Psi$ satisfies the 
$\Delta^2$  condition. Nevertheless, we shall show in Theorem~\ref{equiv weak-comp} that when 
$\Psi\in \Delta^0$, the weak compactness of  $C_\phi$ is equivalent to its compactness. This is obviously false 
(in particular when $L^\Psi$ is reflexive) without any assumption on $\Psi$.
\medskip

\noindent{\bf Proof.} We are going to use the characterization \ref{MC}  for
compact composition operators on $H^2$.\par 
Suppose that the condition on $\mu_\phi$ is not fulfilled. Then there exist 
$\beta\in(0,1)$, $\xi_n\in \T$, and $h_n\in(0,1)$, with 
$h_n\mathop{\longrightarrow}\limits_{n\to+\infty} 0$, such that: 
\begin{displaymath}
\mu_\phi\big(S(\xi_n, h_n)\big) \geq\beta h_n.
\end{displaymath}
We are now going to use the function:
\begin{displaymath}
v_n(z)= \frac{h_n^2}{(1-\overline{a}_nz)^2}\,\raise 0,5mm \hbox{,}
\end{displaymath} 
where:
\begin{displaymath}
a_n = (1-h_n)\xi_n. 
\end{displaymath}
Of course, $v_n$ is actually nothing but $u_{\xi_n, 1 - h_n}$. We have by Corollary~\ref{norme $u_{a,r}$}:
\begin{displaymath}
\|v_n\|_\Psi \leq \frac{1}{\Psi^{-1}(1/h_n)}\cdot
\end{displaymath}
Define $g_n=\Psi^{-1}(1/h_n) v_n$, which is in the unit ball of $HM^\Psi$. We have assumed 
at the beginning of the paper that $x=o\,\big(\Psi(x)\big)$ as $x\to \infty$; hence 
$\Psi^{-1}(x) = o\,(x)$ as $x\to \infty$, and so  $h_n^2\Psi^{-1}(1/h_n)\rightarrow 0$. 
Therefore $(g_n)_n$ converges uniformly to zero on compact subsets of $\D$ and 
$\|g_n\|_1\rightarrow 0$, because $\|g_n\|_1\leq h_n\Psi^{-1}(1/h_n)$. 
Then, in both cases, we should have $\| C_\phi (g_n)\|_\Psi \to 0$. 
Indeed, in case {\it i)}, this follows from Proposition \ref{critere compact}, and in case {\it ii)},  this 
follows from \cite{LLQR}, Theorem 4.
\smallskip

We are going to show that this is not true. Indeed:
\begin{align*}
\int_\T \Psi\Big(\frac{4}{\beta} |g_n \circ\phi|\Big)\,dm
& \geq \int_\D \Psi\Big(\frac{4}{\beta} \Psi^{-1}(1/h_n) |v_n(z)|\Big)\,d\mu_\phi \\
& \geq \int_{S(\xi_n, h_n)} \Psi\Big(\frac{4}{\beta} \Psi^{-1}(1/h_n) |v_n(z)|\Big)\,d\mu_\phi.
\end{align*}
But when $z\in S(\xi_n, h_n)$, one has $|v_n(z)|\geq 1/4$, because
\begin{displaymath}
|1-\overline{a}_n z|\leq |1-\overline{a}_n \xi_n| + |\overline{a}_n ( \xi_n -z)| 
= h_n +  (1 - h_n) h_n  \leq 2h_n.
\end{displaymath}
We obtain that by convexity (since $\beta<1$): 
\begin{align*}
\int_\T \Psi\Big(\frac{4}{\beta} |g_n\circ\phi|\Big)\,dm 
& \geq \int_\T \frac{1}{\beta} \Psi\big( 4 \Psi^{-1}(1/h_n) |v_n(z)| )\,dm \\ 
& \geq \frac{1}{\beta h_n}\mu_\phi\big(S(\xi_n,h_n)\big) \geq 1.
\end{align*}
This implies that
$\|C_\phi(g_n)\|_\Psi \geq \beta/4$ and proves the theorem.\hfill $\square$
\bigskip

\subsection{General measures}

We used several times the criterion \ref{MC}  for compactness on $H^2$ via Carleson measures. 
The fact that this provides such a useful tool leads to wonder if the boundedness and the 
compactness on Hardy-Orlicz spaces can be expressed in such a pleasant manner. 
Theorem~\ref {continuite Carleson} and Theorem~\ref{compacite Carleson} 
are the Orlicz version of Carleson's Theorem for $H^p$ spaces (\cite{Co-McC}, Theorem 2.35).
\medskip

The key for our general characterization is the use of the following functions (see \eqref{def W} for the 
definition of the Carleson's window $W (\xi, h)$).

\begin{definition}
For any positive finite Borel measure $\mu$ on the unit disk $\D$ (or on $\overline \D$), we set, for $h\in(0,1]$:
\begin{align}
\rho_\mu (h) & =\sup_{\xi\in \T}\mu\big( W (\xi,h)\big), \\
K_\mu (h) & = \sup_{0<t \leq h}\frac{\rho_\mu (t)}{t}
\end{align}
\end{definition}

\noindent Hence $\mu\big( W (\xi, t)\big) \leq K_\mu (h) t$ for $t\leq h$. 
\medskip

The measure $\mu$ is a Carleson measure if and only if $K_\mu (h)$ is bounded by a constant $K$, for $h\in(0,1)$ and this happens as soon as $ K_\mu (h_0)$ is finite  for some $h_0\in(0,1)$.
\medskip

\begin{definition}
We say that $\Psi$ satisfies the $\nabla_0$ condition if for some $x_0>0$, $C\ge1$ and  every  $x_0\le x\leq y$, one has:
\begin{equation}\label{nabla0}
\frac{\Psi (2x)}{\Psi (x)} \leq \frac{\Psi (2Cy)}{\Psi (y)}\cdot 
\end{equation}
\end{definition}

This is a condition on the regularity of $\Psi$. It is satisfied if
\begin{displaymath}
\frac{\Psi (2x)}{\Psi (x)} \leq C\,\frac{\Psi (2y)}{\Psi (y)}\cdot 
\end{displaymath}
\medskip
\begin{proposition}\label{equi-nabla}
the following assertions are equivalent
\begin{itemize}
\item [i)] $\Psi$ satisfies the $\nabla_0$ condition.

\item [ii)]  There exists some $x_0>0$ satisfying: for every $\beta>1$, there exists $C_\beta\ge1$ such that
\begin{displaymath}
\frac{\Psi (\beta x)}{\Psi (x)} \leq \frac{\Psi (\beta C_\beta y)}{\Psi (y)}\;\raise 0,5mm \hbox{,}\hbox{ for every  }x_0\le x\leq y. 
\end{displaymath}

\item [iii)]  There exist $x_0>0$, $\beta>1$ and $C_\beta\ge1$ such that
\begin{displaymath}
\frac{\Psi (\beta x)}{\Psi (x)} \leq \frac{\Psi (\beta C_\beta y)}{\Psi (y)}\;\raise 0,5mm \hbox{,}\hbox{ for every  }x_0\le x\leq y. 
\end{displaymath}
\end{itemize}
\end{proposition}
\medskip

\noindent{\bf Proof.} We only have to prove $i)\Rightarrow ii)$, since $iii)\Rightarrow i)$ is similar and $ii)\Rightarrow iii)$ is trivial.

If $\beta\in(1,2]$, it is easy, taking $C_\beta=2C/\beta$. Now, if $\beta\in(2^b,2^{b+1}]$ for some integer $b\ge1$, we write for every $x_0\le x\leq y$:

\begin{displaymath}
\frac{\Psi (\beta x)}{\Psi (x)} \leq \frac{\Psi (2^{b+1}x)}{\Psi (x)}= \frac{\Psi (2^{b+1}x)}{\Psi (2^b x)}\cdots \frac{\Psi (2x)}{\Psi (x)}\cdot
\end{displaymath}
But we have for every integer $j\ge1$: $2^{j-1}x\le(2C)^{j-1}x\le(2C)^{j-1}y$, so:
\begin{displaymath}
\displaystyle\frac{\Psi (2^jx)}{\Psi (2^{j-1} x)}\le\frac{\Psi ((2C)^jy)}{\Psi ((2C)^{j-1}y)}\raise 0,5mm \hbox{,}
\end{displaymath}
and we obtain:
\begin{displaymath}
\frac{\Psi (\beta x)}{\Psi (x)} \leq \frac{\Psi ((2C)^{b+1}y)}{\Psi (y)}\le\frac{\Psi (\beta C_\beta y)}{\Psi (y)}\raise 0,5mm \hbox{,}
\end{displaymath}
where $C_\beta=2C^{b+1}$.\hfill $\square$
\medskip

\noindent{\bf Examples.} It is immediately seen that the following functions satisfy $\nabla_0$: 
$\Psi (x) = x^p$, $\Psi (x) = \exp \big[\big(\log (x+1)\big)^\alpha\big] - 1$, 
$\Psi (x) = \e^{x^\alpha} - 1$, $\alpha \geq1$.\par
\medskip 

Note that when $\Psi \in \nabla_0$, with constant $C=1$, {\it i.e.} $\Psi (\beta x) / \Psi (x)$ is 
increasing for $x$ large enough, then we have the dichotomy: either $\Psi \in \Delta_2$, or 
$\Psi \in \Delta^0$. 
\medskip 
We shall say that $\nabla_0$ is \emph{uniformly satisfied} if there exist $C\geq 1$ and 
$x_0 > 0$ such that, for every $\beta>1$:
\begin{equation}\label{nabla uniforme}
\quad \quad \quad \frac{\Psi (\beta x)}{\Psi (x)} \leq \frac{\Psi (C\beta y)}{\Psi (y)}\quad \text{for} 
\quad x_0 \leq x \leq y,
\end{equation}
\medskip

One has:
\begin{proposition}\label{log convexe} \ \par
1) Condition $\Delta^2$ implies condition $\nabla_0$ uniformly.\par
2) If $\Psi \in \nabla_0$ uniformly, then $\Psi \in \nabla_1$.\par
3) The function $\kappa (x) = \log \Psi (\e^x)$ is convex on $(x_0,+\infty)$ if and only if $\nabla_0$ is satisfied with constant $C=1$.\par
\end{proposition}

We shall say that \emph{$\Psi$ is $\kappa$-convex} when $\kappa$ is convex at infinity. 
Note that $\Psi$ is $\kappa$-convex whenever $\Psi$ is log-convex. In the above examples $\Psi$ is 
$\kappa$-convex; it also the case of $\Psi (x) = x^2 / \log x$, $x\geq \e$; but, on the other hand, if 
$\Psi (x) = x^2 \log x$ for $x \geq \e$, then $\Psi$ is not 
$\kappa$-convex. Nevertheless, for $\beta^2 \leq x \leq y$, one has:
\begin{displaymath}
\frac{\Psi (\beta x)}{\Psi (x)} = \beta^2 \Big( 1 + \frac{\log \beta}{\log x}\Big) \leq \frac{3 \beta^2}{2} 
\leq \frac{3}{2}\frac{\Psi (\beta y)}{\Psi (y)}\,\raise 0,5mm \hbox{,}
\end{displaymath}
and hence $\Psi \in \nabla_0$.\par
We do not know whether $\Psi \in \nabla_0$ uniformly implies that $\Psi$ is equivalent to an Orlicz 
function for which the associated function $\kappa$ is convex.
\medskip

\noindent{\bf Proof.} 1) Since $\Psi \in \Delta^2$, one has $\big[\Psi (u )\big]^2 \leq \Psi (\alpha u)$ 
for some $\alpha >1$ and $x\geq x_0$. We may assume that $\Psi (x_0) \geq 1$. Then, for 
$y\geq x \geq x_0$ and every $\beta>1$:
\begin{displaymath}
\Psi (\beta x) \Psi (y) \leq \big[ \Psi (\beta y)\big]^2 \leq \Psi (\alpha \beta y) \leq \Psi (\alpha \beta y) \Psi (x) ,
\end{displaymath}
which is \eqref{nabla uniforme}.\par
2) Suppose that $\Psi$ satisfies condition $\nabla_0$ uniformly. We may assume 
that $\Psi (x_0) \geq 1$. Let $x_0 \leq u \leq v$; we can write $u = \beta x_0$ for some $\beta  \geq 1$. Then 
condition~\eqref{nabla uniforme} gives:
\begin{displaymath}
\Psi (u) \Psi (v) = \Psi (\beta  x_0) \Psi (v) \leq \Psi (x_0) \Psi (C\beta v) 
\leq \Psi \big( \Psi (x_0) C \beta  v\big) = \Psi (b uv),
\end{displaymath}
with $b = C \Psi (x_0)/ x_0$.\par
3) Assume that $\kappa$ is convex on $(x_0,+\infty)$. For every $\beta >1$, let $\kappa_\beta(t) = \kappa (t \log \beta )=\log(\Psi(\beta^t))$, which is convex on $(x_0/\log(\beta),+\infty)$. 
Taking $y \geq x\geq{\rm e}^{x_0}$, write $x = \beta^\theta$ and $y = \beta^{\theta'}$ with $\theta \leq \theta'$. 
Convexity of $\kappa$ gives, since $\theta'\geq \theta\geq x_0/\log(\beta)$:
\begin{displaymath}
\kappa_\beta  (\theta +1) - \kappa_\beta  (\theta) \leq \kappa_\beta  (\theta' +1) - \kappa_\beta  (\theta'),
\end{displaymath}
which means that:
\begin{displaymath}
\frac{\Psi (\beta x)}{\Psi (x)} \leq \frac{\Psi (\beta y)}{\Psi (y)}\cdot
\end{displaymath}
\par
Assume that \eqref{nabla uniforme} is fulfilled for every $\beta >1$, with $C=1$. 
Then, taking $y =\beta x$, one has:
\begin{displaymath}
\big[\Psi (\beta x )\big]^2 \leq \Psi (x) \Psi (\beta ^2 x).
\end{displaymath}
Let $u < v$ be large enough. Taking $x= u^2$ and $\beta  = v/u$, we get:
\begin{displaymath}
\big[\Psi (uv )\big]^2 \leq \Psi (u^2) \Psi (v^2),
\end{displaymath}
which means that $\kappa$ is convex.\hfill $\square$
\medskip\goodbreak

\noindent{\bf Remark.} The growth and regularity conditions for $\Psi$ can be expressed in the 
following form:
\begin{itemize}
\item $\Psi \in \Delta^0$ iff  
$\kappa (x + \beta') - \kappa (x) \mathop{\longrightarrow}\limits_{x \to \infty} +\infty$, 
for some $\beta'>0$.
\item $\Psi \in \Delta^1$ iff for some $\beta'>0$, one has
$ x + \kappa (x) \leq \kappa (x + \beta')$, for $x$ large enough.
\item $\Psi \in \Delta^2$ iff for some $\alpha'>0$, one has 
$ 2 \kappa (x) \leq \kappa (x + \alpha')$, for $x$ large enough.
\item $\Psi \in \nabla_1$ iff 
$\kappa_B (x) + \kappa_B (y) \leq \kappa_B (x + y)$ for $x,y$ large enough, with  
$B = \e^{-b}$.
\item $\Psi \in \nabla_0$ iff for some $c' \geq 1$ and $A>1$, one has 
$\kappa_A (\theta + 1) - \kappa_A (\theta) \leq \kappa_A (\theta' + c') - \kappa_A (\theta')$ 
for $\theta \leq \theta'$ large enough.
\end{itemize}

\medskip\goodbreak
Before proving the main results of this section, let us collect some basic facts on the compactness of the embedding of $H^{\Psi_1}$ into $L^{\Psi_2}(\mu)$. First:
\goodbreak
\begin{lemme}
Let $\Psi_1$, $\Psi_2$ be two Orlicz functions and $\mu$ a finite Borel measure on $\overline \D$.  Assume that the identity maps $H^{\Psi_1}$ into $L^{\Psi_2}(\mu)$ compactly. Then, $\mu(\T)=0$.
\end{lemme}

\noindent{\bf Proof.} The sequence $(z^n)_{n\ge1}$ is weakly null in $M^{\Psi_1}$ by Riemann-Lebesgue's Lemma. Its image by a compact operator is then norm null. This implies that for every $\eps\in(0,1)$, we have, for $n$ large enough,
\begin{displaymath}
\int_{\overline\D}\Psi_2\bigl(|z|^n\bigr)d\mu\le\eps.
\end{displaymath}
Fatou's Lemma yields $\Psi_2(1)\mu(\T)\le\eps$.\hfill $\square$
\bigskip

Now, we summarize what is true in full generality about compactness for canonical embeddings.

\begin{proposition}\label{fullgene}
Let $\Psi_1$, $\Psi_2$ be two Orlicz functions and $\mu$ a finite Borel measure on $\overline \D$.  The following assertions are equivalent
\begin{itemize}
\item [i)] The identity  maps $H^{\Psi_1}$ into $L^{\Psi_2}(\mu)$ compactly.

\item [ii)] Every sequence in the unit ball of $H^{\Psi_1}$, which is convergent to $0$ uniformly on every compact subset of $\D$, is norm-null in $L^{\Psi_2}(\mu)$.

\item [iii)] The identity maps $H^{\Psi_1}$ into $L^{\Psi_2}(\mu)$ continuously and\par  $\mathop{\lim}\limits_{r\to1^-}\bigl\|I_r\bigr\|=0$, where $I_r(f)=f\ind_{\overline\D\setminus r\D}$.
\end{itemize}
\end{proposition}

\noindent{\bf Proof.} $i)\Rightarrow ii)$: let $(f_n)_{n\ge1}$ be a sequence in the unit ball of $H^{\Psi_1}$, which is uniformly convergent to $0$ on every compact subset of $\D$. In particular, $f_n(z)$ converges to $0$ for every $z\in\D$. This means that $(f_n)_{n\ge1}$ converges to $0$ $\mu$-almost everywhere, since $\mu(\T)=0$ by the preceding lemma. If the conclusion did not hold, we could assume, up to an extraction, that $\underline{\lim}\|f_n\|_{\Psi_2}>0$. Thus by compactness of the embedding, up to a new extraction, $(f_n)_{n\ge1}$ is norm-convergent to some $g\in L^{\Psi_2}(\mu)$. Necessarily $g\not=0$. A subsequence of $(f_n)_{n\ge1}$ would be convergent to $g$ $\mu$-almost everywhere. This gives a contradiction.

$ii)\Rightarrow iii)$: if not, there exist a sequence $(f_n)_{n\ge1}$ in the unit ball of $H^{\Psi_1}$ and $\delta>0$ with $\|f_n\ind_{\overline\D\setminus(1-\frac{1}{n})\D}\|_{\Psi_2}>\delta$, for every $n\ge1$. Let us introduce the sequence $g_n(z)=z^nf_n(z)$ for $z\in\D$. The sequence $(g_n)_{n\ge1}$  lies in the unit ball of $H^{\Psi_1}$ and is convergent to $0$ uniformly on every compact subset of $\D$. But
\begin{displaymath}
\|g_n\|_{\Psi_2}\ge\|z^nf_n\ind_{\overline\D\setminus(1-\frac{1}{n})\D}\|_{\Psi_2}\ge\Bigl(1-\frac{1}{n}\Bigr)^n\|f_n\ind_{\overline\D\setminus(1-\frac{1}{n})\D}\|_{\Psi_2}\ge\Bigl(1-\frac{1}{n}\Bigr)^n\delta.
\end{displaymath}
This contradicts $(ii)$.

$iii)\Rightarrow ii)$ is very easy.

$ii)\Rightarrow i)$ follows from Proposition~\ref{critere compact}.\hfill $\square$

\bigskip
We can now state some deeper characterizations 

\begin{theoreme}\label{continuite Carleson} 
Let $\mu$ be a finite Borel measure on the closed unit disk $\overline \D$ and  
let $\Psi_1$ and $\Psi_2$ be two Orlicz functions. Then:\par
1) If the identity  maps $H^{\Psi_1}$ into $L^{\Psi_2}(\mu)$ continuously, there exists some $A>0$ such that:
\begin{equation}\label{eq:rho}
\rho_\mu (h)\leq \frac{1}{\Psi_2\big(A\Psi_1^{-1}(1/h)\big)}\hskip 2mm \text{for every } h\in (0,1].
\tag{$R$}
\end{equation}
\par
2) If there exists some $A>0$ such that:
\begin{equation}\label{eq:K}
K_\mu (h) \leq \frac{1/h}{\Psi_2\big(A\Psi_1^{-1}(1/h)\big)} \hskip 2mm \text{for every } h\in(0,1],
\tag{$K$}
\end{equation}
then the identity  maps $H^{\Psi_1}$ into $L^{\Psi_2}(\mu)$ continuously.\par

\end{theoreme}

\begin{theoreme}\label{compacite Carleson}
Let $\mu$ be a finite Borel measure on the closed unit disk $\overline\D$ and  
let $\Psi_1$ and $\Psi_2$ be two Orlicz functions. Then:\par
1) If the identity  maps $H^{\Psi_1}$ into $L^{\Psi_2}(\mu)$ compactly, then 
\begin{equation}\label{eq:rho compact}
\rho_\mu (h) =o\,\Big( \frac{1}{\Psi_2\big(A\Psi_1^{-1}(1/h)\big)}\Big) \hskip 2mm \text{as } h\to 0,\hbox{ for every   }A>0.
\tag{$R_0$}
\end{equation}
\par
2) If  $\mu(\T)=0$ and 
\begin{equation}\label{eq:K compact}
K_\mu (h) =o\, \Big(\frac{1/h}{\Psi_2\big(A\Psi_1^{-1}(1/h)\big)}\Big) \hskip 2mm \text{as } h\to 0\, ,\hbox{ for every   }A>0,
\tag{$K_0$}
\end{equation}
then the identity  maps $H^{\Psi_1}$ into $L^{\Psi_2}(\mu)$ compactly.
\par
3) When $\Psi_1=\Psi_2=\Psi$ satisfies condition $\nabla_0$, then the above conditions are 
equivalent: the identity  maps $H^{\Psi_1}$ into $L^{\Psi_2}(\mu)$ compactly if and only if  
condition \eqref{eq:rho compact} is satisfied and if and only if condition 
\eqref{eq:K compact} is satisfied.
\end{theoreme}

\noindent{\bf Remarks.} \par
{\bf 1.} a) If $\rho_\mu (h) \leq C/ \Psi_2\big(A\Psi_1^{-1}(1/h)\big)$ with $C>1$, then 
the convexity of $\Psi_2$ gives $\Psi_2(t/C) \leq \Psi_2(t)/C$ and hence:
\begin{displaymath}
\rho_\mu(h) \leq \frac{1}{\Psi_2\big( \frac{A}{C}\Psi_1^{-1}(1/h)\big)}\cdot
\end{displaymath}
\par
b) If $\rho_\mu(h) \leq 1/ \Psi_2\big(A\Psi_1^{-1}(1/h)\big)$ only for $h\leq h_A$, one can find some 
$C=C_A >0$ such that $\rho_\mu (h) \leq C/ \Psi_2\big(A\Psi_1^{-1}(1/h)\big)$ for every $h\in (0,1]$. 
In fact, $\rho_\mu (h) \leq \mu (\overline \D)$ and 
$1/ \Psi_2\big(A\Psi_1^{-1}(1/h)\big) \geq 1/ \Psi_2\big(A\Psi_1^{-1}(1/h_A)\big)$ for 
$h \geq h_A$; hence $\rho_\mu (h) \leq C/ \Psi_2\big(A\Psi_1^{-1}(1/h)\big)$ with 
$C= \mu (\overline \D)/ \Psi_2\big(A\Psi_1^{-1}(1/h_A)\big)$.\par
The same remark applies for $K_\mu$.\par
\medskip

{\bf 2.} a) In the case where $\Psi_1=\Psi_2 =\Psi$, one has:  
$\Psi\big(A\Psi^{-1}(t)\big)  \leq \Psi\big( \Psi^{-1}(t)\big) = t$ when $A\leq 1$ and 
$\Psi\big(A\Psi^{-1}(t)\big)  \geq \Psi\big( \Psi^{-1}(t)\big) = t$ when $A\geq 1$. On the other hand, 
if $\Psi \in \Delta_2$, one has, for some constant $C=C_A>0$:  
$\Psi (A x)\leq C \Psi (x)$, when $A \geq 1$ and 
$\Psi (A x)\geq (1/C) \Psi (x)$ when $A \leq 1$. Hence, when $\Psi\in \Delta_2$, one has, for every $A>0$:
\begin{displaymath}
\frac{1}{\Psi \big(A\Psi ^{-1}(1/h)\big)} = \frac{1}{\chi_A(1/h)} \approx h
\end{displaymath}
and Theorem \ref{continuite Carleson} is nothing but Carleson's Theorem.\par
b) If $\Psi_1 (x)= x^p$ and $\Psi_2 (x)= x^q$ with $p < q < \infty$, then:
\begin{displaymath}
\Psi_2\big(A \Psi_1^{-1} (t)\big) = A^q t^{q/p},
\end{displaymath}
and condition \eqref{eq:rho} means that $\mu$ is a $\beta$-Carleson measure, with $\beta=q/p$ 
(see \cite{Du}, Theorem 9.4).\par
c)  If, for fixed $A>0$, the function $\displaystyle x\mapsto\frac{\Psi_1 (x)}{\Psi_2 (Ax)}$ is non
increasing, at least for $x$ large enough, conditions \eqref{eq:rho} and \eqref{eq:K} (resp. 
conditions \eqref{eq:rho compact} and \eqref{eq:K compact} below) are clearly equivalent. This is the 
case in the framework of classical Hardy spaces: $\Psi_1 (x) = x^p$ and $\Psi_2 (x) = x^q$, with 
$q\geq p$.\par
When $\Psi_1 =\Psi_2 = \Psi$, this is equivalent, if $A >1$, to the convexity of the function 
$\kappa (x) = \log \Psi (\e^x)$ (see Proposition~\ref{log convexe}).
\medskip

{\bf 3.} a) When $\Psi_1=\Psi_2=\Psi$, the condition $\mu(\T)=0$ is automatically fulfilled (and so can be removed from \eqref{eq:K compact} ). This follows on one hand from the majorization in \eqref{eq:K compact} , which implies that $K_\mu(h)\to0$ (when $h\to0$); and on the other hand from the inequality:
\begin{displaymath}
\mu\big(\overline\D\setminus r\D\big)\le\frac{2\pi}{1-r}\rho_\mu(1-r)\le2\pi K_\mu(1-r).
\end{displaymath}
Indeed, $\overline\D\setminus r\D$ can be covered by less than $\frac{2\pi}{1-r}$ Carleson's windows of ``size'' $1-r$.

b) Nevertheless, the condition $\mu(\T)=0$ cannot be removed in full generality in Theorem \ref{compacite Carleson}. Indeed, if we consider the identity $j$ from $H^4$ into $L^2(\overline\D,\tilde m)$, where $\tilde m$ is $0$ on $\D$ and its restriction to the torus is the normalized Lebesgue measure. It is easily seen that $K(h)$ is bounded and so less than $\displaystyle\frac{1}{A^2h^{1/2}}$, for $h$ small enough. But $j$ is not compact.
\medskip

{\bf 4.} In the case where $\Psi_1=\Psi_2 =\Psi$ and $\mu$ is a Carleson measure, 
then $K_\mu$ is bounded, by say $K \geq 1$, and condition \eqref{eq:K} is satisfied for $A=1/K$, since 
$A\leq 1$ implies, by the convexity of $\Psi$: 
$\Psi\big(A\Psi^{-1}(1/h)\big) \leq A \Psi\big(\Psi^{-1}(1/h)\big) = A/h$. Hence the canonical 
embedding $H^\Psi \hookrightarrow L^\Psi (\mu)$ is continuous. We get hence, by Carleson's Theorem 
(\cite{Du}, Theorem 9.3):
\medskip

\begin{proposition}
Let $\mu$ be a positive finite measure on $\overline\D$. Assume that the canonical embedding 
$j_\mu \colon H^p \to L^p (\mu)$ is continuous for some $ 0 < p < \infty$. Then 
$j_\mu \colon H^\Psi \to L^\Psi (\mu)$ is continuous.
\end{proposition}

Note that this is actually a consequence of the fact that $H^\Psi$ is an interpolation space for $H^1$ and 
$H^\infty$ (see \cite{Be-Sh}, Theorem V.10.8).\par
When $\mu=\mu_\phi$ is the image of the Haar measure $m$ under $\phi^\ast$, where $\phi$ is an 
analytic self-map of $\D$, we know (Proposition \ref{subordination}) that the composition operator 
$C_\phi \colon H^\Psi \to H^\Psi$ is always continuous. This can be read as the continuity of 
$H^\Psi \hookrightarrow L^\Psi(\mu_\phi)$. Hence condition \eqref{eq:rho} 
must be satisfied, for some $A>0$. Note that for $A \leq 1$, $1/\chi_A(1/h) \geq h$, and so 
condition \eqref{eq:rho} is implied by the fact that $\mu_\phi$ is a Carleson measure.
\medskip

{\bf 5.} To have the majorization in condition \eqref{eq:K compact}, it suffices that: for every $A>0$, there exists 
$h_A\in (0,1]$ such that 
\begin{equation}\label{eq:K compact, variante}
\hskip20pt K_\mu (h) \leq \frac{1/h}{\Psi_2\big(A\Psi_1^{-1}(1/h)\big)}
\,\raise 0,5mm\hbox{,} \hskip 2mm \text{for every } h\in (0,h_A].
\end{equation}
In fact, fixing $A>0$ and $\eps\in(0,1)$, we have, by convexity:
\begin{displaymath}
\Psi_2 (Ax) \leq \eps \Psi_2 (\tilde A x), 
\end{displaymath}
with $\tilde A=A/\eps$. Since we have \eqref{eq:K compact, variante} with $\tilde A$, when 
$h$ is small enough (depending on $A$ and $\eps$), we get, for $x=\Psi_1^{-1}(1/h)$,  
condition \eqref{eq:K compact}.
\bigskip

To prove both Theorem \ref{continuite Carleson} and Theorem \ref{compacite Carleson}, we 
shall need some auxiliary results. The following is actually the heart of the classical Theorem of 
Carleson, though it is not usually stated in this form. The maximal (non-tangential) function $M_f$ 
(which is essentially the same as $N_\alpha f$ in the previous section) will be defined by:

\begin{displaymath}
M_f\big(\e^{i\theta}\big) =\sup\{|f(z)|;\; z\in G_\theta\} , 
\end{displaymath}
where
\begin{displaymath}
G_\theta=\{z\in \D \,;\  |\e^{i\theta} - z|<3 (1-|z|)\}.
\end{displaymath}

\begin{theoreme}[Carleson's Theorem]\label{Carleson}
For every $f\in H^1$ and every finite positive measure $\mu$ on the closed unit disk $\overline\D$, one has, 
for every $h\in(0,1]$ and every $t>0$:
\begin{displaymath}
\mu\big(\{z\in \overline{\D}\,;\ |z|> 1-h \text{ and } |f(z)|>t\}\big) \leq 2\pi K_\mu (h)\, m(\{M_f >t\}).
\end{displaymath}
\end{theoreme}

As this theorem is not usually stated in such a way,  we shall give a few words of 
explanations.
\medskip

\noindent{\bf Proof.} For convenience , we shall denote, when $I$ is a subarc of $\T$:
\begin{displaymath}
W(I)=\{z\in \overline{\D}\,;\ |z|>\max(0,1-|I|/2) \text{ and }\frac{z}{|z|}\in I\}.
\end{displaymath}
Obviously, when $|I|\le2$, we have $W(I)=W(\xi,|I|/2)$, where $\xi$ is the center of $I$.

We shall begin by being somewhat sketchy, and refer to 
\cite{Co-McC}, Theorem 2.33, or \cite{Du}, Theorem 9.3, for the details. Following the lines of 
\cite{Co-McC}, page 39, proof of Theorem 2.33, $1) \Rightarrow 2)$: $\{M_f >t\}$ is the disjoint union 
of a countable family of open arcs $I_j$ of $\T$, and on the other hand (see \cite{Co-McC}, page 39), 
$|f (z)| > t$ implies that $z \in W (I_j)$ for some $j$.\par

Now, when $\alpha_j > h$, we have to cover $I_j$ in such a way that we can write 
$I_j\subset J_{j,1} \cup \cdots \cup J_{j,N_j}$ with the arcs 
$J_{j,1}, \ldots, J_{j,N_j}$ satisfying: $W (J_{j,k}) = W (\xi_{j,k}, h)$ for every 
$k=1, \ldots, N_j$ and $2|I_j|\ge\displaystyle\sum_{k=1}^{N_j}|J_{j,k}|$, and some $\xi_{j,k}\in\T$.\par

We then notice that:
\begin{equation}\label{petite mesure}
\mu\big( W (J_{j,k})\big) \leq \frac{1}{2}K_\mu (h)\,|J_{j,k}|\,.
\end{equation}
In fact, $2\alpha=|I|$ if $W (\xi, \alpha) = W (I)$; hence:
\begin{displaymath}
\mu \big( W (J_{j,k})\big) \leq \rho_\mu (h) \leq h K_\mu (h) 
=\frac{1}{2}|J_{j,k}|\,K_\mu (h).
\end{displaymath}
Denoting, for $E \subseteq \overline{\D}$, by $E_h$ the set of points $z\in E$ such that 
$|z| > 1-h$, we therefore have, since: 
\begin{displaymath}
W (I_j)_h \subseteq\bigcup_{1 \leq k \leq N_j}   W (J_{j,k}),
\end{displaymath}
using \eqref{petite mesure}:
\begin{displaymath}
\mu\big( W (I_j)_h \big)\leq \sum_{k=1}^{N_j} \mu \big( W (J_{j,k}) \big)\leq \sum_{k=1}^{N_j} \frac{1}{2} K_\mu (h) \,|J_{j,k}|\leq  K_\mu (h) \,|I_j|.
\end{displaymath}
It follows that:
\begin{align*}
\mu\big(\{z\in \overline{\D}\,;\ |z| > 1 - h &\text{ and } |f(z)|>t\}\big) 
\leq \sum_j \mu \big( (W_j)_h \big) \leq K_\mu (h) \sum_j |I_j| \\ 
& =2 \pi K_\mu (h) \sum_j m(I_j)   = 2\pi K_\mu (h)\, m( \{ M_f > t\}),
\end{align*}
as announced.\hfill $\square$
\medskip

The following estimation will be useful for the study both of boundedness and compactness.

\begin{lemme}\label{estimation Carleson}
Let $\mu$ be a finite Borel measure on $\overline\D$. Let $\Psi_1$ and $\Psi_2$ 
be two Orlicz functions. Suppose that there exists $A>0$ and $h_A\in(0,1)$ such that 
\begin{displaymath}
K_\mu(h)\le\frac{1/h}{\Psi_2\big(A\Psi_1^{-1}(1/h)\big)}\raise2pt\hbox{,} \hbox{  for every }h\in(0,h_A).
\end{displaymath}
Then, for every $f\in H^{\Psi_1}$ such that $\|f\|_{\Psi_1} \leq 1$ and every Borel subset $E$ of $\overline\D$, we have:
\begin{displaymath}
\int_{E} \Psi_2 \Big(\frac{A}{8}|f|\Big)\,d\mu \leq\mu(E)\Psi_2(x_A)+\frac{\pi}{2} \int_\T \Psi_1 (M_f)\,dm
\end{displaymath}
where $x_A=\displaystyle\frac{A}{2}\Psi_1^{-1}(1/h_A)$.
\end{lemme}
\medskip

\noindent{\bf Proof.} For every $s>0$, the inequality 
$|f(z)|> s$ implies that the norm of the evaluation at $z$ is greater than $s$; hence by 
Lemma~\ref{norme evaluation}:
\begin{displaymath}
s < 4 \Psi_1^{-1} \Big( \frac{1}{ 1 - |z|}\Big)\,\raise 0,5mm \hbox{,}
\end{displaymath}
{\it i.e.}:
\begin{displaymath}
|z|> 1- \frac{1}{\Psi_1 (s/4)}\cdot
\end{displaymath}
Carleson's Theorem (Theorem \ref{Carleson}) gives:
\begin{displaymath}
\mu\big(\{|f(z)|> s\}\big) \leq 2\pi K_\mu \Big(\frac{1}{\Psi_1 (s/4)}\Big)\, m(\{M_f > s\})
\end{displaymath}
when $ \Psi_1 (s/4) \geq 1$. Hence:
\begin{displaymath}
\int_{E} \Psi_2\Big(\frac{A}{8}|f|\Big)\,d\mu=\int_0^\infty \Psi_2' (t)\,\mu(\{|f |> 8t/A\}\cap E)\,dt\,.
\end{displaymath}
But our hypothesis means that, when $\Psi_1 (s/4)>1/h_A$:
\begin{displaymath}
\displaystyle K_\mu \Big(\frac{1}{\Psi_1 (s/4)}\Big)\leq\frac{\Psi_1 (s/4)}{\Psi_2 (As/4)}\,\cdot
\end{displaymath}
We have then 
\begin{displaymath}
\mu\big(\{|f(z)|>8t/A\}\big)\leq2\pi\frac{\Psi_1 (2t/A)}{\Psi_2 (2t)}\, m(\{M_f >8t/A\}).
\end{displaymath}
So:
\begin{align*}
\int_{E} \Psi_2\Big(\frac{A}{8}|f|\Big)\,d\mu&\leq\int_0^{x_A} \Psi'_2 (t)\mu(E)\,dt\\
&\hskip20pt +2\pi\int_{x_A}^{+\infty} \Psi'_2 (t)\frac{\Psi_1 (2t/A)}{\Psi_2 (2t)}\, m(\{M_f >8t/A\})\,dt\\
&\leq\Psi_2 (x_A)\mu(E)\\
&\hskip20pt +2\pi\int_{x_A}^{+\infty}\frac{ \Psi'_2 (t)}{\Psi_2 (2t)}\,\Psi_1 (2t/A)\, m(\{M_f >8t/A\})\,dt\,.\\
\end{align*}

For the second integral, note that one has $ \Psi (x) \leq x \Psi' (x) \leq \Psi (2x)$, for any Orlicz function 
$\Psi$. This leads to: 

\begin{align*}
\int_{x_A}^\infty\frac{ \Psi'_2 (t)}{\Psi_2 (2t)}\,\Psi_1 (2t/A)&\,m (\{M_f > 8t/A\})\,dt \\
& \leq \int_0^\infty \frac{\Psi_1 (2t/A)}{t}\,m(\{M_f > 8t/A\})\,dt \\
& \leq \frac{2}{A}\int_0^\infty \Psi'_1 (2t/A)\,m(\{M_f > 8t/A\})\,dt \\
& = \int_0^\infty  \Psi_1' (x)\, m(\{M_f > 4x\})\,dx 
= \int_\T \Psi_1\Big (\frac{1}{4}M_f \Big)\,dm \\
&  \leq \frac{1}{4} \int_\T \Psi_1(M_f)\,dm.
\end{align*}
which leads to the desired result.\hfill $\square$
\bigskip

For the proofs of Theorem \ref{continuite Carleson} and Theorem \ref{compacite Carleson}, we may restrict ourselves to the case of functions $\Psi_1$ and $\Psi_2$ satisfying $\nabla_2$. Indeed, suppose that $\Psi_1$ and $\Psi_2$ are Orlicz functions and define $\widetilde{\Psi_j}(t)=\Psi_j(t^2)$, for $j\in\{1,2\}$. The functions  $\widetilde{\Psi_1}$ and $\widetilde{\Psi_2}$ are Orlicz functions satisfying $\nabla_2$ since, with $\beta=2$, we have for every $t\ge0$:
\begin{displaymath}
\widetilde{\Psi_j}(\beta t)=\Psi_j(4t^2)\ge4\Psi_j(t^2)=2\beta\widetilde{\Psi_j}(t).
\end{displaymath}

Now, we claim that $\mu$ satisfies \eqref{eq:rho}, \eqref{eq:K}, \eqref{eq:rho compact} or \eqref{eq:K compact} for the couple $\big(\Psi_1,\Psi_2\big)$ if and only if $\mu$ satisfies it for the couple $\big(\widetilde{\Psi_1},\widetilde{\Psi_2}\big)$. This is simply due to the fact that for every $A>0$ and  $t\ge0$, we have 
\begin{displaymath}
\widetilde{\Psi_2}\Big(A\widetilde{\Psi_1}^{-1}(t)\Big)=\Psi_2\Big(A^2\Psi_1^{-1}(t)\Big).
\end{displaymath}
Moreover, notice that, writing $f=Bg^2$ (where $B$ is a Blaschke product), we have $f\in H^{\Psi}$ if and only if $g\in H^{\widetilde{\Psi}}$; thus  $\|g\|_{L^{\widetilde{\Psi}}}=\sqrt{\|f\|_{L^{\Psi}}}$. It is then clear that 
\begin{displaymath}
\Big\|Id:H^{\Psi_1}\longrightarrow L^{\Psi_2}(\mu)\Big\|=\Big\|Id:H^{\widetilde{\Psi_1}}\longrightarrow L^{\widetilde{\Psi_2}}(\mu)\Big\|^2,
\end{displaymath}
so that the canonical embedding is bounded (resp. compact) for the couple $\big(\Psi_1,\Psi_2\big)$ if and only if it is so for the couple $\big(\widetilde{\Psi_1},\widetilde{\Psi_2}\big)$, thanks to Proposition~\ref{fullgene}.
\bigskip

\noindent{\bf Proof of Theorem~\ref{continuite Carleson}.} 1) Let $C$ be the norm of the 
canonical embedding $j\colon H^{\Psi_1} \hookrightarrow L^{\Psi_2}(\mu)$, and 
let $\xi \in \T$ and $h\in (0,1)$. 
It suffices to test the continuity  of $j$ on $f= \Psi_1^{-1} (1/h) u_{\xi, 1 -h}$, which is 
in the unit ball of $HM^{\Psi_1}$, by Corollary~\ref{norme $u_{a,r}$}. \par
But, when $z\in W (\xi, h)$ one has, with $a =(1-h) \xi$:
\begin{align*}
| 1 - \bar{a} z| 
& \leq | 1 -\bar{a} \xi | + |\bar{a} \xi - \bar{a} z | 
= h + (1 -h ) \Big[ \Big| \xi - \frac{z}{|z|} \Big| + \Big|  \frac{z}{|z|} - z\Big| \Big] \\
& \leq h + (1 -h ) [ h + (1 - |z|) ] \leq h + (1-h) [ h +h] \leq 3h;
\end{align*}
hence $|u_{\xi, 1 -h}(z)| \geq 1/9$ and $|f(z)| \geq (1/9) \Psi_1^{-1} (1/h)$; therefore:
\begin{displaymath}
1 \geq \int_{\overline \D} \Psi_2\Big(\frac{|f|}{C}\Big)\,d\mu 
\geq \Psi_2\Big(\frac{1}{9 C} \Psi_1^{-1} (1/h)\Big)\,\mu\big(W (\xi, h)\big),
\end{displaymath}
which is \eqref{eq:rho}.\par
2) By Proposition~\ref{Marcinkiewicz}, the maximal (non-tangential) function $M$ is bounded on 
$L^{\Psi_1}(\T)$: there exists a constant $C \geq 1$ such that 
$\|M_f\|_{\Psi_1} \leq C \|f\|_{\Psi_1}$ for every $f\in L^{\Psi_1}(\T)$. We fix $f$ in the unit ball 
of $H^{\Psi_1}$ (note that $\| f /C\|_{\Psi_1}$ remains $\leq 1$) and use Lemma \ref{estimation Carleson}, with $E=\overline\D$ and $f$ replaced by $f/C$ (here $h_A=1$)). Writing $\tilde C=\displaystyle\frac{\pi}{2}+\mu(\overline\D)\Psi_2(x_A)$, we get:
\begin{align*}
\int_{\overline \D} \Psi_2\Big( \frac{A}{8C\tilde C} |f|\Big)\,d\mu 
& \leq\frac{1}{\tilde C}\int_{\overline \D} \Psi_2\Big( \frac{A}{8C} |f|\Big)\,d\mu \\
& \leq\frac{1}{\tilde C}\Bigg(\mu(\overline\D)\Psi_2(x_A)+\frac{\pi}{2}\int_\T \Psi_1\Big(\frac{1}{C} M_f\Big)\,dm\Bigg)\\
&\leq\frac{1}{\tilde C}\Big(\mu(\overline\D)\Psi_2(x_A)+\frac{\pi}{2}\Big)=1,
\end{align*}
which means that $\|f\|_{L^{\Psi_2} (\mu)} \leq\displaystyle\frac{8C\tilde C}{A}\cdot$\hfill$\square$
\bigskip

\noindent{\bf Proof of Theorem~\ref{compacite Carleson}.} 1) Suppose that the embedding is compact, 
but that condition~\eqref{eq:rho compact} is not satisfied. Then there exist $\eps_0\in(0,1)$, $A >0$, a 
sequence of positive numbers $(h_n)_n$ decreasing to $0$, and a sequence of $\xi_n\in\T$, such that:
\begin{displaymath}
\mu \big(W (\xi_n, h_n)\big) \geq 
\frac{\eps_0}{\Psi_2\big(A \Psi_1^{-1}(1/h_n)\big)}\cdot
\end{displaymath}
Consider the sequence of functions
\begin{displaymath}
f_n(z)={\Psi_1}^{-1}(1/h_n)\frac{h_n^2}{(1-\bar{a}_n z)^2} = 
{\Psi_1}^{-1}(1/h_n)\, u_{\xi_n, |a_n|}, 
\end{displaymath}
where $a_n = (1 - h_n)\xi_n$. By Corollary~\ref{norme $u_{a,r}$}, $f_n$ is in the unit ball of
$HM^{\Psi_1}$. Moreover, it is plain that $(f_n)_n$ converges to $0$ uniformly on every compact subset of 
$\D$. By the compactness criterion (Proposition~\ref{critere compact}), $(f_n)_n$ is norm-converging to $0$ 
in $L^{\Psi_2} (\mu)$.\par 
But, as above (proof of Theorem~\ref{continuite Carleson}), for every $n\geq 1$, one has 
$|f_n(z)| \geq (1/9) \Psi_1^{-1}(1/h_n)$ when $z\in W (\xi_n, h_n)$; hence:
\begin{align*}
\int_{\overline{\D}}\Psi_2 \Big(\frac{9A}{\eps_0} |f_n|\Big)\,d\mu
& \geq \Psi_2 \Big(\frac{A}{\eps_0} \Psi_1^{-1}(1/h_n)\Big)\,\mu\big(W (\xi_n, h_n)\big) \\
& \geq \Psi_2 \Big(\frac{A}{\eps_0} \Psi_1^{-1}(1/h_n)\Big)\,
\frac{\eps_0}{\Psi_2\big(A \Psi_1^{-1}(1/h_n)\big)} \geq  1, 
\end{align*}
by the convexity of $\Psi_2$. This implies that $\|f_n\|_{L^{\Psi_2} (\mu)}\geq \eps_0/9A$ and gives a 
contradiction.
\par\smallskip

2) We have to prove that for every $\eps>0$, there exists an $r\in(0,1)$ such that the norm of the injection $I_r:H^{\Psi_1}\longrightarrow L^{\Psi_2}(\overline\D\setminus r\D,\mu)$ is smaller than $\eps$ (see Proposition~\ref{fullgene}).

Let $C \geq 1$ be the norm of the maximal operator, as in the proof of Theorem~\ref{continuite Carleson}: 
$\|M_f\|_{\Psi_1} \leq C \|f\|_{\Psi_1}$ for every $f\in L^{\Psi_1}(\T)$, and 
set $A=16 C/\eps$. Condition \eqref{eq:K compact} gives us $h_A\in(0,1)$ such that:
\begin{displaymath}
K_\mu (h) \leq \frac{1}{2} \frac{1/h} {\Psi_2\big(A\Psi_1^{-1}(1/h)\big)} 
\end{displaymath}
when $h \leq h_A$.

Let $f$ in the unit ball of $H^{\Psi_1}$ and $r\in(0,1)$. By Lemma~\ref{estimation Carleson}:

\begin{align*}
\int_{\overline{\D}\setminus r{\D}} \Psi_2 \Big(\frac{|f|}{\eps}\Big)\,d\mu 
& = \int_{\overline{\D}\setminus r{\D}} \Psi_2 \Big(\frac{A}{16 C} |f|\Big)\,d\mu
\leq\frac{1}{2}\int_{\overline{\D}\setminus r{\D}} \Psi_2 \Big(\frac{A}{8 C} |f|\Big)\,d\mu\\
&\leq\frac{1}{2}\Bigg(\mu\big(\overline{\D}\setminus r{\D}\big)\Psi_2(x_A)+\frac{\pi}{2} \int_\T \Psi_1\bigg(\frac{M_f}{C}\bigg)\,dm\Bigg)\\
&\leq\frac{\pi}{4}+ \Psi_2(x_A)\mu\big(\overline{\D}\setminus r{\D}\big).\\
\end{align*}
As $\mu(\T)=0$, there exists some $r_0\in(0,1)$ such that $\frac{\pi}{4}+ \Psi_2(x_A)\mu\big(\overline{\D}\setminus r{\D}\big)\le1$, for every $r\in(r_0,1)$.
This ends the proof of 2).
\medskip

3) Assume that $\Psi$ satisfies condition $\nabla_0$ and that condition \eqref{eq:rho compact} is fulfilled: for every $A>0$ and every $h\in(0,h_A)$ ($h_A$ small enough), we have:
\begin{displaymath}
\rho_\mu (h) \leq \frac{1}{\Psi \big[A\Psi^{-1}(1/h)\big]}\,\cdot
\end{displaymath}
This implies that:
\begin{displaymath}
K_\mu(h) =\sup_{0<s \leq h} \frac{\rho_\mu (s)}{s }\leq 
\sup_{0<s \leq h} \frac{1/s}{\Psi \big[A\Psi^{-1}(1/s)\big]} 
=\sup_{x \geq \Psi^{-1}(1/h)} \frac{\Psi (x)}{\Psi (A\,x)}\cdot
\end{displaymath}
Fix an arbitrary $\beta>1$ and choose $A=\beta C_\beta>1$, where $C_\beta$ is given by the $\nabla_0$ condition for $\Psi$ and Proposition~\ref{equi-nabla}. We have, for $h$ small enough and $x \geq \Psi^{-1}(1/h)$:
\begin{displaymath}
\frac{\Psi \big[\beta\Psi^{-1} (1/h) \big]}{\Psi \big[ \Psi^{-1} (1/h) \big]} 
\leq \frac{\Psi (\beta C_\beta\,x)}{\Psi (x)}= \frac{\Psi (A\,x)}{\Psi (x)}\,;
\end{displaymath}
we get hence, for $h$ small enough:
\begin{displaymath}
K_\mu(h) \leq  \frac{\Psi \big[\Psi^{-1}(1/h) \big]}{\Psi \big[\beta\Psi^{-1}(1/h)\big]} = 
\frac{1/h} {\Psi\big[\beta\Psi^{-1}(1/h)\big]}\,\raise 0,5mm\hbox{,}
\end{displaymath} 
and condition \eqref{eq:K compact} is fulfilled.
\par
With 1) and 2) previously shown, this finishes the proof.\hfill$\square$
\bigskip\goodbreak
\noindent{\bf Remark.} Actually, the proof of Theorem~\ref{continuite Carleson},~1) gives, for every 
measure $\mu$ on $\overline \D$ and every Orlicz function $\Psi$:
\begin{displaymath}
1 \geq \int_{\overline \D} \Psi \Big( \frac {|u_{\xi, 1- h}|} {\|u_{\xi, 1- h}\|_{L^\Psi (\mu)} }\Big)\,d\mu 
\geq \Psi \Big( \frac{1}{9\, \|u_{\xi, 1- h}\|_{L^\Psi (\mu)} }\Big)\,\mu \big( W (\xi, h) \big),
\end{displaymath}
and hence:
\begin{equation}\label{inegalite Carleson}
\mu \big( W (\xi, h) \big) \leq 
\frac{1} {\Psi \Big(\displaystyle  \frac{\! 1^{\phantom 1}}{9\, \|u_{\xi, 1- h}\|_{L^\Psi (\mu)} }\Big)} \cdot
\end{equation}
In particular, if $\mu =\mu_\phi$ is the image of the Haar measure $m$ under a self-map $\phi$ of $\D$, 
one has:
\begin{equation}\label{inegalite Carleson phi}
\mu \big( W (\xi, h) \big) 
\leq \frac{1} {\Psi \Big( \displaystyle \frac{\! 1^{\phantom 1}}{9\, \|C_\phi (u_{\xi, 1- h})\|_\Psi }\Big) } \cdot
\end{equation}
\medskip
Condition~\eqref{inegalite Carleson} allows to have an upper control for the $\mu$-measure of Carleson 
windows, with $\| u_{\xi, 1-h}\|_{L^\Psi (\mu)}$. It is possible, conversely, to majorize these norms.\par
\medskip
\begin{definition}  
We shall say that a measure $\mu$ on $\overline \D$ is a \emph{$\Psi$-Carleson measure} if there exists some $A>0$ such that
\begin{displaymath}
\quad \mu \big ( W (\xi, h)\big) \leq \frac{1}{\Psi \big(A \Psi^{-1} (1/h) \big)}\,\raise 0,5mm\hbox{,} 
\quad\hbox{for every }\xi \in \T\hbox{ and every }h \in (0,1).
\end{displaymath}
\medskip

We shall say that a measure $\mu$ on $\overline \D$ is a vanishing \emph{$\Psi$-Carleson measure} if, for every $A>0$,
\begin{displaymath}
\mathop{\lim}\limits_{h\to0}\Psi \big(A \Psi^{-1} (1/h)\big).\rho_\mu(h)=0.
\end{displaymath}
Equivalently if, for every $A>0$, there exists $h_A\in(0,1)$ such that
\begin{displaymath}
\quad \mu \big ( W (\xi, h)\big) \leq \frac{1}{\Psi \big(A \Psi^{-1} (1/h) \big)}\,\raise 0,5mm\hbox{,} 
\quad\hbox{for every }\xi \in \T\hbox{ and every }h \in (0,h_A).
\end{displaymath}
\end{definition}

We have the following characterizations:

\begin{proposition}\label{psicarleson}

$1)$ $\mu$ is a $\Psi$-Carleson measure on $\overline{\D}$ if and only if  there exists some constant $C \geq 1$ such that :
\begin{displaymath}
\| u_{\xi, 1-h}\|_{L^\Psi (\mu)} \leq \frac{C}{\Psi^{-1} (1/h)}\raise 0,5mm\hbox{,}  \hbox{ for every }\xi \in \T\hbox{ and every }h \in (0,1).
\end{displaymath}

$2)$ $\mu$ is a vanishing $\Psi$-Carleson measure on $\overline{\D}$ if and only if 
\begin{displaymath}
\mathop{\lim}\limits_{h\to0}\sup_{\xi\in\T}\Psi^{-1} (1/h)\| u_{\xi, 1-h}\|_{L^\Psi (\mu)}=0.
\end{displaymath}
\end{proposition}
\medskip

\noindent{\bf Proof.} The sufficiency (both for 1) and 2)) follows easily from \eqref{inegalite Carleson} in the preceding remark.

The converse is an obvious consequence of the following lemma.

\begin{lemme}\label{lempsicarleson}

Suppose that there exist $A>0$ and $h_0\in(0,1)$ such that: 
\begin{displaymath}
\rho_\mu(h)\le\frac{1}{\Psi\big(A\Psi^{-1}(1/h)\big)}\raise2pt\hbox{,}\hbox{ for every }h\in(0,h_0).
\end{displaymath}

Then there exists $h_1\in(0,1)$ such that: 
\begin{displaymath}
\| u_{\xi, 1-h}\|_{L^\Psi (\mu)}\leq\frac{24}{A\Psi^{-1}(1/h)}\raise2pt\hbox{,}
\end{displaymath}
for every $\xi\in\T$ and every $h\in(0,h_1)$.
\end{lemme}

\noindent{\bf Proof of the lemma.} This is inspired from \cite{Ga}, Chapter~VI, Lemma~3.3, page 239. 
We may assume that $ h \leq h_0/4\leq1/4$. 

First, writing $a = (1-h)\xi$ (where $\xi\in\T$), we observe that, when $|z-\xi|\ge bh$ for a $b>0$, we have:
\begin{align*}
| 1 - \bar{a} z |^2&=1+|a|^2|z|^2-2|a|\Re(\bar\xi z)\\
&=|a| |\xi-z|^2+(1-|a|)+|a|^2|z|^2-|a||z|^2\\
&\ge |a|b^2h^2+(1-|a|)^2\ge(|a|b^2+1)h^2.
\end{align*}
So, we have $| u_{\xi, 1-h} (z) |\le\displaystyle\frac{1}{|a|b^2+1}\le\min(1,2/b^2)$, when $|z-\xi|\ge bh$.

Now, define, for every $n\in\N$ and $\xi\in\T$:
\begin{displaymath}
S_n = S ( \xi, 2^{n+1} h) =\{ z \in \overline{\D}\,;\ |z - \xi | < 2^{n+1} h\}\subset W(\xi, 2.2^{n+1} h).
\end{displaymath}
Our observation implies that $| u_{\xi, 1-h} (z) |\leq\min(1,2/4^n)$, for every $z\in\overline{\D}\setminus S_{n-1}$. For $z \in S_0$, one has simply $|u_{\xi, 1-h} (z) | \leq 1$. 

There exists an integer $N$ such that $2^{N+2}h\le h_0<2^{N+3}h$. 

Let us compute:
\begin{align*}
\int_{\overline{\D}} \Psi \big( M\,|u_{\xi, 1-h}|\big)\,d\mu 
& = \int_{S_0} \Psi \big( M\,|u_{\xi, 1-h}|\big)\,d\mu 
+ \sum_{n=1}^{N} \int_{S_n \setminus S_{n-1}} \hskip -6mm \Psi \big( M\,|u_{\xi, 1-h}|\big)\,d\mu \\
&\hskip3,7cm +\int_{\overline\D\setminus S_N} \hskip -6mm \Psi \big( M\,|u_{\xi, 1-h}|\big)\,d\mu \\
& \leq \Psi (M) \mu (S_0) + \sum_{n=1}^N \Psi \Big( \frac{2M}{4^n} \Big) \mu (S_n) +\Psi \Big( \frac{2M}{4^N} \Big) \mu (\overline\D)\\
&\leq\sum_{n=0}^N \frac{1}{2^{n+1}}\Psi\Big( \frac{4M}{2^n}\Big) \mu (S_n)+\Psi \Big( \frac{2M}{4^N} \Big) \mu (\overline\D).
\end{align*}

But for $n\le N$, we have $2.2^{n+1} h\le 2^{N+2}h\le h_0$, so the hypothesis gives:
\begin{displaymath}
\mu(S_n)\le\frac{1}{\Psi\big(A\Psi^{-1}(1/2^{n+2}h)\big)}\cdot
\end{displaymath}

Take now:
\begin{equation}\label{M}
M = \frac{A}{24} \Psi^{-1} \Big(\frac{1}{h}\Big) \cdot
\end{equation}
We have, using that $\Psi^{-1}\big(\frac{1}{2^{n+2}h}\big)\ge\frac{1}{2^{n+2}}\Psi^{-1}\big(\frac{1}{h}\big)$, 
\begin{align*}
\int_{\overline{\D}} \Psi \big( M\,|u_{\xi, 1-h}|\big)\,d\mu 
&\leq \sum_{n=0}^N \frac{1}{2^{n+1}}+\mu (\overline\D)\Psi\bigg(\frac{A}{12.4^N}\Psi^{-1} \Big(\frac{1}{h}\Big)\bigg)\\
&\leq  \frac{1}{2}+\mu (\overline\D)\Psi\bigg(\frac{16A.h^2}{3h_0^2}\Psi^{-1} \Big(\frac{1}{h}\Big)\bigg),
\end{align*}
because $\displaystyle\frac{1}{4^N}\le\bigg(\frac{8h}{h_0}\bigg)^2$.

We can choose $h_1$ small enough to have:
\begin{displaymath}
\mu (\overline\D)\Psi\bigg(\frac{16A.h^2}{3h_0^2}\Psi^{-1} \Big(\frac{1}{h}\Big)\bigg)\le\frac{1}{2}
\end{displaymath}
for every $h\in(0,h_1)$, since $\mathop{\lim}\limits_{h\to0}h^2\Psi^{-1} \Big(\frac{1}{h}\Big)=0$.

We get for such $h$
\begin{displaymath}
\int_{\overline{\D}} \Psi \big( M\,|u_{\xi, 1-h}|\big)\,d\mu 
\leq1,
\end{displaymath}
so that
\begin{displaymath}
\| u_{\xi, 1-h} \|_{L^\Psi (\mu)} \leq \frac{1}{M} = \frac{24}{A} \frac{1} {\Psi^{-1} (1/h) } 
\, \raise 0,5mm \hbox{,}
\end{displaymath}
as it was announced. \hfill $\square$
\bigskip

\noindent{\bf Examples and counterexamples.} 
\smallskip

We are going to give some examples showing that we do not have the reverse implications in general in Theorem~\ref{continuite Carleson} and Theorem~\ref{compacite Carleson}.
\medskip

{\bf 1.} {\sl Condition \eqref{eq:rho} is not sufficient in general to have a continuous embedding.} Let $\Psi(x)={\rm e}^x-1$ (note that this Orlicz function even fulfills the $\Delta^1$ condition~!). Note that $\Psi(A\Psi^{-1}(1/h))\sim h^{-A}$, when $h\to0$.

{\bf a.} Let $\nu$ be a probability measure on $\T$, supported by a compact set $L$ of Lebesgue measure zero, such that $\nu(I)\le|I|^{1/2}$, for each $I$. We can associate to $\nu$ the measure on $\overline{\D}$ defined by $\tilde\nu(E)=\nu(E\cap\T)$. Then the identity map from $H^\Psi$ to $L^\Psi(\tilde\nu)$ is not even defined. Nevertheless the condition  \eqref{eq:rho} is clearly fulfilled with $A=1/2$.

{\bf b.} Now, we exhibit a similar example (less artificial) on the open disk. Let $\nu$ be as previously. By a standard argument: for every integer $n$, there exists a function $g_n$ in the unit ball of the disk algebra such that $|g_n|=1$ on $L$ and $\|g_n\|_\Psi\le4^{-n}$. As $L$ is compact, there exists some $r_n\in(1/2,1)$ such that $|g_n(r_nz)|\ge1/2$ for every $z\in L$. Now, define the measure $\mu$ by:
\begin{displaymath}
\mu(E)=\sum_{n=1}^\infty\frac{1}{2^n}\nu_n(E),
\end{displaymath}
where:
\begin{displaymath}
\nu_n(E)=\nu\big(\{z\in\T|\, r_nz\in E\}\big).
\end{displaymath}

If $W$ is a Carleson window of ``size'' $h$ then, for each $n\ge1$, we have:
\begin{displaymath}
\nu\big(\{z\in\T|\, r_nz\in W\}\big)\le\nu(W\cap\T)\le (2h)^{1/2}. 
\end{displaymath}
Hence, $\mu(W)\le(2h)^{1/2}$ and the condition \eqref{eq:rho} is  fulfilled.

Nevertheless, the identity from  $H^\Psi$ to $L^1(\mu)$ is not continuous: $\|g_n\|_\Psi\le4^{-n}$; but 
\begin{displaymath}
\|g_n\|_1\ge\frac{1}{2^n}\int_{r_n\T}|g_n|\,d\nu_n\ge\frac{1}{2^n}\int_L|g_n(r_nw)|\,d\nu(w)\ge\frac{1}{2^{n+1}}\cdot
\end{displaymath}

\medskip
{\bf 2.} {\sl Condition \eqref{eq:K} is not necessary in general to have a continuous embedding.} When $\Psi$ satisfies $\Delta_2$, the identity from $H^\Psi$ to $L^\Psi(\mu)$ is continuous if and only if $\mu$ is a Carleson measure. So the conditions \eqref{eq:rho}, \eqref{eq:K} and the continuity are equivalent in this case. Actually, when $\Psi$ does not satisfy $\Delta_2$, we construct below a measure $\mu$ on $\D$ such that the identity from $H^\Psi$ to $L^\Psi(\mu)$ is continuous and order bounded , but $\mu$ is not a Carleson measure ({\it a fortiori} does not verify \eqref{eq:K}). Note that the measure $\mu$ is then $\Psi$-Carleson but not Carleson. Here is the example:

We have assumed that $\Psi$ does not satisfy $\Delta_2$; so there exists an increasing sequence $(a_n)_{n\ge1}$ such that $\displaystyle\frac{\Psi(a_n)}{n}$ is increasing and $\displaystyle\frac{\Psi(2a_n)}{\Psi(a_n)}\ge n2^n$. Now, define the discrete measure
\begin{displaymath}
\mu=\sum_{n=1}^\infty\bigg(\frac{n}{\Psi(2a_n)}-\frac{n+1}{\Psi(2a_{n+1})}\bigg)\delta_{x_n},
\end{displaymath}
where:
\begin{displaymath}
x_n=1-\frac{1}{\Psi(2a_n)}\cdot
\end{displaymath}

As $\displaystyle\mu\big([x_N,1]\big)=\displaystyle\frac{N}{\Psi(2a_N)}\raise2pt\hbox{,}$ the measure $\mu$ is not Carleson: it should be bounded by $c(1-x_N)=\displaystyle\frac{c}{\Psi(2a_N)}\raise2pt\hbox{,}$ where $c$ is some constant.

We know that for every $f$ in the unit ball of $H^\Psi$ and every $x\in(0,1)$, we have $|f(x)|\le 4\Psi^{-1}\bigg(\displaystyle\frac{1}{1-x}\bigg)$: see Lemma~\ref{norme evaluation}. So we only have to see that $g\in L^\Psi(\mu)$, where $g(x)=\Psi^{-1}\bigg(\displaystyle\frac{1}{1-x}\bigg)$. Indeed, we have:

\begin{align*}
\int_\D\Psi\bigg(\frac{|g|}{2}\bigg)\,d\mu & =\sum_{n=1}^\infty\Big(\frac{n}{\Psi(2a_n)}-\frac{n+1}{\Psi(2a_{n+1})}\Big)\Psi\Big(\frac{|g(x_n)|}{2}\Big)\\
 &\le\sum_{n=1}^\infty\frac{n}{\Psi(2a_n)}\Psi\Big(\frac{1}{2}\Psi^{-1}\Big(\frac{1}{1-x_n}\Big)\Big)\\
 &\le\sum_{n=1}^\infty n\frac{\Psi(a_n)}{\Psi(2a_n)}\le\sum_{n=1}^\infty\frac{1}{2^n}\\
 &\le1
\end{align*}
so $\|g\|_\Psi\le2$.

\medskip
{\bf 3.} {\sl Condition \eqref{eq:K compact} is not necessary in general to have a compact embedding.} We can find, for every Orlicz function $\Psi$ not satisfying $\nabla_0$, a measure $\mu$ such that the identity from $H^\Psi$ to $L^\Psi(\mu)$ is compact but \eqref{eq:K compact} is not satisfied. Indeed: since $\Psi\notin\nabla_0$, we can select two increasing sequences $(x_n)_{n\ge1}$ and $(y_n)_{n\ge1}$, with $1\le x_n\le y_n\le x_{n+1}$ and $\Psi(x_n)>1$, and such that $\lim x_n=+\infty$ and:
\begin{displaymath}
\frac{\Psi(2x_n)}{\Psi(x_n)} \ge \frac{\Psi(2^ny_n)}{\Psi(y_n)}\cdot
\end{displaymath}

Define $\displaystyle r_n=1-\frac{1}{\Psi(y_n)}$ and the discrete measure:
\begin{displaymath}
\mu=\sum_{n=1}^\infty\frac{1}{\Psi(2^ny_n)}\,\delta_{r_n}.
\end{displaymath}
The series converge since $\Psi(2^ny_n)\ge2^n$. 

Thanks to Lemma~\ref{norme evaluation}, we have, for every $f$ in the unit ball of $H^\Psi$ and every $n\ge1$:
\begin{displaymath} 
|f(r_n)|\le 4\Psi^{-1}\bigg(\displaystyle\frac{1}{1-r_n}\bigg)=4y_n. 
\end{displaymath}
Given $r>r_1$, there exists an integer $N\ge1$ such that $r_N<r\le r_{N+1}$. Then, for every $f$ in the unit ball of $H^\Psi$,  we have $\displaystyle\|f\|_{L^\Psi(\overline{\D}\setminus r\D,\mu)}\le2^{-N+2}$, since:
\begin{align*}
\int_{\overline{\D}\setminus r\D}\Psi\bigg(\frac{|f|}{2^{-N+2}}\bigg)\,d\mu
&=\sum_{n=N+1}^\infty\frac{1}{\Psi(2^ny_n)}\Psi\bigg(\frac{|f(r_n)|}{2^{-N+2}}\bigg)\\
&\le\sum_{n=N+1}^\infty\frac{\Psi(2^Ny_n)}{\Psi(2^ny_n)}\\
&\le\sum_{n=N+1}^\infty\frac{1}{2^{n-N}}=1.
\end{align*}
This implies that $\displaystyle\mathop{\lim}\limits_{r\to1}\sup_{\|f\|_\Psi\le1}\|f\|_{L^\Psi(\overline{\D}\setminus r\D,\mu)}=0$. By Proposition~\ref{fullgene}, the identity from $H^\Psi$ to $L^\Psi(\mu)$ is compact.

On the other hand, writing $h_n=\displaystyle\frac{1}{\Psi(x_n)}$ and $t_n=\displaystyle\frac{1}{\Psi(y_n)}\raise2pt\hbox{,}$ we have:

\begin{align*}
K_\mu(h_n)\ge\frac{\mu([1-t_n,1])}{t_n}&=\Psi(y_n)\sum_{m=n}^\infty\frac{1}{\Psi(2^my_m)}\ge\frac{\Psi(y_n)}{\Psi(2^ny_n)}\\
&\ge\frac{\Psi(x_n)}{\Psi(2x_n)}=\frac{1/h_n}{\Psi\big(2\Psi^{-1}(1/h_n)\big)}\raise2pt\hbox{,}
\end{align*}
and this shows that \eqref{eq:K compact} is not satisfied.

\medskip
{\bf 4.} {\sl Condition \eqref{eq:rho compact} is not sufficient in general to have a compact embedding.} We can find an Orlicz function $\Psi$ and a vanishing $\Psi$-Carleson measure ({\it i.e.} \eqref{eq:K compact} is satisfied) $\mu$ such that the identity from $H^\Psi$ to $L^\Psi(\mu)$ is not compact.

We shall use the Orlicz function introduced in \cite{LLQR}. The key properties of this function  $\Psi$ are: 
\begin{itemize}

\item [i)] For every $x>0$, $\Psi(x)\ge x^3/3$.

\item [ii)] For every integer $k\ge1$, $\Psi(k!)\le (k!)^3$.

\item [iii)] For every integer $k\ge1$, $\Psi(3(k!))> k.(k!)^3$.
\end{itemize}

Once again, the job is done by a discrete measure. 

Define $x_k=k!$; $y_k=\displaystyle\frac{(k+1)!}{k^{1/3}}$; $r_k=1-\displaystyle\frac{1}{\Psi(y_k)}$ and $\rho_k=1-\displaystyle\frac{1}{\Psi(x_k)}\cdot$ Of course, $x_2<y_2<x_3<\cdots$.

Let $\nu$ be the discrete measure defined by: 
\begin{displaymath}
\nu=\sum_{k=2}^\infty\nu_k,
\end{displaymath} 
where:
\begin{displaymath} 
\nu_k=\frac{1}{\Psi\big((k+1)!\big)}\sum_{a^{k^2}=1}\delta_{r_ka}.
\end{displaymath} 
Observe that $\displaystyle\|\nu_k\|\le\frac{k^2}{\Psi\big((k+1)!\big)}\le\frac{3.k^2}{(k+1)!^3}$ so that the series converges. Note that $\nu$ is supported in the union of the circles of radii $r_k$ and not in a subset of the segment $[0,1]$ as in the preceding counterexamples.

In order to show that \eqref{eq:rho compact} is satisfied, it is clearly sufficient to prove that, when $\displaystyle\frac{1}{\Psi(y_k)}\le h<\displaystyle\frac{1}{\Psi(y_{k-1})}$ (with $k\ge3$), we have:
\begin{displaymath} 
\rho_\nu(h)\le\frac{1}{\Psi\Big(\displaystyle\frac{k^{1/3}}{2}\Psi^{-1}(1/h)\Big)}\cdot
\end{displaymath} 

Supposing then $\displaystyle\frac{1}{\Psi(y_k)}\le h<\displaystyle\frac{1}{\Psi(y_{k-1})}\raise2pt\hbox{,}$ we have $\Psi^{-1}(1/h)\le y_k$ so
\begin{displaymath} 
\Psi\Big(\frac{k^{1/3}}{2}\Psi^{-1}(1/h)\Big)\le\frac{1}{2}\Psi\big((k+1)!\big).
\end{displaymath} 
Therefore, it is sufficient to establish that $\displaystyle\rho_\nu(h)\le\frac{2}{\Psi\big((k+1)!\big)}\cdot$

A Carleson window $W(\xi,h)$ (where $\xi\in\T$) can contain at most one $k^2$-root of the unity, since 
$\displaystyle2h<\frac{2}{\Psi(y_{k-1})}\le\frac{6}{y_{k-1}^3}\le\frac{6(k-1)}{(k!)^3}\le\frac{2\pi}{k^2}\cdot$ This implies that
\begin{displaymath} 
\nu_k\big(W(\xi,h)\big)\le\frac{1}{\Psi\big((k+1)!\big)}\cdot
\end{displaymath} 
Nevertheless, when $j<k$, the window $W(\xi,h)$ cannot meet any circle of radius $r_j$ (centered at the origin), so $\nu_j\big(W(\xi,h)\big)=0$. We obtain:

\begin{align*}
\nu\big(W(\xi,h)\big)&=\sum_{j=k}^\infty\nu_j\big(W(\xi,h)\big)\le\frac{1}{\Psi\big((k+1)!\big)}+\sum_{j>k}\frac{j^2}{\Psi\big((j+1)!\big)}\\
&\le\frac{1}{\Psi\big((k+1)!\big)}+\sum_{j=k+1}^\infty\frac{3j^2}{(j+1)!^3}\\
&\le\frac{1}{\Psi\big((k+1)!\big)}+\frac{3}{(k+1)!^3}\sum_{s=1}^\infty\bigg(\frac{1}{k+1}\bigg)^s\\
&\le\frac{2}{\Psi\big((k+1)!\big)}\cdot
\end{align*}
This proves that \eqref{eq:rho compact} is satisfied.

Let us introduce the function $\displaystyle f_k(z)=x_k u_{1,\rho_k}\big(z^{k^2}\big)=x_k\bigg(\frac{1-\rho_k}{1-\rho_kz^{k^2}}\bigg)^2$. 

It lies in the unit ball of $H^\Psi$ by Corollary~\ref{norme $u_{a,r}$}:

\begin{displaymath} 
\|f_k\|_\Psi=x_k\|u_{1,\rho_k}\|_\Psi\le\frac{x_k}{\Psi^{-1}\big(\frac{1}{1-\rho_k}\big)}=1.
\end{displaymath} 

An easy computation gives $r_k^{k^2}\ge\rho_k$, for every $k\ge2$. So, for every $a\in\T$ with $a^{k^2}=1$, we have:
\begin{displaymath} 
f_k(ar_k)\ge x_k\bigg(\frac{1-\rho_k}{1-\rho_k^2}\bigg)^2\ge\frac{1}{4}x_k.
\end{displaymath} 
So 

\begin{align*}
\int_{\overline{\D}\setminus r_{k-1}\D} \Psi(12|f_k|)\,d\nu&\ge\int_{\overline{\D}\setminus r_{k-1}\D}\Psi(12|f_k|)\,d\nu_k\ge\frac{k^2}{\Psi((k+1)!)}\Psi(3x_k)\\
&>\frac{k^2}{\Psi((k+1)!)}\big(k.(k!)^3\big)\ge1.
\end{align*}
Therefore, we conclude that $\displaystyle\sup_{\|f\|_\Psi\le1}\|f\|_{L^\Psi(\overline{\D}\setminus r_k\D,\mu)}\ge\frac{1}{12}$, though $r_k\to1$. By Proposition~\ref{fullgene}, the identity from $H^\Psi$ to $L^\Psi(\mu)$ is not compact.
\bigskip\goodbreak

\subsection{Characterization of the compactness of composition operators}

For composition operators, compactness can be characterized in terms of $\Psi$-Carleson measures, as 
stated in the following result.

\begin{theoreme}\label{equiv Carleson}
For every analytic self-map $\phi\colon \D \to \D$ and every Orlicz function $\Psi$, the composition  
operator $C_\phi \colon H^\Psi \to H^\Psi$ is compact if and only if one has:
\begin{equation}
\rho_\mu (h) =o\,\Big( \frac{1}{\Psi\big(A\Psi^{-1}(1/h)\big)}\Big) \hskip 2mm \text{as } h\to 0,\hbox{ for every }A>0.
\tag{$R_0$}
\end{equation}
\end{theoreme}

In other words, if and only if $\mu_\phi$ is a \emph{vanishing $\Psi$-Carleson measure}.\par
\medskip

In order to get this result, we shall show that in Theorem~\ref{compacite Carleson}, 
conditions~\eqref{eq:rho compact} and \eqref{eq:K compact} are equivalent for the pull-back measure 
$\mu_\phi$ induced by $\phi$. This is the object of the following theorem.\par

\begin{theoreme}\label{equiv windows}
There exists a constant $k_1>0$ such that, for every analytic self-map $\phi \colon\D \to\D$, one has:  
\begin{equation}\label{eq equiv windows}
\mu_\phi \big( S (\xi,\eps h)\big) \leq k_1\,\eps\, \mu_\phi \big (S(\xi, h) \big)\,,
\end{equation}
for every $h\in (0,1-|\phi(0)|)$, and every $\eps\in (0,1)$.
\end{theoreme}

Note that we prefer here to work with the sets:
\begin{displaymath}
\qquad S(\xi,h) =\{z\in\overline\D \,;\ |z-\xi|< h\}\,,\quad \xi\in \T,\  0 < h <1,
\end{displaymath}
instead of the Carleson windows $W (\xi, h)$. Recall also that the pull-back measure $\mu_\phi$ is defined 
by~\eqref{pull-back}.\par
\medskip

We are going to postpone the proof  of Theorem~\ref{equiv windows}, and shall give before some consequences.

\subsubsection{Some consequences}

An immediate consequence of Proposition~\ref{psicarleson} and Theorem~\ref{equiv Carleson} is the following

\begin{theoreme}\label{equiv W-comp}
Let $\phi:\D\to\D$ be analytic and $\Psi$ be an Orlicz function. 

The operator $C_\phi$ on $H^\Psi$ is compact if and only if 
\begin{equation}
\hskip 1cm \sup_{\xi\in \T} \|C_\phi (u_{\xi, 1-h})\|_\Psi
=o\,\bigg(\frac{1}{\Psi^{-1}\big(1/h\big)}\bigg),
\hskip 3mm \text{as}\ h\to 0.\tag*{\ref{petit o} }
\end{equation}
\end{theoreme}

We deduce:

\begin{theoreme}\label{equiv weak-comp}
Let $\phi:\D\to\D$ be analytic. 

$1)$ Assume that the Orlicz function $\Psi$ satisfies condition $\Delta^0$. Then, the operator $C_\phi$ on $H^\Psi$ is weakly compact if and only if it is compact.

$2)$ Assume that the Orlicz function $\Psi$ satisfies condition $\nabla_2$. Then, the operator $C_\phi$ on $H^\Psi$ is a Dunford-Pettis operator if and only if it is compact.
\end{theoreme}

Recall that (see Theorem~\ref{corollaire Delta2}), under condition $\Delta^2$ for $\Psi$, 
the weak compactness of the composition operator $C_\phi$ is equivalent to its compactness, and even to 
$C_\phi$ being order bounded into $M^\Psi (\T)$. However, we shall see below, in 
Theorem~\ref{suite exemple}, that there exist Orlicz functions $\Psi \in \Delta^0$ (and even 
$\Psi \in \Delta^1 $) for which $C_\phi$ is compact, but not order bounded into $M^\Psi (\T)$.
\medskip

\noindent{\bf Proof.} In both cases, the result follows from Theorem~\ref{equiv W-comp}, since condition \ref{petit o} is satisfied. Indeed, if $C_\phi \colon H^\Psi \to H^\Psi$ is weakly compact and $\Psi \in \Delta^0$, we use Theorem~\ref{faiblement compact}. If $C_\phi \colon H^\Psi \to H^\Psi$ is a Dunford-Pettis operator, this is due to Proposition~\ref{DP}.\hfill$\square$
\medskip

Now, we have:

\begin{theoreme}\label{suite exemple}
There exist an Orlicz function $\Psi$ satisfying $\Delta^1$, 
and an analytic self-map $\phi \colon \D \to \D$ such that the composition operator 
$C_\phi\colon H^\Psi \to H^\Psi$ is not order bounded into $M^\Psi(\T)$, though it is compact.
\end{theoreme}

\noindent{\bf Remark.} It follows that our assumption that $\Psi \in \Delta^2$ in 
Theorem \ref{corollaire Delta2} is not only a technical one, though it might perhaps be weakened.
\bigskip

\noindent{\bf Proof.} Let:
\begin{displaymath}
\Psi(x) = \left\{ 
\begin{array}{lll}
\exp\big( (\log x)^2\big) & \text{if} & x\geq \sqrt{\e}\,,\\
\e^{-1/4} x & \text{if} & 0 \leq x \leq \sqrt{\e}.
\end{array}
\right.
\end{displaymath}
It is plain that $\Psi \in \Delta^1\cap \nabla_0$.\par
Moreover, for every $A>0$, one has, for $h$ small enough:
\begin{displaymath}
\frac{1/h}{\Psi\big(A\Psi^{-1}(1/h)\big)} =\exp\Big[- (\log A)^2 - 2(\log A)\sqrt{\log(1/h)}\Big].
\end{displaymath}
\par
Consider now $\phi=\phi_2$, the analytic self-map of $\D$ constructed in Theorem \ref{meme module}. 
Then $C_\phi\colon H^\Psi \to H^\Psi$ is not order bounded into $M^\Psi(\T)$, by Theorem 
\ref{impli-order bounded}, since, otherwise, $C_{\phi_1} \colon H^\Psi \to H^\Psi$ would  
also be order bounded into $M^\Psi(\T)$, which is easily seen to be not the case (we may also argue as follows: 
$C_{\phi_1} \colon H^\Psi \to H^\Psi$ would be compact, and hence, by Theorem \ref{implique Carleson}, 
$C_{\phi_1}$ would be compact from $H^2$ into $H^2$, which is false).\par
On the other hand, we have proved that $\rho_{\mu_\phi} (h)=O\,\big(h^{3/2}\big)$.
So the conclusion follows from Theorem \ref{compacite Carleson}, 3) 
and the fact that for every $c>0$:
\begin{displaymath}
-(\log A)^2 - 2(\log A)\sqrt{\log(1/h)} \geq -c\frac{1}{2}\log(1/h)
\end{displaymath} 
when $h$ is small enough.\hfill$\square$
\bigskip

\subsubsection{Preliminary results} 

We shall use the \emph{radial maximal function} $N$, defined for every harmonic function $u$ on $\D$ by:
\begin{equation}
\qquad (Nu) (\xi) = \sup_{0 \leq r <1} | u (r\xi) |\,,\quad \xi \in \T.
\end{equation}

Recall that for every positive harmonic function $u\colon \D \to \C$ whose boundary values $u^\ast$ are in 
$L^1(\T)$ ({\it i.e.} $u\in h^1$), one has, for every $\xi \in \T$:
\begin{equation}\label{ineg radiale}
Nu (\xi ) \leq Mu^\ast (\xi) \leq \pi Nu (\xi)\,,
\end{equation}
where $Mu^\ast$ is the Hardy-Littlewood maximal function of $u^\ast$ (see \cite{ABR}, Theorem~6.31, and 
\cite{Rudin}, Theorem 11.20 and Exercise~19, Chapter~11).
\medskip

We shall denote by $\Pi$ the right half-plane:
\begin{displaymath}
\Pi= \{ z \in \C \,;\ {\rm Re}\,z >0\}
\end{displaymath}
and by ${\mathscr C}$ the cone:
\begin{displaymath}
{\mathscr C} =\{z \in \C\,;\ -\pi/6 < {\rm Arg}\, z < \pi/6\}.
\end{displaymath}

The next result follows from Kolgomorov's Theorem, saying that  the Hilbert transform is a weak
$(1,1)$ operator (in applying this theorem to the positive harmonic function $2\,{\rm Re}\,g = g +\bar{g}$, 
noting that $\| {\rm Re}\,g \|_1 = {\rm Re}\,g (0) = g(0)$).

\begin{lemme}\label{Kolmo}
There exists a constant $c>0$ such that, if $g \colon \D \to \Pi$ is an analytic function with $g(0) >0$,  then:
\begin{displaymath}
\qquad m(\{ |g^*| > \lambda\}) \le c\, \frac{g(0)}{\lambda} \raise 0,5mm\hbox{,} 
\quad\hbox{for all $ \lambda > 0$.}
\end{displaymath}
\end{lemme}

Applying Lemma~\ref{Kolmo} to $g(z) =\bigl( f(z)+ \overline{f(0)}\bigr)^3$, and taking into account
that $|w_1+w_2|\ge |w_1|$, if  $w_1$, $w_2\in {\mathscr C}$, we get:

\begin{lemme}\label{cons Kolmo}
Let $f \colon \D\to {\mathscr C}$ be an analytic function with values in the cone ${\mathscr C}$, and write $f =u+iv$, 
with $u,v$ real-valued. Then:
\begin{displaymath}
\qquad m(\{ |f^*| > \lambda\}) \le 8c\, \Big(\frac{u(0)}{\lambda}\Big)^3\!,  \quad 
\hbox{for all $\lambda > 0$.}
\end{displaymath}
\end{lemme}

The next proposition is one of the keys. We postpone its proof.
\goodbreak

\begin{proposition}\label{prop-cle}
There exists a constant $k_2>0$ such that, for every analytic function 
$f= u+iv  \colon \D \to {\mathscr C}$ with values in the cone ${\mathscr C}$, one has:
\begin{equation}\label{ineg prop-cle}
\hskip 2mm m(\{ |f^*| > \lambda\} \cap I) \le k_2 \left(\frac{\alpha}{\lambda} \right)^3 m(I),\quad 
\hbox{for all $\lambda > 0$,}
\end{equation}
where $I$ is the arc $I=\{\e^{it} \,;\ a< t < b\}$, with $a,b\in \R$  and $\alpha>0$ satisfying 
$b-a< \pi/2$, $\alpha\ge Nu(\e^{ia})$, and $\alpha\ge Nu(\e^{ib})$.
\end{proposition}

As a corollary we obtain:

\begin{proposition}\label{cons prop-cle}
Let $f \colon \D \to {\mathscr C}$ be an analytic function, and write $f=u +iv$, as in 
Proposition~\ref{prop-cle}. If $\alpha>0$ satisfies $m(\{ Nu > \alpha\}) < 1/4$, then:
\begin{equation}\label{ineq cons prop-cle}
m(\{ |f^*| > \lambda\})\le k_2 \left(\frac{\alpha}{\lambda}\right)^3 m(\{ Nu > \alpha\}),
\quad \text{for all $\lambda \ge \frac{2\alpha}{\sqrt 3}\,\cdot$}
\end{equation}
\end{proposition}

\noindent{\bf Proof. } The set $\{ Nu > \alpha\}$ is open, and one can decompose it into a disjoint union of 
open arcs $\{ I_j\}_j$. Each arc has measure $m (I_j)\le m(\{ Nu >\alpha\}) < 1/4 $, and so is an arc of 
length less than $\pi/2$. We can then apply Proposition~\ref{prop-cle} and we obtain:
\begin{displaymath}
\qquad m(\{ |f^*| > \lambda\}\cap I_j)\le k_2 \left(\frac{\alpha}{\lambda} \right)^3 m(I_j),
\quad \hbox{for every $j$.}
\end{displaymath}
Summing up all these inequalities we get:
\begin{displaymath}
m\bigl(\{ |f^*| > \lambda\}\cap\{ Nu > \alpha\}\bigl) \le
k_2 \left(\frac{\alpha}{\lambda}\right)^3 m(\{ Nu > \alpha\})\,.
\end{displaymath}
The proposition follows since $|f^*|\le \frac{2}{\sqrt 3} u^* \le \frac{2}{\sqrt 3}\, Nu$, and
then $\{ |f^*| > \lambda\}$ is contained in $\{ Nu > \alpha\}$, for 
$ \lambda \ge \frac{2\alpha}{\sqrt 3}\cdot$ \hfill $\square$
\medskip

We shall need one more result.

\begin{proposition}\label{demi-alpha}
There exists a constant $k_3>0$ such that for every analytic function $f \colon \D \to {\mathscr C}$ with 
values in the cone ${\mathscr C}$, one has, writing $f =u + iv$: 
\begin{displaymath}
\qquad m(\{ Mu^* > \alpha\}\le k_3\, m(\{ u^* > \alpha/2\},
\quad \hbox{for every $\alpha > 0$. }
\end{displaymath}
\end{proposition}

In order to prove it, we shall first prove the following lemma.

\begin{lemme}\label{carre}
There exists a constant $k_4>0$ such that for every analytic function $f \colon \D \to {\mathscr C}$ with 
values in the cone ${\mathscr C}$, one has, writing $f =u + iv$: 
\begin{displaymath}
\qquad M\big((u^*)^2\big)(\xi)\le k_4\bigl(Mu^*(\xi)\bigr)^2\,,
\quad \hbox{for all $\xi\in \T$. }
\end{displaymath}
\end{lemme}

\noindent{\bf Proof. } Observe that $u^2 - v^2$ is a positive harmonic function since it is the real part of 
$f^2$. $f^2$ belongs to $H^1$ (see \cite{Du}, Theorem 3.2), so $u^2 - v^2 \in h^1$ and we can use 
inequalities~\eqref{ineg radiale}. We also have, since $-\pi/3 < {\rm Arg}\,(f^2) < \pi/3$, that 
$u^2 \ge 3 v^2$, and so $u^2 - v^2\ge 2u^2/3$ and $|f|\le \frac{2}{\sqrt 3}u$. We get:
\begin{align*}
\hskip 2cm M\big((u^*)^2\big)
&\le \frac{3}{2} M\big((u^*)^2 - (v^*)^2 \big) 
\le\frac{3\pi}{2} N(u^2 - v^2) \cr
&\le \frac{3\pi}{2} N(f^2) = \frac{3\pi}{2} (Nf)^2
\le \frac{3 \pi}{2} \left( \frac{2}{\sqrt 3} Nu\right)^2 \cr
&\le 2\pi \,(Mu^*)^2\,. \hskip 6,1cm \square
\end{align*}
\goodbreak

\noindent{\bf Proof of Proposition~\ref{demi-alpha}. } Write $A=\{ Mu^* > \alpha\}$ and 
$B=\{ u^* >  \alpha/2\}$. For every $\xi\in A$, there exists an open arc $I_\xi$ centered at $\xi$ such that:
\begin{displaymath}
p = \frac{1}{m(I_\xi)} \int_{I_\xi} u^* \,dm > Mu^*(\xi)/2\,,\quad\hbox 
{and}\quad  p>\alpha\,.
\end{displaymath}
We have, using Lemma~\ref{carre}:
\begin{displaymath}
\frac{1}{m(I_\xi)} \int_{I_\xi} (u^*)^2 \,dm \le M\big((u^*)^2\big) (\xi)
\le k_4  \bigl(Mu^*(\xi)\bigr)^2\le 4k_4\, p^2 .
\end{displaymath}
Let $L$ be the set $L =\{ u^*>p/2\}\cap I_\xi$. We have:
\begin{displaymath}
\frac{1}{m (I_\xi)} \int_{I_\xi \setminus L} \frac{u^\ast}{\! p}\,dm \leq \frac{1}{2}\,;
\end{displaymath}
hence, using the Cauchy-Schwarz inequality:
\begin{align*}
\frac{1}{2}
&\le \frac{1}{m (I_\xi)} \int_{I_\xi} \frac{u^\ast}{\! p}\,dm 
- \frac{1}{m (I_\xi)} \int_{I_\xi \setminus L} \frac{u^\ast}{\! p}\,dm 
= \frac{1}{m(I_\xi)}\int_L  \frac{u^*}{\! p}\,dm \\
& \le \sqrt{\frac{m(L)}{m(I_\xi)}} \left(\frac{1}{m(I_\xi)} 
\int_{I_\xi} \Bigl(\frac{u^*}{\! p}\Bigr)^2 \,dm  \right)^{1/2}
\le 2\,\sqrt{\frac{k_4 m(L)}{m(I_\xi)}}\,\cdot
\end{align*}
Therefore, $16k_4\,m(L)\ge m(I_\xi)$, and, since $L \subseteq B\cap I_\xi $, we have, for every $\xi\in A$, 
an arc $I_\xi$ containing $\xi$ such that $16k_4\,m( B\cap I_\xi)\ge m(I_\xi)$.
Applying the Hardy-Littlewood covering lemma, we then obtain:
\begin{displaymath}
\hskip 1,4cm m(A) \le 3\sum_{j=1}^n m(I_{\xi_j})
\le 3\times 16\times k_4 \times  \sum_{j=1}^n m(I_{\xi_j}\cap B) \le k_3 m(B).\hskip 1,3cm \square
\end{displaymath}
\medskip

\noindent{\bf Proof of Proposition~\ref{prop-cle}.} Composing $f$ with a suitable rotation, we can suppose 
that $a=-\delta $ and $b=\delta$, for $0<\delta<\pi/4$. Let us call $I^-$ and $I^+$ the arcs
\begin{displaymath}
I^- =\{\e^{it} \,;\  -\delta< t< 0\}\,, \qquad
I^+=\{\e^{it} \,;\  0< t< \delta\}\,.
\end{displaymath}
We shall prove that:
\begin{equation}\label{I-plus}
m(\{ |f^*| > \lambda\}\cap I^+) \le k_2 \left(\frac{\alpha}{\lambda} \right)^3 m(I^+)\,,
\end{equation}
using just the fact that $Nu(\e^{i\delta})\le \alpha$.\par 
In the same way one can prove:
\begin{equation}\label{I-moins}
m(\{ |f^*| > \lambda\}\cap I^-) \le k_2 \left(\frac{\alpha}{\lambda} \right)^3 m(I^-)\,,
\end{equation}
using just that $Nu(\e^{-i\delta})\le \alpha$.\par
Then, summing up \eqref{I-plus} and \eqref{I-moins}, Proposition~\ref{prop-cle} will follow.
\medskip

Let $Q$ be the right half-disc  
\begin{displaymath}
Q =\{ z\in \D \,;\ {\rm Re}\, z > 0\},
\end{displaymath}
and denote by $\psi$ the (unique) homeomorphism from $\overline\D$ to $\overline Q$, which is a conformal 
mapping from $\D$ onto $Q$ and sends $1$ to $1$, $i$ to $i$, and $-i$ to $-i$.\par
We can construct $\psi$ as the composition of a Moebius transformation $T$,  with the square root function 
and then with $T^{-1}$. Namely, let:
\begin{displaymath}
Tz = -i \frac{z+i}{z - i}\,;
\end{displaymath}
$T$ maps $\D$ onto the upper-half plane, sending $-i$ into $0$, $-1$ into $-1$, 
$1$ into $1$, and also $0$ into $i$, and $i$ into $\infty$.
The square root function maps the upper-half plane into the first quadrant and $T^{-1}$:
\begin{displaymath}
T^{-1}z = \frac{z - i}{1 - iz}
\end{displaymath}
maps this quadrant onto the half-disk  $Q$.\par
It is not difficult to see that $\psi(-1)=0$, $\psi(0)=\sqrt 2 -1$,  
and that there exist $\rho \in (0,\pi/2)$ such that
$\psi(\e^{i\rho}) = \e^{\pi i/4}$ and $\psi(\e^{-i\rho}) = \e^{-\pi i/4}$
(we must have $\e^{i\rho} = 1/3 + (\sqrt 8)i/3$; hence $\rho =\arctan (\sqrt 8)$).\par
\smallskip

Let $J$ be the arc: 
\begin{displaymath}
J=\{\e^{it} \,;\  -\rho < t < \rho\}.
\end{displaymath}
The map $\psi$ is regular on $J$, and so there exist two constants $\gamma_1$ and $\gamma_2>0$ such 
that for every Borel subset $E$ of $J$, one has:
\begin{displaymath}
\gamma_1\, m(E) \le m\bigl( \psi(E)\bigr) \le \gamma_2\, m(E) .
\end{displaymath}
If now $\beta\in (0,1)$, we put:
\begin{displaymath}
\psi_\beta(z) =(\psi(z))^\beta. 
\end{displaymath}
Then, it is easy to see that for every Borel subset $E$ of $J$, one has:
\begin{displaymath}
m\bigl( \psi_\beta(E)\bigr) =\beta \, m\bigl( \psi(E)\bigr),
\end{displaymath}
and so
\begin{displaymath}
\gamma_1\beta\, m(E) \le m\bigl( \psi_\beta(E)\bigr) \le \gamma_2\beta\, m (E).
\end{displaymath}
\medskip

In order to prove \eqref{I-plus}, consider the function $F \colon \D \to {\mathscr C}$ defined by:
\begin{displaymath}
\qquad F(z) = f\bigl(\e^{i\delta}\psi_\beta (z)\bigr)\,,\quad
\hbox{where $\beta =4\delta/\pi $.}
\end{displaymath}
Then: 
\begin{displaymath}
{\rm Re}\,\big(F(0)\big) =u\bigl((1-\sqrt 2)^\beta \e^{i\delta}\bigr)
\le Nu(\e^{i \delta}) \le\alpha. 
\end{displaymath}
Let us call $\chi$ the map:
\begin{displaymath}
\chi(z) = \e^{i\delta}\psi_\beta(z). 
\end{displaymath}
It is clear that $I^+$ is contained in $\chi (J)$. If $A=\{ |f^*| > \lambda\}\cap I^+$,
then $E=\chi^{-1} (A)$ is a Borel subset of $J$, and for every $\xi\in  E$, one has $|F^*(\xi)|>\lambda$. 
Then:
\begin{align*}
m(A) 
& = m\bigl( \chi (E)\bigr) =m\bigl( \psi_\beta (E)\bigr) 
\le \gamma_2 \beta\, m(E) \\
& \le 8\gamma_2 \frac{\delta}{2\pi} m(\{ |F^*| > \lambda\}) \cr
& = 8\gamma_2\, m(I^+)\,m(\{ |F^*| > \lambda\}) \cr
\noalign{\hbox{and using Lemma~\ref{cons Kolmo} for $F$,}}
& \le 8\gamma_2 \, m(I^+)  \times 8c_3 \left(\frac{{\rm Re}\, F(0)} {\lambda}\right)^3 \\
& \le 64\gamma_2 c_3  \left(\frac{\alpha}{\lambda}\right)^3 m(I^+)
=k_2 \left(\frac{\alpha}{\lambda}\right)^3 m(I^+)  \,.
\end{align*}

The proof of \eqref{I-plus} is finished, and Proposition~\ref{prop-cle} follows.\hfill $\square$

\bigskip\goodbreak

\subsubsection{Proof of Theorem~\ref{equiv windows}}

Using the fact $h\mapsto \mu_\phi\bigl(S(\xi,h)\bigl)$ is nondecreasing, it is enough to prove that
there exist $h_0> 0$, and $\varepsilon_0 > 0$, such that \eqref{eq equiv windows} is true for 
$0 < h < h_0 (1-|\phi(0)|)$, and $0 <\varepsilon <\varepsilon_0$, because changing the constant 
$k_1$, if necessary, the theorem will follow.\par
\smallskip

We can also suppose that $\xi=1$.\par
\smallskip

The real part of $1/\big(1-\phi(z) \big)$ is positive, in fact greater than 1/2, for every $z\in \D$. Take 
$0< h< h_0$, and consider the analytic function $f$ defined by:
\begin{displaymath}
f(z)=\left(\frac{h}{1-\phi(z)}\right)^{1/3} ,
\end{displaymath}
where the cubic root is taken in order that, for every $z\in \D$, $f(z)$ belongs to the cone
\begin{displaymath}
{\mathscr C} =\{z\in\C : -\pi/6 < {\rm Arg}\,(z) <\pi/6\,\}.
\end{displaymath}
Clearly $\mu_\phi\bigl( S(1,h) \bigl) = m(\{ |f^*|>1\})$ and
$\mu_\phi \bigl( S(1,\varepsilon h)\bigl) =m(\{|f^*|>1/\root 3 \of  \varepsilon \})$. We also have:
\begin{displaymath}
|f(0)| \le \left( \frac{h}{1 - |\phi(0)|}\right)^{1/3} < \root 3 \of {h_0}.
\end{displaymath}

We shall write $f=u+iv$ where $u$ and $v$ are real-valued harmonic functions. Observe that:
\begin{displaymath}
\qquad |v(z)|< \frac{1}{\sqrt 3}\, u(z), \quad \text{for every $z\in \D$.}
\end{displaymath}
It is known that $f\in H^p$, for every $p<3$ (see \cite{Du}, Theorem 3.2), and so $u$ and $v$ are the
Poisson integrals of $u^*$ and $v^*$ (in particular $u, v \in h^1$).\par
\medskip

We are looking for a control of
$m(\{ |f^*|>1/ \root 3 \of \varepsilon \})$ by $\varepsilon$ times $m(\{ |f^*|>1\})$. 
Proposition~\ref{cons prop-cle} provides this control replacing $m(\{ |f^*|>1\})$ by  
$m(\{ Nu > 2\})$:
\begin{displaymath}
m(\{ |f^*|>1/ \root 3 \of \varepsilon \})\le8k_2\eps m(\{ Nu > 2\}),
\end{displaymath}
when $m(\{ Nu > 2\})<1/4$.

As $f$ is valued in ${\mathscr C}$, $|f^*|$ is controlled by $u^*$. Then 
what we need in fact is a control of the measure of level sets of $Nu$ by the measure of level sets of $u^*$.
This will be done by using Proposition~\ref{demi-alpha}.\par
\smallskip

Indeed, by Proposition~\ref{demi-alpha}, we have:
\begin{displaymath}
m(\{Nu > 2\}) \le m(\{Mu^* > 2\}) \le k_4\, m (\{u^* > 1\}).
\end{displaymath}
We know that $||u^*||_1 = u(0) \le |f(0)| \le h_0^{1/3}$. Then, choosing $h_0$ small enough, we can have
$k_4 h_0^{1/3}< 1/4$, and so $m(\{Nu > 2\}) < 1/4$; therefore, we can use Proposition~\ref{cons prop-cle}. 
Moreover we also have $\{u^* > 1\} \subseteq \{|f^*| > 1\}$. Taking 
$ \varepsilon_0 < 3\sqrt 3 / 64$, we have, for $0<\varepsilon<\varepsilon_0$, by 
Proposition~\ref{cons prop-cle}:
\begin{align*}
m(\{ |f^*|> 1/\root 3 \of \varepsilon \}) 
& \le 8 k_2\, \varepsilon \, m(\{Nu > 2\})
\le 8 k_4k_2\, \varepsilon\,  m(\{u^* > 1\}) \\
& \le k_1\,\varepsilon\, m ( \{|f^*| > 1\})\,,
\end{align*}
taking $k_1= 8k_4k_2$. \hfill $\square$

\bigskip\goodbreak

\section{Bergman spaces}

\subsection{Bergman-Orlicz spaces}

\begin{definition}
Let $d{\mathscr A}(z) =\frac{dx\,dy}{\pi}$ ($z=x+iy$) be the normalized Lebesgue measure on $\D$. 
The \emph{Bergman-Orlicz space} ${\mathscr B}^\Psi$ denotes the space of analytic functions 
$f\colon \D \to \C$ which are in the Orlicz space $L^\Psi(\D, d{\mathscr A})$. The 
\emph{Bergman-Morse-Transue space} is the subspace 
${\mathscr B}M^\Psi = {\mathscr B}^\Psi\cap M^\Psi(\D,d{\mathscr A}).$
\end{definition}

${\mathscr B}^\Psi$, equipped with the induced norm of $L^\Psi(\D,d{\mathscr A})$, is a Banach 
space, as an obvious consequence of the following lemma, analogous to 
Lemma \ref{norme evaluation}. 

\begin{lemme}\label{norme evaluation Bergman}
For every $a\in \D$, the norm of the evaluation functional $\delta_a$, which maps 
$f\in {\mathscr B}^\Psi$ to $f(a)$, is:
\begin{displaymath}
\|\delta_a\| \approx \Psi^{-1}\Big(\frac{1}{(1-|a|)^2}\Big)\,\cdot
\end{displaymath}
\end{lemme}

\noindent{\bf Proof.} For every analytic function $g\colon \D \to \C$, the mean-value 
property gives:
\begin{displaymath}
g(0)=\int_\D g(z)\,d{\mathscr A}(z).
\end{displaymath}
Hence if $\phi_a\colon \D\to \D$ denotes the analytic automorphism 
\begin{displaymath}
\phi_a(z)= \frac{z -a}{1 -\bar{a}z}\,\raise 0,5mm\hbox{,}
\end{displaymath}
whose inverse is $\phi_a^{-1}=\phi_{-a}$, one has, for every 
$f\in {\mathscr B}^\Psi$, using the change of 
variable formula:
\begin{align*}
f(a) = f\circ \phi_{-a}(0) &=\int_\D f\circ \phi_{-a}(z)\,d{\mathscr A}(z) 
= \int_\D f(w)\,|\phi_a'(w)|^2\,d{\mathscr A}(w) \\
& = \int_\D f(w) H_a(w)\,d{\mathscr A}(w)\,,
\end{align*}
where:
\begin{displaymath}
H_a(w)= |\phi_a'(w)|^2 =\frac{(1-|a|^2)^2}{|1 - \bar{a}w|^4}\,\cdot
\end{displaymath}
The kernel $H_a$ plays for ${\mathscr B}^\Psi$ the role that the Poisson kernel $P_a$ plays for 
$H^\Psi$: the \emph{analytic reproducing kernel} $K_a$ for ${\mathscr B}^2$ being 
$K_a(z)= \frac{1}{(1-\bar{a}z)^2}\,\cdot$\par
We therefore have (using \cite{Rao}, Proposition 4, page 61, or \cite{Be-Sh}, Theorem 8.14):
\begin{displaymath}
|f(a)| \leq 2\|f\|_\Psi \|H_a\|_\Phi\,,
\end{displaymath}
which proves the continuity of $\delta_a$. To estimate its norm, we are going to majorize 
$\|H_a\|_\Phi$, with the help of Lemma \ref{norme Psi}. Let us notice that, on the one hand,  
$\|H_a\|_1=1$ (take $f=\ind$ in the above identity); and on the other hand:
\begin{displaymath}
\|H_a\|_\infty= \frac{(1 - |a|^2)^2}{(1 - |a|)^4} = \frac{(1 +|a|)^2}{(1 -|a|)^2} 
\,\cdot
\end{displaymath}
We get, setting $b= \|H_a\|_\infty$, and using Lemma \ref{norme Psi} 
for $\|\ \|_\Phi$: 
\begin{displaymath}
\| H_a\|_\Phi \leq \frac{b}{\Phi^{-1}(b)}\,\cdot
\end{displaymath}
But $b\leq \Phi^{-1}(b)\Psi^{-1}(b)$ (see \cite{Rao}, Proposition 1 (ii), page 14). Hence 
$\| H_a\|_\Phi \leq \Psi^{-1}(b) $. Now:
\begin{displaymath}
b\leq \frac{4}{(1 - |a|)^2}\,\cdot
\end{displaymath}
We have $\Psi^{-1}(4t) \leq 4\Psi^{-1}(t)$ for all $t>0$. 
It follows that 
\begin{displaymath}
\| H_a\|_\Phi \leq C\Psi^{-1}\Big(\frac{1}{(1-|a|)^2}\Big)\,\cdot
\end{displaymath}
Since $\|\delta_a\| \leq 2\|H_a\|_\Phi$, we get the upper bound in 
Lemma \ref{norme evaluation Bergman}.\par\smallskip

For the lower bound, we simply observe that $H_a=|G_a|$, where $G_a(z)=\displaystyle{(1-|a|^2)^2 \over (1-\overline{a} z)^4 }$, and by Lemma \ref{norme Psi}:
\begin{align*}
\|\delta_a \| & \geq \frac{|G_a(a)|}{\|G_a\|_\Psi} =\frac{|H_a(a)|}{\|H_a\|_\Psi} \geq  
\frac{\displaystyle\frac{1}{(1 -|a|^2)^2}}{b/\Psi^{-1}(b)} =  
\frac{\Psi^{-1}(b)}{b (1 -|a|^2)^2} \\
& = \frac{\Psi^{-1}(b)}{(1 + |a|)^4}\geq \frac{1}{16} \Psi^{-1}(b) 
\geq \frac{1}{16} \Psi^{-1}\Big(\frac{1}{(1-|a|)^2}\Big)\,,
\end{align*}
since $b \geq 1/(1 - |a|)^2$.\hfill$\square$

\begin{proposition}\label{densitepoly}

We have the following properties

\begin{itemize}

\item [i)] ${\mathscr B}M^\Psi$ is the closure of $H^\infty(\D)$ in $L^\Psi (\D, {\mathscr A})$ and actually the algebraic polynomials are dense in ${\mathscr B}M^\Psi$. 

\item [ii)] On the unit ball of ${\mathscr B}^\Psi$, the weak-star topology $\sigma\big(L^\Psi(\D,{\mathscr A}),M^\Phi(\D,{\mathscr A})\big)$ coincides with the topology of convergence on compact subsets of $\D$.

\item [iii)] ${\mathscr B}^\Psi$ is closed in $L^\Psi(\D,{\mathscr A})$ for the weak-star topology.

\item [iv)] If $\Psi\in\nabla_2$, ${\mathscr B}^\Psi$ is (isometric to) the bidual of ${\mathscr B}M^\Psi$.
\end{itemize}
\end{proposition}

\noindent{\bf Proof.} For the first point, let us fix $f\in{\mathscr B}M^\Psi$. 
Setting $f_r(z)= f(rz)$ for  $z\in \D$ and $0\leq r <1$, it suffices to show that 
$\| f_r - f\|_\Psi \mathop{\longrightarrow}\limits_{r\to 1} 0$, since, being analytic in the disk 
$r \overline{\D}\subset\D$,   $f_r$ can be uniformly approximated on $\overline{\D}$ 
by its Taylor series. But the norm of $M^\Psi$ is absolutely continuous (see \cite{Rao}, Theorem 14, 
page 84) and therefore, for every $\eps >0$, there is some $R>0$, with $1/3\leq R<1$, such that 
$\| f \ind_{\D \setminus R\D} \|_\Psi \leq \eps$; hence:
\begin{displaymath}
\int_{\D \setminus R\D} \Psi\Big(\frac{|f|}{4\eps}\Big)\,d{\mathscr A} \leq 
\int_{\D \setminus R\D} \frac{1}{4}\Psi\Big(\frac{|f|}{\eps}\Big)\,d{\mathscr A} \leq \frac{1}{4}\cdot
\end{displaymath}
When $r\geq \frac{2R}{R+1} \geq 1/2$, we therefore have:
\begin{displaymath}
\int_{\D \setminus \frac{1+R}{2}\D} \Psi\Big(\frac{|f_r|}{4\eps}\Big)\,d{\mathscr A} \leq 1,
\end{displaymath}
and, by convexity of $\Psi$:
\begin{align*}
\int_\D \Psi\Big(\frac{| f_r - f|}{8\eps}\Big)\,d{\mathscr A} 
& \leq  \int_{\frac{1+R}{2}\D}  \Psi\Big(\frac{| f_r - f|}{8\eps}\Big) \,d{\mathscr A} \\ 
& \hskip 18mm + \int_{\D \setminus \frac{1+R}{2}\D} \frac{1}{2} \Big[  \Psi\Big(\frac{|f_r|}{4\eps}\Big) 
+ \Psi\Big(\frac{|f|}{4\eps}\Big)\Big]\,d{\mathscr A} \\
& \leq  \int_{\frac{1+R}{2}\D}  \Psi\Big(\frac{| f_r - f|}{8\eps}\Big) \,d{\mathscr A} \\
& \hskip 10mm + \frac{1}{2} 
\int_{\D \setminus \frac{1+R}{2}\D} \Psi\Big(\frac{| f_r |}{4\eps}\Big)\,d{\mathscr A}  
+ \frac{1}{2} \int_{\D\setminus R\D }\Psi\Big(\frac{| f|}{4\eps}\Big)\,d{\mathscr A} \\
& \leq 1,
\end{align*}
for $r$ close enough to $1$ since $f_r - f$ tends to $0$ uniformly on $\frac{1 + R}{2}\D$. Hence, for some 
$r_0 <1$, one has $\| f_r - f \|_\Psi \leq 8\eps$ for $r_0\leq r <1$. This was the claim.

{\sl ii)} It suffices to use a sequential argument since the topologies are metrizable on balls (the space $M^\Phi(\D,{\mathscr A})$ is separable). Assume that $f\in{\mathscr B}^\Psi$ (with $\|f\|_\Psi\le1$) is the weak-star limit of a sequence of analytic functions $f_n\in{\mathscr B}^\Psi$ (with $\|f_n\|_\Psi\le1$). Testing this with the function $h_k(z)=(k+1)\overline{z}^k$, we obtain that the $k^{\rm th}$ Taylor coefficient $a_k(n)$ of $f_n$ converges to the $k^{\rm th}$ Taylor coefficient $a_k$ of $f$: 
\begin{displaymath} 
a_k(n)=\int_\D f_n  h_k  \, d{\mathscr A} \longrightarrow \int_\D f  h_k \, d{\mathscr A} =a_k \; ,\hbox{ for every }k\ge0.
\end{displaymath}

Fix a compact $K\subset\D$, there exists an $r\in(0,1)$ such that $K\subset r\D$. We have:
\begin{displaymath}
\sup_{z\in K}|f_n(z)-f(z)|\le\sum_{k\ge0}|a_k(n)-a_k|r^k\longrightarrow0
\end{displaymath}
by the dominated convergence Theorem (observe that $| a_k -a_k(n)|\le 2(k+1)$ for every $k\ge 0$ and every $n\ge0$).

For the converse, suppose now that $f_n\in{\mathscr B}^\Psi$ (with $\|f_n\|_\Psi\le1$) converges uniformly on every compact subsets of $\D$ to $f\in{\mathscr B}^\Psi$ (with $\|f\|_\Psi\le1$). Fixing $g\in M^\Phi(\D,{\mathscr A})$ and $\eps>0$, there exists an $r\in(0,1)$ such that $\|g\ind_{\D\setminus r\D}\|_\Phi\le\eps$. We have:

\begin{displaymath}
\bigg|\int_\D(f_n-f)gd{\mathscr A}\bigg|\le\bigg|\int_{r\D}(f_n-f)gd{\mathscr A}\bigg|+2\eps\le\sup_{z\in r\D}|f_n(z)-f(z)|.\|g\|_1+2\eps.
\end{displaymath}

By hypothesis $\displaystyle\sup_{z\in r\D}|f_n(z)-f(z)|\longrightarrow0$. The conclusion follows.

{\sl iii)} By the classical Theorem of Banach-Dieudonn\'e, it is sufficient to prove that the balls are weak-star closed (equivalently  weak-star compact) and by separability of $M^\Phi(\D,{\mathscr A})$, the weak-star topology is metrizable on balls. The previous fact ii) shows that it is equivalent to prove that the unit ball of ${\mathscr B}^\Psi$ is compact  for the topology of convergence on compact subsets. But this is easy: indeed, if $f_n$ in the unit ball of ${\mathscr B}^\Psi$. This is a normal family thanks to Lemma~\ref{norme evaluation Bergman}. A subsequence converges to an analytic function $f$ on compact subsets of $\D$ and the Fatou Lemma implies that $f$ actually lies in the unit ball. 

{\sl iv)} Assume now that $\Psi$ satisfies $\nabla_2$. Since $(M^\Psi)^{\ast\ast}=L^\Psi(\D,{\mathscr A})$, we have $({\mathscr B}M^\Psi)^{\ast\ast}=\overline{{\mathscr B}M^\Psi}^{w^\ast}$, in the space $L^\Psi(\D,{\mathscr A})$. Hence it suffices to show that $\overline{{\mathscr B}M^\Psi}^{w^\ast}={\mathscr B}^\Psi$. We already know that ${\mathscr B}^\Psi$ is weak-star closed. Now, let $f\in{\mathscr B}^\Psi$. Obviously, $f_r\in{\mathscr B}M^\Psi$ for every $r\in(0,1)$, where $f_r(z)=f(rz)$. Moreover $\|f_r\|_\Psi\le\|f\|_\Psi$. But this is clear that $f_r$ is uniformly convergent to $f$ on compact subsets of $\D$, when $r\rightarrow1^-$. By ii), the conclusion follows. \hfill$\square$

\medskip

In the previous proof, we can see the points ii) and iii) in a slighlty different way: the unit ball $B$ of ${\mathscr B}^\Psi$ is compact for the topology $\tau$ of uniform convergence on compact subsets. Then observe that the identity from $B$, equipped with the topology $\tau$, to $B$, equipped with the weak-star topology, is continuous (the sequential argument is sufficient by metrizability). This implies, since the weak-star topology is separated, that $B$ is weak-star compact (hence closed) and that the topologies coincide. Now, by Banach-Dieudonn\'e, the space ${\mathscr B}^\Psi$ is weak-star closed.
\bigskip

\subsection{Compact composition operators on Bergman-Orlicz spaces}

We shall begin with showing that, as in the Hardy-Orlicz case, every symbol $\phi$ defines a 
bounded composition operator.

\begin{proposition}
Every analytic self-map $\phi \colon \D \to \D$ induces a bounded composition operator 
$C_\phi \colon {\mathscr B}^\Psi \to {\mathscr B}^\Psi$. Moreover, 
$C_\phi f \in {\mathscr B}M^\Psi$ for every $f \in {\mathscr B}M^\Psi$; hence, when $\Psi\in\nabla_2$, the former 
operator is the bi-transposed of $C_\phi \colon {\mathscr B}M^\Psi \to {\mathscr B}M^\Psi$.
\end{proposition}

\noindent{\bf Proof.} It suffices to follow the lines of Proposition \ref{subordination}, and to 
integrate the integrals written there between $0$ and $1$, with respect to the measure $2rdr$.
\hfill$\square$
\bigskip

Before stating and proving the main theorem, we are going to prove the following auxiliary 
result, interesting in itself, and which reinforces an example of J. H. Shapiro (\cite{Shap}, 
Example, page 185).

\begin{proposition}\label{Blaschke}
There exists a Blaschke product $B$ having angular derivative at no point of $\T=\partial\D$, 
in the following sense:
\begin{equation}\label{Blaschke inegalite}
(\forall \eps >0)\ (\exists c_\eps>0)\hskip 5mm 
1 - |B(z)| \geq c_\eps (1 -|z|)^\eps\,,\hskip 2mm \forall z\in \D.
\end{equation}
\end{proposition}

\noindent{\bf Proof.} We shall take:
\begin{displaymath}
B(z)= \prod_{n=1}^{+\infty} \frac{|z_n|}{z_n}\frac{z_n -z}{1 -\bar{z}_n z}\,
\raise 0,5mm\hbox{,}
\end{displaymath}
where:
\begin{displaymath}
\{z_n\,;\ n\geq 1\} = \bigcup_{n\geq 1} A_n, 
\end{displaymath}
with:
\begin{displaymath}
A_n=\{r_n\omega_n^j \,;\ \omega_n=\e^{2\pi i/p_n}\,,\ 0\leq j \leq p_n-1\}\,,
\end{displaymath}
where $(r_n)_{n\geq 1}$ is a (strictly) increasing sequence with $0 < r_n <1$, and the integers $p_n$ 
will have to be adjusted, satisfying the Blaschke condition:
\begin{displaymath}
\sum_{n=1}^{+\infty} (1 -|z_n|) = \sum_{n=1}^{+\infty} p_n(1-r_n) <+\infty.
\end{displaymath}
\par

One has:
\begin{displaymath}
|B(z)|^2  = \prod_{n=1}^{+\infty} \Big| \frac{z_n -z}{1 -\bar{z}_n z}\Big|^2 
= \prod_{n=1}^{+\infty} \bigg[1 - 
\frac{(1 -|z|^2)(1 - |z_n|^2)}{|1 -\bar{z}_n z|^2}\bigg] 
\leq \exp\big(- S(z)\big)\,,
\end{displaymath}
where:
\begin{displaymath}
S(z)= \sum_{n=1}^{+\infty} \frac{(1 -|z|^2)(1 - |z_n|^2)}{|1 -\bar{z}_n z|^2}\,\cdot
\end{displaymath}

We now proceed to minorize $S(z)$. For this purpose, we shall need the following simple 
lemma, whose proof will be temporarily postponed.

\begin{lemme}\label{petit lemme}
For every positive integer $p$ and every $a\in\D$, one has, setting 
$\omega=\e^{2\pi i/p}$:
\begin{displaymath}
\frac{1}{p}\sum_{k=0}^{p-1} \frac{1}{|1 - a\omega^k|^2} = 
\frac{1 -|a|^{2p}}{1 -|a|^2}\frac{1}{|1- a^p|^2} \geq \frac{1}{4} p\,|a|^p.
\end{displaymath}
\end{lemme} 

Then, setting $r=|z|$, we have:
\begin{align*}
S(z) &\geq (1-r)\sum_{n=1}^{+\infty} (1-r_n) \sum_{k=0}^{p_n-1} 
\frac{1}{|1 - r_n \omega_n^{-k} z|^2} \\
& \geq \frac{1-r}{4}\, \sum_{n=1}^{+\infty} p_n^{\,2} (1-r_n) (r_n r)^{p_n}.
\end{align*}
We shall take:
\begin{displaymath}
p_n =\big[(1-r_n)^{\eps_n -1}\big]+1\,,
\end{displaymath}
with: 
\begin{displaymath}
\eps_n= \frac{1}{\sqrt n}\hskip 5mm \text{and} \hskip 5mm r_n=1 -\frac{1}{2^n}
\end{displaymath}
and where $[\ ]$ stands for the integer part. More explicitly:
\begin{displaymath}
p_n =\big[2^{n - \sqrt n}\big]+1.
\end{displaymath}
If $r\geq 1/2$, let $N\geq 1$ be such that $r_N < r \leq r_{N+1}$. One has:
\begin{displaymath}
S(z) \geq \frac{1-r}{4}\,(1-r_N)^{2\eps_N-1}\, r_N^{\,2p_N} 
\geq \frac{1}{8} (1-r_N)^{2\eps_N}\, r_N^{\,2p_N}\,,
\end{displaymath}
since $1-r \geq 1-r_{N+1} = (1-r_N)/2$. Moreover:
\begin{displaymath}
p_N \leq 2(1-r_N)^{\eps_N -1} = 2.2^{N (1-\eps_N)} \leq 2.2^N\,,
\end{displaymath}
so that:
\begin{displaymath}
S(z) \geq \frac{1}{8}(1-r)^{2\eps_N}\,(1- 2^{-N})^{2^{N+1}} \geq c(1-r)^{2\eps_N}\,,
\end{displaymath}
where $c$ is a positive numerical constant.\par
Hence, setting:
\begin{displaymath}
\eps(z)= 2\eps_N \hskip 5mm \text{for}\ |z|\geq 1/2\ \text{and}\ r_N < r \leq r_{N+1}\,,
\end{displaymath}
one has:
\begin{displaymath}
\eps(z)\mathop{\longrightarrow}_{|z|\mathop{\to}\limits^{<}1} 0\,,
\end{displaymath}
and we get:
\begin{displaymath}
1- |B(z)|^2 \geq 1 - \e^{-S(z)} \geq 1- \e^{-c(1-|z|)^{\eps(z)}} \geq c'(1-|z|)^{\eps(z)}\,,
\end{displaymath}
where $c'$ is another positive numerical constant. This gives condition 
(\ref{Blaschke inegalite}), since $1 -|B(z)| \geq \displaystyle\frac{1- |B(z)|^2}{2}\,\cdot$\par
\smallskip
Finally, the Blaschke condition is satisfied, since:
\begin{displaymath}
p_n(1-r_n) \leq 2(1-r_n)^{\eps_n}= 2. 2^{-\sqrt n}.
\end{displaymath}
This ends the proof of Proposition \ref{Blaschke}.\hfill$\square$
\par\medskip

\noindent{\bf Proof of Lemma \ref{petit lemme}.} Let $G_p$ be the finite group of $p^{th}$ roots 
of unity, equipped with its normalized Haar measure. For $u\colon G_p \to \C$ and 
$0\leq k \leq p-1$, we denote by $\hat u(k)$ the $k^{th}$ Fourier coefficient of $u$, 
{\it i.e.}:
\begin{displaymath}
\hat u(k)= \frac{1}{p} \sum_{z\in G_p} u(z)z^{-k}.
\end{displaymath}
Then, the Plancherel-Parseval formula for $G_p$ reads:
\begin{displaymath}
\sum_{k=0}^{p-1} |\hat u(k)|^2 =\frac{1}{p} \sum_{z\in G_p} |u(z)|^2.
\end{displaymath}
Applying this to 
\begin{displaymath}
u(z)=\frac{1}{1-az}=\sum_{\scriptscriptstyle \stackrel{l\geq 0} {\scriptscriptstyle 0\leq k\leq p-1}} 
a^{lp+k} z^k 
= \sum_{k=0}^{p-1} \hat u(k) z^k\,,
\end{displaymath}
with
\begin{displaymath}
\hat u(k)= \sum_{l\geq 0} a^{lp+k} =\frac{a^k}{1- a^p}\,\raise 0,5mm\hbox{,}
\end{displaymath}
we get:
\begin{displaymath}
\frac{1}{p}\sum_{k=0}^{p-1} \frac{1}{|1- a \omega^k|^2} = 
\sum_{k=0}^{p-1} \frac{|a|^{2k}}{|1- a^p|^2} = 
\frac{1- |a|^{2p}}{(1-|a|^2)|1- a^p|^2}\,\cdot
\end{displaymath}
To finish, we note that $|1- a^p| \leq 2$, and that, by the arithmetico-geometric inequality, 
we have, with $x=|a|^2$:
\begin{align*}
\frac{1 - |a|^{2p}}{1 - |a|^2} &= 1+ x + \cdots + x^{p-1} \\
&\geq p\big(x^{1+2+\cdots +(p-1)}\big)^{1/p} = p\, x^{\frac{p-1}{2}} \geq p\, x^{p/2} 
= p\, |a|^p.
\end{align*}
\vbox{}\hfill$\square$

\begin{theoreme}\label{Bergman-compact} 
If the composition operator $C_\phi \colon {\mathscr B}^\Psi \to {\mathscr B}^\Psi$ is compact, then  
\begin{equation}\label{condition Bergman-compact}
\frac{\Psi^{-1}\bigg[{\displaystyle \frac{1}{(1 - |\phi(a)|)^2}}\bigg]}
{\Psi^{-1}\bigg[{\displaystyle \frac{1}{(1 - |a|)^2}}\bigg]} 
\mathop{\longrightarrow}_{|a| \mathop{\to}\limits^< 1} 0.
\end{equation}
This condition is sufficient if $\Psi\in \Delta^2$.\par
\end{theoreme}

Before giving the proof of this theorem, let us note that in case where  
$\Psi=\Psi_2$, with $\Psi_2(x)=\e^{x^2}-1$, it reads: \emph{the composition operator 
$C_\phi \colon {\mathscr B}^{\Psi_2} \to {\mathscr B}^{\Psi_2}$ is compact if and only if:}
\begin{equation}\label{condition Bergman-compact Psi2}
(\forall \eps >0)\ (\exists c_\eps >0)\hskip 1cm 
1 - |\phi(z)| \geq c_\eps (1 -|z|)^\eps\,,\hskip 2mm \forall z\in \D.
\end{equation}
Indeed, $\Psi_2\in \Delta^2$, and:
\begin{displaymath}
\frac{\Psi_2^{-1}\bigg[{\displaystyle \frac{1}{(1 - |\phi(a)|)^2}}\bigg]}
{\Psi_2^{-1}\bigg[{\displaystyle \frac{1}{(1 - |a|)^2}}\bigg]} = 
\frac{\sqrt{\displaystyle \log \frac{1}{(1 - |\phi(a)|)^2}}} 
{\sqrt{\displaystyle \log \frac{1}{(1 - |a|)^2}}}\,;
\end{displaymath}
hence $C_\phi$ is compact if and only if 
\begin{displaymath}
\frac{\displaystyle \log \frac{1}{(1 - |\phi(a)|)^2}} 
{\displaystyle \log \frac{1}{(1 - |a|)^2}} \mathop{\longrightarrow}_{|a|\to 1} 0.
\end{displaymath}
Then, for every $\eps>0$, we can find some $c_\eps>0$ such that, for all $a\in \D$:
\begin{displaymath}
\log\frac{1}{1 - |\phi(a)|} \leq \eps\,\log \frac{1}{(1 - |a|)^2} + c_\eps\,;
\end{displaymath}
which is equivalent to (\ref{condition Bergman-compact Psi2}).\hfill$\square$
\medskip

That allows to have compact composition operators on ${\mathscr B}^{\Psi_2}$ which are not 
compact on $H^{\Psi_2}$. However it is very likely that this is the case for every 
$\Psi\in \Delta^2$, but we have not tried to see this full generality.

\begin{theoreme}
There exist symbols $\phi\colon \D \to \D$ such that the composition 
operators $C_\phi$ is compact from ${\mathscr B}^{\Psi_2}$ into itself, but not compact 
from $H^{\Psi_2}$ into itself, and even $C_\phi$ is  an isometry onto its image.
\end{theoreme}

A similar example is well-known for the Hilbert spaces ${\mathscr B}^2$ and $H^2$ (see \cite{Shap}, 
pages 180--186).
\medskip

\noindent{\bf Proof.} Let $B$ be a Blaschke product verifying the condition of Proposition 
\ref{Blaschke}. We introduce $\phi(z)=zB(z)$. The function still verifies $1 - |\phi(z)| \geq C_\eps (1 -|z|)^\eps$, since $|\phi(z)| \leq |B(z)|$ on $\D$. From Theorem \ref{Bergman-compact}, it follows, since $\Psi_2\in \Delta^2$, that $C_\phi\colon {\mathscr B}^{\Psi_2} \to {\mathscr B}^{\Psi_2}$ is compact.\par
We are now going to see that $C_\phi\colon H^{\Psi_2} \to H^{\Psi_2}$ is an isometry. Indeed, recall the following well-known fact, that we already used (see \cite{Nord}, Theorem 1): since $\phi$ is an inner function, the image $\phi(m)$ of the Haar measure $m$ of $\T$ under $\phi$ is equal to $P_a.m$, where $a=\phi(0)$ and $P_a$ is the Poisson kernel at $a$. Here $\phi(0) = 0$ so $\phi(m)=m$. It follows that for every $f\in H^{\Psi_2}$, one has, for $C>0$:
\begin{displaymath}
\int_\T \Psi_2 \Big( \frac {|f\circ \phi|}{C} \Big)\,dm =\int_\T \Psi_2 \Big( \frac{|f|}{C}\Big)\,dm 
\end{displaymath}

so that $\displaystyle\| f \|_{\Psi_2}=\|f \circ\phi\|_{\Psi_2}$.\hfill$\square$
\bigskip

We shall need also the following lemma, which completes Lemma \ref{norme evaluation Bergman}.

\begin{lemme}\label{valeur evaluation Bergman-Morse}
For every $f\in {\mathscr B}M^\Psi$, one has:
\begin{displaymath}
\hskip 1,5cm f(a) = o\,\bigg(\Psi^{-1}\Big(\frac{1}{(1 -|a|)^2}\Big)\bigg)\, \hskip 3mm 
\text{as}\hskip 3mm |a| \mathop{\longrightarrow}^{<}1.
\end{displaymath}
\end{lemme}

\noindent{\bf Proof.} This is obvious for the monomials $e_n\colon z\mapsto z^n$ since 
$|e_n(a)|\leq 1$, whereas 
$\Psi^{-1} \big(1 / (1 -|a|)^2\big) \mathop{\longrightarrow}\limits_{|a|\to 1} +\infty$. Since 
the evaluation $\delta_a$ is bounded on ${\mathscr B}M^\Psi$ and 
$\big\| \delta_a /\big(\Psi^{-1} (1/ 1 -|a|^2)\big) \big\| =\mathop{O\,} (1) $, it suffices to use that the 
polynomials are dense in ${\mathscr B}M^\Psi$; but this was already proved in Proposition~\ref{densitepoly}.\par\hfill$\square$
\bigskip

\noindent{\bf Proof of Theorem \ref{Bergman-compact}.} If 
$C_\phi \colon {\mathscr B}^\Psi \to {\mathscr B}^\Psi$ is compact, 
then so is the restriction $C_\phi \colon {\mathscr B}M^\Psi \to {\mathscr B}M^\Psi$ and its adjoint 
$C_\phi^\ast = C_\phi \colon \big({\mathscr B}M^\Psi\big)^\ast \to \big({\mathscr B}M^\Psi\big)^\ast$. 
Since $C_\phi^\ast(\delta_a) =\delta_{\phi(a)}$, Lemma \ref{valeur evaluation Bergman-Morse} gives 
$\delta_a / \|\delta_a\| \mathop{\longrightarrow}\limits_{|a| \to 1}^{w^\ast} 0$. Compactness of 
$C_\phi^\ast$ now leads us to 
\begin{displaymath}
C_\phi^\ast\Big( \frac{\delta_a }{\|\delta_a\|} \Big) \mathop{\longrightarrow}_{|a| \to 1}^{\|\ \|} 0.
\end{displaymath}
 That gives (\ref{condition Bergman-compact}), in view of Lemma \ref{norme evaluation Bergman}.\par
\smallskip

Conversely, assume that (\ref{condition Bergman-compact}) is verified. Observe first that, since 
$\Psi\in \Delta^2$, one has:
\begin{equation}\label{petite inegalite}
\hskip 1,5cm \Psi^{-1}(x^2) \leq \alpha \Psi^{-1} (x) \hskip 5mm \text{for $x$ large enough.}
\end{equation}
Indeed, let $x_0>0$ be such that $\Psi(\alpha x)\geq \big(\Psi(x)\big)^2$ for $x\geq x_0$. For 
$x\geq y_0=\sqrt {\Psi (\alpha x_0)}$, with $y=\Psi^{-1}(x^2)$, one has  
$x^2 =\Psi (y)\geq  \big(\Psi (y/\alpha)\big)^2$, 
{\it i.e.} $x\geq \Psi (y/\alpha)$, and hence $\Psi^{-1} (x^2)= y \leq \alpha \Psi^{-1}(x)$.\par
Therefore condition (\ref{condition Bergman-compact}), which reads: 
\begin{displaymath}
\Psi^{-1} \Big(\frac{1}{(1 - |\phi (z)|)^2}\Big) = o\,\bigg( \Psi^{-1} \Big(\frac{1}{(1 - |z|)^2}\Big) \bigg)\,,
\end{displaymath}
reads as well, because of (\ref{petite inegalite}):
\begin{equation}\label{nouvelle petite inegalite}
\Psi^{-1}\Big(\frac{1}{1 - |\phi (z)|}\Big) = o\,\bigg( \Psi^{-1}\Big(\frac{1}{1 - |z|}\Big) \bigg)\,,\hskip 5mm 
\text{as } |z|\to 1.
\end{equation}
We have to prove that (\ref{nouvelle petite inegalite}) implies the compactness of 
$C_\phi \colon {\mathscr B}^\Psi \to {\mathscr B}^\Psi$. So, by 
Proposition \ref{critere compacite general}, we have to prove that: for every sequence 
$(f_n)_n$ in the unit ball of ${\mathscr B}^\Psi$ which converges uniformly on compact sets of $\D$, 
one has $\| f_n \circ \phi\|_\Psi \mathop{\longrightarrow}\limits_{n\to \infty} 0$.\par
But (\ref{condition Bergman-compact}) and (\ref{petite inegalite}) imply that, for some $C>0$:
\begin{equation}
|f_n (z)| \leq C\, \Psi^{-1}\Big(\frac{1}{1 - |z|}\Big)\,\raise 0,5mm \hbox{,}\hskip 5mm \forall z\in \D\,.
\end{equation}
Let $\eps >0$ and set $\eps_0 = \eps / \alpha C$. Due to (\ref{nouvelle petite inegalite}), we can find some 
$r$ with $0 < r <1$ such that:
\begin{equation}\label{2*}
\left\{
\begin{array}{l}
 \displaystyle  \quad \sqrt{1-r} \leq \frac{1}{8}\,; \hskip 1cm 
\Psi^{-1} \Big(\frac{1}{1-r}\Big) \geq \alpha x_0\,; \\
 \\
 \Psi^{-1} \Big(\displaystyle \frac{1}{1 - |\phi(z)|}\Big) 
\leq \eps_0 \Psi^{-1} \Big(\displaystyle \frac{1}{ 1 - |z|}\Big)\hskip 2mm 
 \text{if}\ z\in \D \setminus r\D\,.
\end{array}
\right.
\end{equation}
\par
Then, since $(f_n)_n$ converges uniformly on $r\D$, we have, for $n$ large enough ($n\geq n_0$):
\begin{displaymath}
\int_{r\D} \Psi\Big(\frac{| f_n \circ \phi|}{\eps}\Big)\,d{\mathscr A}(z) \leq \frac{1}{2}\,\cdot
\end{displaymath}
On the other hand, by (\ref{2*}):
\begin{align*}
\int_{\D \setminus r\D} \Psi\Big(\frac{| f_n \circ \phi|}{\eps}\Big)\,d{\mathscr A}(z) 
& \leq \int_{\D \setminus r\D} \
\Psi\bigg(\frac{C}{\eps}\Psi^{-1}\Big(\frac{1}{1 - | \phi(z) |}\Big)\bigg)\,d{\mathscr A}(z) \\
& \leq  \int_{\D \setminus r\D} 
\Psi\bigg(\frac{\eps_0 C}{\eps}\Psi^{-1}\Big(\frac{1}{1 - | z |}\Big)\bigg)\,d{\mathscr A}(z) \\
& = \int_{\D \setminus r\D} \Psi\bigg(\frac{1}{\alpha} \Psi^{-1} \Big(\frac{1}{1 - |z|}\Big)\bigg)
\,d{\mathscr A}(z) \\ 
& \leq \int_{\D \setminus r\D} \sqrt{\Psi\bigg(\Psi^{-1}\Big(\frac{1}{1 - |z|}\Big)\bigg)}
\,d{\mathscr A}(z) \\
& \hskip 20mm \text{since}\ \Psi^{-1}\Big(\frac{1}{1-r}\Big) \geq\alpha x_0\,,\\
& = \int_{\D \setminus r\D} \frac{1}{\sqrt{1 - |z|}}\,d{\mathscr A}(z) 
= 2\int_r^1 \frac{\rho\,d\rho}{\sqrt{1 -\rho}} \\
& \leq 2\int_r^1 \frac{\,d\rho}{\sqrt{1 -\rho}} 
= 4 \sqrt{1-r} \leq \frac{1}{2}\,\cdot
\end{align*}
Putting together these two inequalities, we get, for $n\geq n_0$:
\begin{displaymath}
\int_\D \Psi\Big(\frac{| f_n \circ \phi|}{\eps}\Big)\,d{\mathscr A}(z) \leq \frac{1}{2} + \frac{1}{2} =1,
\end{displaymath}
and hence: $\| f_n \circ \phi\|_\Psi \leq \eps$, which ends the proof of Theorem \ref{Bergman-compact}.
\hfill$\square$

\bigskip\bigskip
\noindent{\bf Acknowledgement.} Part of this work was made when the fourth-named author 
was \emph{Professeur invit\'e de l'Universit\'e d'Artois} in May-June 2005 and during 
a visit of this author in Lille and in Lens in March 2006. He wishes to thank both departements for their hospitality. A part was made too during a stay of the first named author in the University of Sevilla in September 2006. He is gratefull to the {\it departamento de An\'alisis Matem\'atico}.

\bigskip

\vbox{\noindent{\it 
P. Lef\`evre and D. Li, Universit\'e d'Artois,
Laboratoire de Math\'ematiques de Lens EA 2462, 
F\'ed\'eration CNRS Nord-Pas-de-Calais FR 2956, 
Facult\'e des Sciences Jean Perrin,
Rue Jean Souvraz, S.P.\kern 1mm 18,\\ 
62\kern 1mm 307 LENS Cedex,
FRANCE \\ 
pascal.lefevre@euler.univ-artois.fr \hskip 3mm  \\ 
daniel.li@euler.univ-artois.fr
\smallskip

\noindent
H. Queff\'elec,
Universit\'e des Sciences et Techniques de Lille, 
Labo\-ratoire Paul Painlev\'e U.M.R. CNRS 8524, 
U.F.R. de Math\'ematiques,\par\noindent
59\kern 1mm 655 VILLENEUVE D'ASCQ Cedex, 
FRANCE \\ 
queff@math.univ-lille1.fr
\smallskip

\noindent
Luis Rodr{\'\i}guez-Piazza, Universidad de Sevilla, Facultad de 
Matematicas, Dpto de An\'alisis Matem\'atico, Apartado de Correos 1160,\par\noindent 
41\kern 1mm 080 SEVILLA, SPAIN \\ 
piazza@us.es\par}
}

\end{document}